\documentclass[]{nmeauth}

\usepackage{subfigure}
\usepackage{endfloat}

\newcommand{\real}{\mathbb{R}} 
\newcommand{\integer}{\mathbb{N}}
\newcommand{\vectornorm}[1]{\left\|#1\right\|}

\usepackage[numbers,sort&compress]{natbib}
\usepackage[lined,boxed]{algorithm2e}

\begin{document}
\NME{0}{0}{00}{00}{00}

\runningheads{M. Arnst, R. Ghanem, E. Phipps and J. Red-Horse}{Stochastic modeling of coupled problems}

\received{4 January 2010}
\norevised{}
\noaccepted{}

\title{Dimension reduction in stochastic modeling of coupled problems}

\author{M. Arnst\affil{1}$^{,}$\affil{2}\corrauth, R. Ghanem\affil{1}, E. Phipps\affil{3} and J. Red-Horse\affil{3}}

\address{\affilnum{1} 210 KAP Hall, University of Southern California, Los Angeles, CA 90089, USA.\\
                 \affilnum{2} B52/3, Universit\'{e} de Li\`{e}ge, Chemin des Chevreuils 1, B-4000 Li\`{e}ge, Belgium.\\
                 \affilnum{3} Sandia National Laboratories\footnotemark[2], P.O. Box 5800, Albuquerque, NM 87185, USA.}

\corraddr{210 KAP Hall, University of Southern California, Los Angeles, CA 90089, USA.}
\footnotetext[2]{Sandia National Laboratories is a multi-program laboratory managed and operated by Sandia Corporation, a wholly owned subsidiary of Lockheed Martin Corporation, for the U.S. Department of Energy's National Nuclear Security Administration under contract DE-AC04-94AL85000.}

\begin{abstract}
Coupled problems with various combinations of multiple physics, scales, and domains are found in numerous areas of science and engineering.
A key challenge in the formulation and implementation of corresponding coupled numerical models is to facilitate the communication of information across physics, scale, and domain interfaces, as well as between the iterations of solvers used for response computations.
In a probabilistic context, any information that is to be communicated between subproblems or iterations should be characterized by an appropriate probabilistic representation. 
Although the number of sources of uncertainty can be expected to be large in most coupled problems, our contention is that exchanged probabilistic information often resides in a considerably lower dimensional space than the sources themselves.  
This work thus presents an investigation into the characterization of the exchanged information by a reduced-dimensional representation and, in particular, by an adaptation of the Karhunen-Lo\`{e}ve decomposition.
The effectiveness of the proposed dimension-reduction methodology is analyzed and demonstrated through a multiphysics problem relevant to nuclear engineering.
\end{abstract}

\keywords{uncertainty quantification, coupled problems, multiphysics, polynomial chaos}

\section{Introduction}
The modeling and simulation of coupled systems governed by multiple physical processes that may exist simultaneously across multiple scales and domains are critical tools for addressing numerous challenges encountered in many areas of science and engineering. 
However, models are, by definition, only approximations of their target scenarios and are thus prone to modeling errors.
Additionally, parametric uncertainties may exist owing to various limitations in manufacturing and experimental methods. 
Uncertainty quantification~(UQ) thus constitutes a key requirement for achieving realistic predictive simulations.

Probability theory provides a rigorous mathematical framework for UQ, which permits a unified treatment of modeling errors and parametric uncertainties.
The first step in a probabilistic UQ analysis typically involves using methods from mathematical statistics~\citep{cramer1946,kullback1968} to characterize the uncertain features associated with a model as one or more random variables, random fields, random matrices, or random operators.  
The second step is to map this probabilistic representation of inputs through the system model into a probabilistic representation of responses.
This can be achieved in several ways, which include Monte Carlo sampling techniques~\citep{robert2005} and stochastic expansion methods.
The latter typically involve the computation of a representation of the predictions as a polynomial chaos~(PC) expansion.
Several approaches are available to calculate the coefficients in this expansion, such as embedded projection~\citep{ghanem2003,soize2004}, nonintrusive projection~\citep{soize2004}, and collocation~\citep{ghanem1998,ghanem1998b,ghiocel2002,xiu2005,babuska2007}.

A key challenge in the formulation and implementation of a coupled model is to facilitate the communication of information across physics, scale, and domain interfaces, as well as between the iterations of solvers used for response computations.
This information can comprise physical properties, energetic quantities, or solution patches, among other quantities.
Although the number of sources of uncertainty can be expected to be large in most coupled problems, we believe that the exchanged information often resides in a considerably lower dimensional space than the sources themselves.
Exchanged information can be expected to have a \textit{low effective stochastic dimension} in multiphysics problems when this information consists of a solution field that has been smoothed by a forward operator, and in multiscale problems when this information is obtained by summarizing fine-scale quantities into a coarse-scale representation.

In this work, we thus investigate the effectiveness of dimension-reduction techniques for the representation of the exchanged information.
We propose to represent the exchanged information by an adaptation of the Karhunen-Lo\`{e}ve~(KL) decomposition as this information passes from subproblem to subproblem and from iteration to iteration.
When the exchanged information has a low effective stochastic dimension, this representation allows a reduction in the number of requisite stochastic degrees of freedom to be achieved while maintaining accuracy, thus paving the way for a solution in a reduced-dimensional space, which in turn reduces the computational cost.
It should be noted that in references~\citep{doostan2007,maute2009,nouy2007,acharjee2006,guedri2006,sachdeva2006}, the integration of dimension-reduction techniques in algorithms for solving stochastic partial differential equations has already been demonstrated; however, this paper contributes by highlighting the role that dimension-reduction techniques can play in solving coupled problems.

The organization of this paper is as follows.
First, in Sec.~\ref{sec:sec4}, we outline the proposed methodology.
Then, in Secs.~\ref{sec:kl} and~\ref{sec:klb}, we describe a version of the KL decomposition that is well-adapted to construct a reduced-dimensional representation of exchanged information.
In Sec.~\ref{sec:sec5}, we focus on the implementation of our dimension-reduction methodology.
Finally, in Secs.~\ref{sec:sec6a} and~\ref{sec:sec6b}, we demonstrate the effectiveness of the proposed dimension-reduction methodology through an illustration problem for which we also provide numerical results.

\section{Dimension-reduction methodology}\label{sec:sec4}

\subsection{Model problem}
This paper is devoted to the solution of a stochastic coupled model of the following form:
\begin{equation}
\label{eq:coupling12}\begin{aligned}
\boldsymbol{f}(\boldsymbol{u},\boldsymbol{x},\boldsymbol{\xi})=\boldsymbol{0},&&\boldsymbol{y}=\boldsymbol{h}(\boldsymbol{u},\boldsymbol{\xi}),&&\boldsymbol{f}:\real^{r}\times\real^{s_{0}}\times\real^{m}\rightarrow\real^{r},&&\boldsymbol{h}:\real^{r}\times\real^{m}\rightarrow\real^{r_{0}},\\
\boldsymbol{g}(\boldsymbol{y},\boldsymbol{v},\boldsymbol{\zeta})=\boldsymbol{0},&&\boldsymbol{x}=\boldsymbol{k}(\boldsymbol{v},\boldsymbol{\zeta}),&&\boldsymbol{g}:\real^{r_{0}}\times\real^{s}\times\real^{n}\rightarrow\real^{s},&&\boldsymbol{k}:\real^{s}\times\real^{n}\rightarrow\real^{s_{0}}.
\end{aligned}
\end{equation}
To avoid certain technicalities involved in infinite-dimensional representations, we assume that these equations are discretized representations of a stochastic model that couples two physics, two scales, two domains, or a combination of these subproblems.
For instance, these equations may be obtained from the spatial discretization of a steady-state problem, or they may be the equations obtained at a single time step after the spatial and temporal discretization of an evolution problem.
Further, we assume that the data of the first subproblem, which enter this subproblem as coefficients or loadings or both, depend on a finite number of uncertain real parameters denoted as $\xi_{1},\ldots,\xi_{m}$, and that the data of the second subproblem depend on a finite number of uncertain real parameters denoted as~$\zeta_{1},\ldots,\zeta_{n}$.
Lastly, we model these sources of uncertainty as random variables and collect them into vectors, $\boldsymbol{\xi}=(\xi_{1},\ldots,\xi_{m})$ and~$\boldsymbol{\zeta}=(\zeta_{1},\ldots,\zeta_{n})$, which are assumed to be defined on a probability triple $(\Theta,\mathcal{T},P)$ and considered to have values in $\real^{m}$ and $\real^{n}$, respectively.
A probability triple is characterized by its constituents: $\Theta$, a sample space of possible outcomes; $\mathcal{T}$, a collection of subsets of $\Theta$ known as events; and $P$, a probability measure.

The stochastic coupled model~(\ref{eq:coupling12}) is a general bidirectionally coupled model. 
The \textit{solution variables}, $\boldsymbol{u}$, of the first subproblem, $\boldsymbol{f}$, depend on the solution variables, $\boldsymbol{v}$, of the second subproblem, $\boldsymbol{g}$, through \textit{coupling variables}, $\boldsymbol{x}$; likewise, $\boldsymbol{v}$ depends on~$\boldsymbol{u}$ through~$\boldsymbol{y}$.

Thus, to solve this stochastic coupled model, we require to find the random variables~$\boldsymbol{u}$ and~$\boldsymbol{v}$ defined on~$(\Theta,\mathcal{T},P)$ with values in~$\real^{r}$ and~$\real^{s}$ such that~(\ref{eq:coupling12}) is satisfied
under the assumption that the stochastic coupled model is well-posed in that it admits a unique and stable solution.

For example, if the stochastic coupled model were a fluid-structure interaction model, Eq.~(\ref{eq:coupling12}) could represent the fluid and the structural model; the solution variables $\boldsymbol{u}$ and~$\boldsymbol{v}$ could collect the solution fields required to describe the states of the fluid and the structure; and the coupling variables $\boldsymbol{x}$ and $\boldsymbol{y}$ could be the traces of the velocity field of the structure and the pressure field of the fluid on the deforming fluid-structure interface.

\subsection{Partitioned iterative solution}
Because a coupled model usually characterizes its response only in an implicit manner, the numerical solution of a coupled model typically requires an iterative method. 
One strategy could be to develop an operator-specific iterative method and associated solver for the coupled model as a whole.
However, here, we assume that iterative methods and solvers already exist for each subproblem, and we therefore consider a hybrid iterative method that reuses the aforementioned separate solvers as steps in a global iterative method built around them to obtain a solution to the coupled model.
The former approach is commonly labeled \textit{monolithic} in the literature, whereas the latter is referred to as \textit{partitioned}.
Clearly, the key advantage of the partitioned method is that it enables immediate reuse of legacy software that may already be available to solve subproblems within their own physical and geometric domains, where one can account for suboperator specifics such as mesh details and problem scales.  

Let us assume that each of the aforementioned separate iterative methods is based on the reformulation of the associated subproblem as a fixed-point problem:
\begin{equation}
\label{eq:couplingA12}\begin{aligned}
&\boldsymbol{u}=\boldsymbol{a}(\boldsymbol{u},\boldsymbol{x},\boldsymbol{\xi}),&&\boldsymbol{y}=\boldsymbol{h}(\boldsymbol{u},\boldsymbol{\xi}),&&\boldsymbol{a}:\real^{r}\times\real^{s_{0}}\times\real^{m}\rightarrow\real^{r},&&\boldsymbol{h}:\real^{r}\times\real^{m}\rightarrow\real^{r_{0}},\\
&\boldsymbol{v}=\boldsymbol{b}(\boldsymbol{y},\boldsymbol{v},\boldsymbol{\zeta}),&&\boldsymbol{x}=\boldsymbol{k}(\boldsymbol{v},\boldsymbol{\zeta}),&&\boldsymbol{b}:\real^{r_{0}}\times\real^{s}\times\real^{n}\rightarrow\real^{s},&&\boldsymbol{k}:\real^{s}\times\real^{n}\rightarrow\real^{s_{0}}.
\end{aligned}
\end{equation}
It should be noted that these equations can be obtained by setting~$\boldsymbol{a}(\boldsymbol{u},\boldsymbol{v},\boldsymbol{\xi})=\boldsymbol{u}-\boldsymbol{f}(\boldsymbol{u},\boldsymbol{v},\boldsymbol{\xi})$ and~$\boldsymbol{b}(\boldsymbol{u},\boldsymbol{v},\boldsymbol{\zeta})=\boldsymbol{v}-\boldsymbol{g}(\boldsymbol{u},\boldsymbol{v},\boldsymbol{\zeta})$, but that alternative reformulations, such as those involving a direct solution of the subproblems or of their linear approximations, are often better adapted.
We then consider the solution of the stochastic coupled model by a \textit{Gauss-Seidel} iterative method using suitable initial values~$\boldsymbol{u}^{0}$, $\boldsymbol{v}^{0}$, and~$\boldsymbol{x}^{0}=\boldsymbol{k}(\boldsymbol{v}^{0},\boldsymbol{\zeta})$ as follows:
\begin{equation}
\label{eq:couplingAGS12}\begin{aligned}
&\boldsymbol{u}^{\ell}=\boldsymbol{a}\big(\boldsymbol{u}^{\ell-1},\boldsymbol{x}^{\ell-1},\boldsymbol{\xi}\big),&&\quad\quad\quad\boldsymbol{y}^{\ell}=\boldsymbol{h}(\boldsymbol{u}^{\ell},\boldsymbol{\xi}),\\
&\boldsymbol{v}^{\ell}=\boldsymbol{b}\big(\boldsymbol{y}^{\ell},\boldsymbol{v}^{\ell-1},\boldsymbol{\zeta}\big),&&\quad\quad\quad\boldsymbol{x}^{\ell}=\boldsymbol{k}(\boldsymbol{v}^{\ell},\boldsymbol{\zeta}).
\end{aligned}
\end{equation}
This is not the only partitioned iterative method available; however, for simplicity, we employ only this method in this work.
It should be noted that although we implement the proposed methodology using the Gauss-Seidel iterative method, one can readily use the proposed methodology with other iterative methods such as Jacobi, relaxation, and Newton methods.

Further, it should be noted that the successive approximations determined by the iterative method~(\ref{eq:couplingAGS12}) can be constructed as random variables of the following form:
\begin{equation}
\label{eq:usefulnessbefore12}
\begin{aligned}
&{\boldsymbol{u}}{}^{\ell}(\theta)\equiv{\boldsymbol{u}}{}^{\ell}\big(\boldsymbol{\xi}(\theta),\boldsymbol{\zeta}(\theta)\big),\\
&{\boldsymbol{v}}{}^{\ell}(\theta)\equiv{\boldsymbol{v}}{}^{\ell}\big(\boldsymbol{\xi}(\theta),\boldsymbol{\zeta}(\theta)\big);
\end{aligned}
\end{equation}
i.e., $\boldsymbol{u}{}^{\ell}$ and~$\boldsymbol{v}{}^{\ell}$ can be constructed as transformations of the input random variables~$\boldsymbol{\xi}=(\xi_{1},\ldots,\xi_{m})$ and~$\boldsymbol{\zeta}=(\zeta_{1},\ldots,\zeta_{n})$.
The random variables~$\boldsymbol{u}{}^{\ell}$ and~$\boldsymbol{v}{}^{\ell}$ thus exist in a solution space of stochastic dimension $m+n$.

Finally, we acknowledge that analytical conditions under which the sequence of approximations generated by (\ref{eq:couplingAGS12}) is ensured to converge can be established.  
However, a detailed treatment of this issue is beyond the scope of this paper, and we only note that the Banach contraction-mapping theorem and the attracting fixed-point theorem~\citep{zeidler1985} allow one to establish convergence properties under global Lipschitz continuity and local differentiability conditions.

\subsection{Dimension reduction}\label{sec:dimensionreduction}
We believe that exchanged information often resides in a considerably lower dimensional space than the input sources of uncertainty themselves.
Therefore, we investigate the effectiveness of dimension-reduction techniques for the representation of the exchanged information.
Rather than exchanging the coupling variables $\boldsymbol{x}^{\ell}$ and $\boldsymbol{y}^{\ell}$ and the solution variables $\boldsymbol{u}^{\ell}$ and $\boldsymbol{v}^{\ell}$ in their original form, we propose to approximate these random variables by a truncated KL decomposition as they pass from subproblem to subproblem and from iteration to iteration.
We specifically consider the solution of the stochastic coupled model by a Gauss-Seidel iterative method that accommodates a dimension-reduction technique as follows:
\begin{equation}
\label{eq:couplingAGSred12}\begin{aligned}
&\hat{\boldsymbol{u}}{}^{\ell}=\boldsymbol{a}\big(\hat{\boldsymbol{u}}{}^{\ell-1,e},\hat{\boldsymbol{x}}{}^{\ell-1,e},\boldsymbol{\xi}\big),&&\quad\quad\quad\hat{\boldsymbol{y}}{}^{\ell}=\boldsymbol{h}(\hat{\boldsymbol{u}}{}^{\ell},\boldsymbol{\xi}),\\
&\hat{\boldsymbol{v}}{}^{\ell}=\boldsymbol{b}\big(\hat{\boldsymbol{y}}{}^{\ell,d},\hat{\boldsymbol{v}}^{\ell-1,d},\boldsymbol{\zeta}\big),&&\quad\quad\quad\hat{\boldsymbol{x}}{}^{\ell}=\boldsymbol{k}(\hat{\boldsymbol{v}}{}^{\ell},\boldsymbol{\zeta}),
\end{aligned}
\end{equation}
where $\boldsymbol{q}^{\ell,d}=[\hat{\boldsymbol{y}}{}^{\ell,d};\hat{\boldsymbol{v}}{}^{\ell-1,d}]$ and $\boldsymbol{r}^{\ell,e}=[\hat{\boldsymbol{u}}{}^{\ell-1,e};\hat{\boldsymbol{x}}{}^{\ell-1,e}]$ are truncated KL decompositions of $\boldsymbol{q}^{\ell}=[\hat{\boldsymbol{y}}{}^{\ell};\hat{\boldsymbol{v}}{}^{\ell-1}]$ and $\boldsymbol{r}^{\ell}=[\hat{\boldsymbol{u}}{}^{\ell-1};\hat{\boldsymbol{x}}{}^{\ell-1}]$, respectively, which read as follows:
\begin{equation}
\label{eq:KLqqqq}\begin{aligned}
&\boldsymbol{q}^{\ell,d}=\overline{\boldsymbol{q}}{}^{\ell}+\sum_{j=1}^{d}\sqrt{\lambda_{j}^{\ell}}\eta_{j}^{\ell}\boldsymbol{\phi}^{j,\ell},\\
&\boldsymbol{r}^{\ell,e}=\overline{\boldsymbol{r}}{}^{\ell}+\sum_{j=1}^{e}\sqrt{\kappa_{j}^{\ell}}\iota_{j}^{\ell}\boldsymbol{\psi}^{j,\ell}.
\end{aligned}
\end{equation}
These decompositions are described in detail in later sections.  
\clearpage 

It should be noted that~(\ref{eq:KLqqqq}) provides a combined reduced-dimensional representation of~$\hat{\boldsymbol{y}}{}^{\ell}$ and~$\hat{\boldsymbol{v}}{}^{\ell-1}$ in terms of a single set of reduced random variables~$\boldsymbol{\eta}^{\ell}=(\eta_{1}^{\ell},\ldots,\eta_{d}^{\ell})$ and a combined reduced-dimensional representation of~$\hat{\boldsymbol{u}}{}^{\ell-1}$ and~$\hat{\boldsymbol{x}}{}^{\ell-1}$ in terms of a single set of reduced random variables~$\boldsymbol{\iota}^{\ell}=(\iota_{1}^{\ell},\ldots,\iota_{e}^{\ell})$.
However, this is not the only construction of a reduced-dimensional representation that could be considered.  
One can readily use the proposed methodology with other dimension-reduction techniques such as those involving the construction of a separate reduced-dimensional representation of the coupling and solution variables, with each representation having its own reduced random variables.  

Truncation of KL decompositions most often results in approximation errors, and we thus use a hat superscript to distinguish the successive approximations determined by~(\ref{eq:couplingAGSred12}) from those determined by~(\ref{eq:couplingAGS12}).
Further, we say that exchanged information has a \textit{low effective stochastic dimension} when its KL decomposition (\ref{eq:KLqqqq}) can be truncated after a few terms while maintaining sufficient accuracy. 
Random variables can be expected to have a low effective stochastic dimension when their components exhibit significant statistical correlation.

Owing to the dimension reduction, the successive approximations determined by the iterative method~(\ref{eq:couplingAGSred12}) can be constructed as random variables of the following form:
\begin{equation}
\label{eq:usefulness12}\begin{aligned}
&\hat{\boldsymbol{u}}{}^{\ell}(\theta)\equiv\hat{\boldsymbol{u}}{}^{\ell}\big(\boldsymbol{\xi}(\theta),\boldsymbol{\iota}^{\ell}(\theta)\big),\\
&\hat{\boldsymbol{v}}{}^{\ell}(\theta)\equiv\hat{\boldsymbol{v}}{}^{\ell}\big(\boldsymbol{\eta}^{\ell}(\theta),\boldsymbol{\zeta}(\theta)\big);
\end{aligned}
\end{equation}
i.e., $\hat{\boldsymbol{u}}{}^{\ell}$ can be constructed as a transformation of the input random variables~$\boldsymbol{\xi}=(\xi_{1},\ldots,\xi_{m})$ and the reduced random variables~$\boldsymbol{\iota}^{\ell}=(\iota_{1}^{\ell},\ldots,\iota_{e}^{\ell})$ and~$\hat{\boldsymbol{v}}{}^{\ell}$ can be constructed as a transformation of~$\boldsymbol{\eta}^{\ell}=(\eta_{1}^{\ell},\ldots,\eta_{d}^{\ell})$ and~$\boldsymbol{\zeta}=(\zeta_{1},\ldots,\zeta_{n})$.
The random variable~$\hat{\boldsymbol{u}}{}^{\ell}$ thus exists in a space of stochastic dimension~$m+e$ and $\hat{\boldsymbol{v}}{}^{\ell}$ exists in a space of stochastic dimension~$d+n$.

Finally, although our notations do not express a potential dependence of $d$ and $e$ on $\ell$, it should be noted that the reduced dimensions can be allowed to depend on the iteration. 

\subsection{Effectiveness of the proposed dimension-reduction methodology}
The key feature of the proposed methodology is that it enables a solution of the subproblems in a reduced-dimensional space when the exchanged information has a low effective stochastic dimension.
Specifically, a solution in a reduced-dimensional space is enabled when the reduced dimensions can be selected such that $d<m$ and $e<n$ while maintaining sufficient accuracy; refer to~(\ref{eq:usefulnessbefore12}) and~(\ref{eq:usefulness12}).
This benefit is of particular significance for implementations of stochastic coupled models using stochastic expansion methods.
These methods suffer from a curse of dimensionality in that their computational cost increases quickly with an increase in the stochastic dimension.
The proposed methodology addresses the curse of dimensionality by mitigating the increase in stochastic dimension when information is exchanged.
Although this prospect is the main motivation for this work, we limit ourselves in this paper to the description and analysis of the dimension reduction itself, and we defer the presentation of algorithms that exploit the dimension reduction to achieve computational gains to a later paper.

The proposed methodology also has the potential of providing insight into the structure and evolution of the information that is exchanged between the subproblems and the iterations.

Finally, 
it should be noted that the proposed methodology can readily be adapted to meet various requirements of specific applications.  
One could, for instance, use the KL decomposition to represent only those exchanged random variables that are of low effective stochastic dimension and adopt an alternative probabilistic representation for the remaining variables.

\section{Karhunen-Loeve decomposition}\label{sec:kl}

Here, we recall the use of the KL decomposition to construct a reduced-dimensional representation of a random variable~$\boldsymbol{q}$ that is defined on a probability triple~$(\Theta,\mathcal{T},P)$, takes values in a Euclidean space~$\real^{w}$, and is of the second order:
\begin{equation}
\int_{\Theta}\vectornorm{\boldsymbol{q}}^{2}dP<+\infty,
\end{equation}
where~$\vectornorm{\cdot}$ denotes the Euclidean norm.

\subsection{Interpretation}
In the context of the proposed methodology, we think of~$\boldsymbol{q}$ as a random variable that passes from one subproblem to another or from one iteration to the next.
With reference to~(\ref{eq:KLqqqq}), the random variable~$\boldsymbol{q}$ could collect the components of $\hat{\boldsymbol{y}}$ and $\hat{\boldsymbol{v}}$, with $w=r_{0}+s$, or the components of $\hat{\boldsymbol{u}}$ and $\hat{\boldsymbol{x}}$, with $w=r+s_{0}$, at a specific iteration. 

\subsection{Second-order descriptors}
The mean vector~$\overline{\boldsymbol{q}}$ and the covariance matrix~$\boldsymbol{C}_{\boldsymbol{q}}$ of the second-order random variable~$\boldsymbol{q}$ are defined as the $w$-dimensional vector and square matrix such that
\begin{align}
\label{eq:meanu}\overline{\boldsymbol{q}}&=\int_{\Theta}\boldsymbol{q}dP,\\
\label{eq:covaru}\boldsymbol{C}_{\boldsymbol{q}}&=\int_{\Theta}(\boldsymbol{q}-\overline{\boldsymbol{q}})(\boldsymbol{q}-\overline{\boldsymbol{q}})^{\mathrm{T}}dP.
\end{align}

\subsection{Reduced-dimensional representation}
Because~$\boldsymbol{C}_{\boldsymbol{q}}$ is a symmetric and positive semidefinite matrix, the solution of the eigenproblem
\begin{equation}
\boldsymbol{C}_{\boldsymbol{q}}\boldsymbol{\phi}^{j}=\lambda_{j}\boldsymbol{\phi}^{j},\label{eq:eigenproblem}
\end{equation}
provides a set of~$w$ eigenvalues~$\lambda_{1}\geq\lambda_{2}\geq\ldots\geq\lambda_{w}\geq 0$ in addition to~$w$ eigenvectors~$\boldsymbol{\phi}^{1},\ldots,\boldsymbol{\phi}^{w}$, which constitute an orthonormal basis of $\real^{w}$ such that
\begin{equation}
(\boldsymbol{\phi}^{i})^{\mathrm{T}}\boldsymbol{\phi}^{j}=\delta_{ij},\label{eq:normalitykl1}
\end{equation}
where $\delta_{ij}$ is the Kronecker delta, which is equal to 1 if $i=j$ and 0 otherwise.
The KL decomposition of~$\boldsymbol{q}$ is then given by
\begin{equation}
\boldsymbol{q}=\overline{\boldsymbol{q}}+\sum_{j=1}^{w}\sqrt{\lambda_{j}}\eta_{j}\boldsymbol{\phi}^{j},\label{eq:KLvector}
\end{equation}
where the~$\eta_{j}$ are random variables defined on~$(\Theta,\mathcal{T},P)$, with values in~$\real$, such that
\begin{equation}
\eta_{j}=\frac{1}{\sqrt{\lambda_{j}}}(\boldsymbol{q}-\overline{\boldsymbol{q}})^{\mathrm{T}}\boldsymbol{\phi}^{j}.
\end{equation}
Further, the~$\eta_{j}$ are zero-mean and uncorrelated:
\begin{align}
\label{eq:eta1}&\int_{\Theta}\eta_{j}dP=0,\\
\label{eq:eta2}&\int_{\Theta}\eta_{i}\eta_{j}dP=\delta_{ij}.
\end{align}
The truncation of~(\ref{eq:KLvector}) after~$d$ terms provides a reduced-dimensional representation as follows:
\begin{equation}
\boldsymbol{q}^{d}=\overline{\boldsymbol{q}}+\sum_{j=1}^{d}\sqrt{\lambda_{j}}\eta_{j}\boldsymbol{\phi}^{j},\label{eq:KLred}
\end{equation}
in which the truncation error, because of the orthonormality properties~(\ref{eq:normalitykl1}) and~(\ref{eq:eta2}), satisfies
\begin{equation}
\int_{\Theta}\vectornorm{\boldsymbol{q}-\boldsymbol{q}^{d}}^{2}dP=\sum_{j=d+1}^{w}\lambda_{j}.\label{eq:KLaccuracy}
\end{equation}
Equality~(\ref{eq:KLaccuracy}) indicates that the accuracy of~(\ref{eq:KLred}) can be improved systematically by increasing the number of retained terms.
Further, for any orthonormal basis~$\{\boldsymbol{e}^{1},\ldots,\boldsymbol{e}^{w}\}$ of~$\real^{w}$, the random variable~$\boldsymbol{q}$ can be expanded as~$\boldsymbol{q}=\overline{\boldsymbol{q}}+\sum_{j=1}^{w}(\boldsymbol{q}-\overline{\boldsymbol{q}})^{\mathrm{T}}\boldsymbol{e}^{j}\boldsymbol{e}^{j}$.
The error~$\boldsymbol{\varepsilon}_{d}=\sum_{j=d+1}^{w}(\boldsymbol{q}-\overline{\boldsymbol{q}})^{\mathrm{T}}\boldsymbol{e}^{j}\boldsymbol{e}^{j}$ introduced owing to the truncation of this expansion after~$d$ terms satisfies
\begin{align}
\notag\int_{\Theta}\vectornorm{\boldsymbol{\varepsilon}_{d}}^{2}dP&=\int_{\Theta}\bigg(\sum_{j=d+1}^{w}(\boldsymbol{q}-\overline{\boldsymbol{q}})^{\mathrm{T}}\boldsymbol{e}^{j}\boldsymbol{e}^{j}\bigg)^{\mathrm{T}}\bigg(\sum_{j=d+1}^{w}(\boldsymbol{q}-\overline{\boldsymbol{q}})^{\mathrm{T}}\boldsymbol{e}^{j}\boldsymbol{e}^{j}\bigg)dP\\
\notag&=\int_{\Theta}\sum_{j=d+1}^{w}\Big((\boldsymbol{q}-\overline{\boldsymbol{q}})^{\mathrm{T}}\boldsymbol{e}^{j}\Big)^{2}dP\\
&=\sum_{j=d+1}^{w}(\boldsymbol{e}^{j})^{\mathrm{T}}\boldsymbol{C}_{\boldsymbol{q}}\boldsymbol{e}^{j}.
\end{align}
This expression indicates that the reduced-dimensional representation~(\ref{eq:KLred}) is optimal because from among all decompositions with~$d$ orthonormal basis vectors, this representation minimizes the mean-square norm of the approximation error, as shown below:
\begin{equation}
\int_{\Theta}\vectornorm{\boldsymbol{q}-\boldsymbol{q}^{d}}^{2}dP=\min_{\substack{\{\boldsymbol{e}^{1},\ldots,\boldsymbol{e}^{d}\}\in\real^{w}\\(\boldsymbol{e}^{i})^{\mathrm{T}}\boldsymbol{e}^{j}=\delta_{ij}}}\int_{\Theta}\bigg\|\boldsymbol{q}-\bigg(\overline{\boldsymbol{q}}+\sum_{j=1}^{d}(\boldsymbol{q}-\overline{\boldsymbol{q}})^{\mathrm{T}}\boldsymbol{e}^{j}\boldsymbol{e}^{j}\bigg)\bigg\|^{2}dP.\label{eq:optimality}
\end{equation}

\section{Proposed adaptation of the Karhunen-Loeve decomposition}\label{sec:klb}
When our dimension-reduction methodology is applied to a stochastic coupled model that has been discretized in space and time, this methodology repeatedly requires the reduction of random variables with values in a finite-dimensional Euclidean space.
The standard KL decomposition, recalled in the previous section, allows random variables with values in a finite-dimensional Euclidean space to be reduced in a manner that is optimal in the natural mean-square norm.
However, it is not always suitable to use a dimension-reduction technique that aims at maintaining accuracy in the natural mean-square norm; for instance, when the samples of the random variables to be reduced have space or time derivatives that also carry important pieces of information, it may be preferable to use an alternative dimension-reduction technique that also gives weight to maintaining accuracy in these derivatives.
The significance of the function values, derivatives, and other hallmarks associated with the random variables to be reduced can be expected to be determined by the function-analytic structure that the stochastic model exhibited prior to its discretization.
In this section, we present an adaptation of the KL decomposition to address this issue.
Our adaptation allows a random variable that solves a space-time discretized stochastic model to be reduced in a manner that maintains consistency with the function-analytic structure exhibited by the stochastic model before discretization.   

\subsection{Methodology}
In the research area of stochastic modeling and analysis, the KL decomposition is known best as a tool for the reduction of stochastic processes that describe uncertain \textit{fields of coefficients and boundary conditions} of stochastic partial differential equations.
In such applications of the KL decomposition, the classical description of stochastic processes as indexed collections of random variables is used.
In this work, we rather intend to reduce the \textit{solution} of stochastic models or to reduce coupling variables that depend on this solution through a specific mapping.
However, the solution of many stochastic models, including that of many stochastic partial differential equations, is not always amenable to a description as a classical stochastic process and is often better described as a random variable that takes its values in a function space.
As we have already mentioned, the structure of this function space can be expected to determine the significance of the function values, the derivatives, and the various other hallmarks of the solution that needs to be reduced.
A detailed explanation of the distinction between classical stochastic processes and function-space-valued random variables can be found in~\citep{kree1983,redhorse2009} and the references therein. 
This distinction has repercussions on the manner in which the second-order descriptors are defined and the properties of these second-order descriptors can be exploited to construct a reduced-dimensional representation.    
Thus, the approach adopted in this section is as follows.
First, we describe the construction of a KL-type decomposition of a random variable with values in a function space.
Then, we describe a construction of a KL-type decomposition of a random variable with values in a discrete approximation of this function space such that consistency with the structure of this function space is maintained.

\subsection{KL-type decomposition of function-space-valued random variables}\label{sec:klfunction}
Let~$\boldsymbol{q}$ be a random variable defined on a probability triple $(\Theta,\mathcal{T},P)$ with values in a separable Hilbert space~$H$.   
Let~$\langle\cdot,\cdot\rangle_{H}$ denote the inner product on~$H$ and~$\vectornorm{\cdot}_{H}=\sqrt{\langle\cdot,\cdot\rangle_{H}}$ denote the norm induced by the inner product.
Let the random variable~$\boldsymbol{q}$ be of the second order:
\begin{equation}
\int_{\Theta}\vectornorm{\boldsymbol{q}}^{2}_{H}dP<+\infty.\label{eq:secondorderb}
\end{equation}
This level of abstraction is very general: in addition to the standard Euclidean spaces of vectors and matrices, examples of~$H$ include spaces of square-integrable functions, spaces of sequences, and certain Sobolev spaces, among many other possibilities. 

\subsubsection{Second-order descriptors.}\label{subsec:kl3}
The manner in which second-order descriptors and their properties are defined for random variables with values in an infinite-dimensional Hilbert space is more complicated than the manner in which they are defined for random variables with values in a finite-dimensional Euclidean space.
All definitions used here are consistent with those given in references~\citep{vakhania1987,kree1983,grenander2008,redhorse2009}. 
Now, the mean is defined as the linear function~$m_{\boldsymbol{q}}$ from~$H$ into~$\real$ such that
\begin{equation}
m_{\boldsymbol{q}}(\boldsymbol{p})=\int_{\Theta}\langle\boldsymbol{q},\boldsymbol{p}\rangle_{H}dP,\quad\forall\boldsymbol{p}\in H,\label{eq:meanq1}
\end{equation}
and the covariance is defined as the bilinear function~$c_{\boldsymbol{q}}$ from~$H\times H$ into~$\real$ such that
\begin{equation}
c_{\boldsymbol{q}}(\boldsymbol{p},\boldsymbol{r})=\int_{\Theta}\big(\langle\boldsymbol{q},\boldsymbol{p}\rangle_{H}-m_{\boldsymbol{q}}(\boldsymbol{p})\big)\big(\langle\boldsymbol{q},\boldsymbol{r}\rangle_{H}-m_{\boldsymbol{q}}(\boldsymbol{r})\big)dP,\quad\forall\boldsymbol{p},\boldsymbol{r}\in H.\label{eq:meancovar1}
\end{equation}
We observe that these definitions are consistent with the function-analytic membership of the random variable~$\boldsymbol{q}$, in that the membership of the realizations of~$\boldsymbol{q}$ in the Hilbert space~$H$ implies that~$\boldsymbol{q}$ is analyzed most naturally through its projection onto fixed elements or basis vectors of~$H$.
Indeed, the mean function associates to any fixed element~$\boldsymbol{p}$ of~$H$ the mean of the random variable~$\langle\boldsymbol{q},\boldsymbol{p}\rangle_{H}$ obtained by projecting~$\boldsymbol{q}$ onto~$\boldsymbol{p}$, and the covariance function associates to any fixed pair of elements~$\boldsymbol{p}$ and~$\boldsymbol{r}$ of~$H$ the covariance of the random variables~$\langle\boldsymbol{q},\boldsymbol{p}\rangle_{H}$ and~$\langle\boldsymbol{q},\boldsymbol{r}\rangle_{H}$.
Clearly, the covariance function is symmetric:
\begin{equation}
c_{\boldsymbol{q}}(\boldsymbol{p},\boldsymbol{r})=c_{\boldsymbol{q}}(\boldsymbol{r},\boldsymbol{p}),\quad\forall\boldsymbol{p},\boldsymbol{r}\in H.
\end{equation}
Further, it is also positive:
\begin{equation}
c_{\boldsymbol{q}}(\boldsymbol{p},\boldsymbol{p})\geq 0,\quad\forall\boldsymbol{p}\in H.
\end{equation}
Equation~(\ref{eq:secondorderb}) and H\"{o}lder's inequality indicate that $m_{\boldsymbol{q}}$ is continuous.
Thus, by Riesz's representation theorem, a unique vector~$\overline{\boldsymbol{q}}{}$ exists in~$H$, i.e. the mean vector, such that
\begin{equation}
\langle\overline{\boldsymbol{q}},\boldsymbol{p}\rangle_{H}=m_{\boldsymbol{q}}(\boldsymbol{p}),\quad\forall\boldsymbol{p}\in H.\label{eq:meanq2}
\end{equation}
Likewise, equation~(\ref{eq:secondorderb}) and H\"{o}lder's inequality indicate that $c_{\boldsymbol{q}}$ is continuous; and therefore, by Riesz's representation theorem, there exists~\citep{reed1980} a unique linear continuous operator~$\mathcal{C}_{\boldsymbol{q}}$ from~$H$ into~$H$, i.e. the covariance operator, such that
\begin{equation}
\left\langle\mathcal{C}_{\boldsymbol{q}}(\boldsymbol{p}),\boldsymbol{r}\right\rangle_{H}=c_{\boldsymbol{q}}(\boldsymbol{p},\boldsymbol{r}),\quad\forall\boldsymbol{p},\boldsymbol{r}\in H. \label{eq:meancovar2}
\end{equation}
Let~$\{\boldsymbol{e}^{j}\}_{j=1}^{\infty}$ be any complete orthonormal basis of~$H$.  
Because the Hilbert-Schmidt norm,
\begin{align}
\notag\sqrt{\sum_{i,j=1}^{\infty}\left\langle\mathcal{C}_{\boldsymbol{q}}(\boldsymbol{e}^{i}),\boldsymbol{e}^{j}\right\rangle_{H}^{2}}&=\sqrt{\sum_{i,j=1}^{\infty}\left(\int_{\Theta}\langle\boldsymbol{q}-\overline{\boldsymbol{q}}{},\;\boldsymbol{e}^{i}\rangle_{H}\langle\boldsymbol{q}-\overline{\boldsymbol{q}}{},\;\boldsymbol{e}^{j}\rangle_{H}dP\right)^{2}}\\
\notag&\leq\sum_{j=1}^{\infty}\int_{\Theta}\langle\boldsymbol{q}-\overline{\boldsymbol{q}}{},\;\boldsymbol{e}^{j}\rangle_{H}^{2}dP\\
&=\int_{\Theta}\vectornorm{\boldsymbol{q}}_{H}^{2}dP-\vectornorm{\overline{\boldsymbol{q}}}_{H}^{2},
\end{align}
is bounded owing to~(\ref{eq:secondorderb}), the covariance operator is~\citep{reed1980} a Hilbert-Schmidt operator.
Moreover, because the covariance function is symmetric, the covariance operator is~\citep{kree1983} self-adjoint:
\begin{equation}
\mathcal{C}_{\boldsymbol{q}}=\mathcal{C}_{\boldsymbol{q}}^{\mathrm{t}},
\end{equation}
where the operator~$\mathcal{C}_{\boldsymbol{q}}^{\mathrm{t}}$ from~$H$ into~$H$ is the adjoint operator of~$\mathcal{C}_{\boldsymbol{q}}$ such that
\begin{equation}
\langle\mathcal{C}_{\boldsymbol{q}}(\boldsymbol{p}),\boldsymbol{r}\rangle_{H}=\langle\boldsymbol{p},\mathcal{C}_{\boldsymbol{q}}^{\mathrm{t}}(\boldsymbol{r})\rangle_{H},\quad\forall\boldsymbol{p},\boldsymbol{r}\in H.
\end{equation}
Because the covariance function is positive, the covariance operator is also positive:
\begin{equation}
\langle\mathcal{C}_{\boldsymbol{q}}(\boldsymbol{p}),\boldsymbol{p}\rangle_{H}\geq 0,\quad\forall\boldsymbol{p}\in H.
\end{equation}

\subsubsection{Reduced-dimensional representation.}\label{subsec:kl4}
Because~$\mathcal{C}_{\boldsymbol{q}}$ is a Hilbert-Schmidt operator and thus compact and because~$\mathcal{C}_{\boldsymbol{q}}$ is self-adjoint and positive, the solution of the eigenproblem
\begin{equation}
\mathcal{C}_{\boldsymbol{q}}(\boldsymbol{\phi}^{j})=\lambda_{j}\boldsymbol{\phi}^{j}\label{eq:eigenproblemmmeme}
\end{equation}
provides a denumerable set of eigenvalues~$\{\lambda_{j}\}_{j=1}^{\infty}$ in addition to a corresponding denumerable set of eigenmodes~$\{\boldsymbol{\phi}^{j}\}_{j=1}^{\infty}$; the eigenvalues are positive and square-summable in that~$\sum_{j=1}^{\infty}\lambda_{j}^{2}<+\infty$ and thus~$\lambda_{j}$ tends to 0 as~$j$ goes to infinity; and the eigenmodes constitute a Hilbertian basis of~$H$ such that
\begin{equation}
\langle\boldsymbol{\phi}^{i},\boldsymbol{\phi}^{j}\rangle_{H}=\delta_{ij}.\label{eq:propkl1}
\end{equation}
The KL decomposition of the random variable $\boldsymbol{q}$ is then expressed as follows:
\begin{equation}
\boldsymbol{q}=\overline{\boldsymbol{q}}+\sum_{j=1}^{\infty}\sqrt{\lambda_{j}}\eta_{j}\boldsymbol{\phi}^{j},\label{eq:hkl}
\end{equation}
where the~$\eta_{j}$ are random variables defined on~$(\Theta,\mathcal{T},P)$, with values in~$\real$, such that
\begin{equation}
\eta_{j}=\frac{1}{\sqrt{\lambda_{j}}}\langle\boldsymbol{q}-\overline{\boldsymbol{q}}{},\boldsymbol{\phi}^{j}\rangle_{H},
\end{equation}
and are zero-mean and uncorrelated, as given by~(\ref{eq:eta1})--(\ref{eq:eta2}).
By the Hilbert-space orthogonal decomposition theorem~\citep{reed1980}, the decomposition converges strongly:
\begin{equation}
\lim_{d\rightarrow\infty}\int_{\Theta}\bigg\|\boldsymbol{q}-\bigg(\overline{\boldsymbol{q}}+\sum_{j=1}^{d}\sqrt{\lambda_{j}}\eta_{j}\boldsymbol{\phi}^{j}\bigg)\bigg\|_{H}^{2}dP=0.
\end{equation}

\subsection{KL-type decomposition of function-space-valued random variables after discretization}
Now, let us consider a random variable~$\boldsymbol{q}^{w}$ that is still defined on the probability triple $(\Theta,\mathcal{T},P)$ but takes, this time, its values in a finite-dimensional subspace~$H^{w}$ of~$H$.
Let this random variable~$\boldsymbol{q}^{w}$ be represented by an expansion of the following form:
\begin{equation}
\boldsymbol{q}^{w}=\sum_{j=1}^{w}q_{j}\boldsymbol{n}^{j},\label{eq:qnj}
\end{equation}
where the vectors~$\boldsymbol{n}^{1},\ldots,\boldsymbol{n}^{w}$ constitute a basis in~$H^{w}$, and let it be of the second order:
\begin{equation}
\int_{\Theta}\vectornorm{\boldsymbol{q}^{w}}_{H}^{2}dP<+\infty.\label{eq:secondorderbbb}
\end{equation}
Let~$\boldsymbol{q}=(q_{1},\ldots,q_{w})$ now be the random vector with values in~$\real^{w}$ which collects the random coordinates in expansion (\ref{eq:qnj}).
In the remainder of this section, we will propose an adaptation of the KL decomposition which allows a reduced-dimensional representation of~$\boldsymbol{q}$ to be obtained in a manner that is consistent with Hilbertian projections in $H$.

\subsubsection{Second-order descriptors.}
Similarly to the KL decomposition of Sec.~\ref{sec:kl}, the first step in the construction of our adaptation of the KL decomposition consists in determining the mean vector~$\overline{\boldsymbol{q}}$ and the covariance matrix~$\boldsymbol{C}_{\boldsymbol{q}}$ of~$\boldsymbol{q}$, as defined by~(\ref{eq:meanu})--(\ref{eq:covaru}). 

\subsubsection{Reduced-dimensional representation.}
Whereas the KL decomposition of Sec.~\ref{sec:kl} was obtained by solving the eigenproblem determined solely by the covariance matrix, our adaptation of the KL decomposition is rather obtained by solving the generalized eigenproblem
\begin{equation}
\boldsymbol{W}^{\mathrm{T}}\boldsymbol{C}_{\boldsymbol{q}}\boldsymbol{W}\boldsymbol{\phi}^{j}=\lambda_{j}\boldsymbol{W}\boldsymbol{\phi}^{j},\label{eq:eigenproblemeee}
\end{equation}
which features not only the covariance matrix~$\boldsymbol{C}_{\boldsymbol{q}}$ but also the Gram matrix,
\begin{equation}
\boldsymbol{W}=\begin{bmatrix}
\langle\boldsymbol{n}^{1},\boldsymbol{n}^{1}\rangle_{H} & \ldots & \langle\boldsymbol{n}^{1},\boldsymbol{n}^{w}\rangle_{H}\\
\vdots & & \vdots \\
\langle\boldsymbol{n}^{w},\boldsymbol{n}^{1}\rangle_{H} & \ldots & \langle\boldsymbol{n}^{w},\boldsymbol{n}^{w}\rangle_{H}
\end{bmatrix},
\end{equation}
of the basis vectors~$\boldsymbol{n}^{1},\ldots,\boldsymbol{n}^{w}$ as a weighting matrix; the Gram matrix is the~$w$-dimensional, symmetric, positive definite matrix that collects the inner products of these basis vectors. 
Because~$\boldsymbol{C}_{\boldsymbol{q}}$ is symmetric and positive semidefinite and~$\boldsymbol{W}$ is symmetric and positive definite, the solution of~(\ref{eq:eigenproblemeee}) provides a set of~$w$ eigenvalues~$\lambda_{1}\geq\lambda_{2}\geq\ldots\geq\lambda_{w}\geq 0$ in addition to~$w$ eigenvectors~$\boldsymbol{\phi}^{1},\ldots,\boldsymbol{\phi}^{w}$, which constitute, this time, a~$\boldsymbol{W}$-weighted orthonormal basis of $\real^{w}$:
\begin{equation}
(\boldsymbol{\phi}^{i})^{\mathrm{T}}\boldsymbol{W}\boldsymbol{\phi}^{j}=\delta_{ij}.\label{eq:orthonromalkfa3}
\end{equation}
It should be noted that the eigenvalues and eigenmodes of the generalized eigenproblem~(\ref{eq:eigenproblemeee}) usually differ from those of the eigenproblem~(\ref{eq:eigenproblem}), even though we use the same notations.
As an alternative to~(\ref{eq:KLred}), we then obtain a reduced-dimensional representation~$\boldsymbol{q}^{d}$ of~$\boldsymbol{q}$ as:
\begin{equation}
\boldsymbol{q}^{d}=\overline{\boldsymbol{q}}+\sum_{j=1}^{d}\sqrt{\lambda_{j}}\eta_{j}\boldsymbol{\phi}^{j},\label{eq:reducedw}
\end{equation}
where the~$\eta_{j}$ are, this time, random variables on~$(\Theta,\mathcal{T},P)$, with values in~$\real$, such that
\begin{equation}
\eta_{j}=\frac{1}{\sqrt{\lambda_{j}}}(\boldsymbol{q}-\overline{\boldsymbol{q}})^{\mathrm{T}}\boldsymbol{W}\boldsymbol{\phi}^{j},\label{eq:constr8}
\end{equation}
and are zero-mean and uncorrelated, as given by~(\ref{eq:eta1})--(\ref{eq:eta2}).
Again, it should be noted that the reduced random variables of the decomposition~(\ref{eq:reducedw}) usually differ from those of the decomposition~(\ref{eq:KLred}), even though we use the same notations.
Because of the orthonormality properties, the truncation error incurred by~$\boldsymbol{q}^{d}$ satisfies: 
\begin{equation}
\int_{\Theta}\vectornorm{\boldsymbol{q}-\boldsymbol{q}^{d}}^{2}_{\boldsymbol{W}}dP=\sum_{j=d+1}^{w}\lambda_{j},\label{eq:keypropertyyy}
\end{equation}
where~$\sqrt{\int_{\Theta}\vectornorm{\boldsymbol{p}}_{\boldsymbol{W}}^{2}dP}=\sqrt{\int_{\Theta}\boldsymbol{p}^{\mathrm{T}}\boldsymbol{W}\boldsymbol{p}\,dP}$ for any second-order random variable~$\boldsymbol{p}$ on~$(\Theta,\mathcal{T},P)$ with values in~$\real^{w}$.
It should be noted that whereas the truncation error was gauged by the natural mean-square norm in~(\ref{eq:KLaccuracy}), it is gauged here by the $\boldsymbol{W}$-weighted mean-square norm.

Moreover, following a reasoning that is similar to the one exhibited in Sec.~\ref{sec:kl}, it can readily be shown that~$\boldsymbol{q}^{d}$ is optimal because from among all decompositions with~$d$ orthonormal basis vectors, this representation minimizes the mean-square norm of the approximation error:
\begin{equation}
\int_{\Theta}\vectornorm{\boldsymbol{q}-\boldsymbol{q}^{d}}_{\boldsymbol{W}}^{2}dP=\min_{\substack{\{\boldsymbol{e}^{1},\ldots,\boldsymbol{e}^{d}\}\in\real^{w}\\(\boldsymbol{e}^{i})^{\mathrm{T}}\boldsymbol{W}\boldsymbol{e}^{j}=\delta_{ij}}}\int_{\Theta}\bigg\|\boldsymbol{q}-\bigg(\overline{\boldsymbol{q}}+\sum_{j=1}^{d}(\boldsymbol{q}-\overline{\boldsymbol{q}})^{\mathrm{T}}\boldsymbol{W}\boldsymbol{e}^{j}\boldsymbol{e}^{j}\bigg)\bigg\|_{\boldsymbol{W}}^{2}dP.\label{eq:optimalityyy}
\end{equation}
Again, it should be noted that whereas the reduced-dimensional representation~(\ref{eq:KLred}) achieved optimality in~(\ref{eq:optimality}) in the natural mean-square norm, the reduced-dimensional representation~(\ref{eq:reducedw}) achieves optimality here in the $\boldsymbol{W}$-weighted mean-square norm.

The interest of introducing the Gram matrix of the basis vectors is that the space of second-order random variables with values in~$\real^{w}$ equipped with the norm $\sqrt{\int_{\Theta}\vectornorm{\cdot}_{\boldsymbol{W}}^{2}dP}$ and the space of second-order random variables  with values in~$H^{w}$ equipped with the norm $\sqrt{\int_{\Theta}\vectornorm{\cdot}_{H}^{2}dP}$ share the same structure.
Indeed, the function that maps any second-order random variable~$\boldsymbol{p}$ with values in~$\real^{w}$ onto a corresponding second-order random variable~$\boldsymbol{p}^{w}=\sum_{j=1}^{w}p_{j}\boldsymbol{n}^{j}$ with values in~$H^{w}$ is a structure-preserving linear bijection between these spaces such that
\begin{equation}
\sqrt{\int_{\Theta}\vectornorm{\boldsymbol{p}}_{\boldsymbol{W}}^{2}dP}=\sqrt{\int_{\Theta}\vectornorm{\boldsymbol{p}^{w}}_{H}^{2}dP},\quad\quad\quad\forall\boldsymbol{p}\in L^{2}_{P}(\Theta,\real^{w}).
\end{equation}
As a conclusion, when the random coordinates~$\boldsymbol{q}=(q_{1},\ldots,q_{w})$ of a random variable~$\boldsymbol{q}^{w}=\sum_{j=1}^{w}q_{j}\boldsymbol{n}^{j}$ with values in a function space~$H$ must be reduced, the use of the Gram matrix~$\boldsymbol{W}$ of the basis vectors~$\boldsymbol{n}^{1},\ldots,\boldsymbol{n}^{w}$ as a weighting matrix in the construction of the KL decomposition allows a reduced-dimensional representation~$\boldsymbol{q}^{d}$ of~$\boldsymbol{q}$ to be obtained which is consistent with Hilbertian projections in $H$ and therefore optimal in a norm consistent with the norm of~$H$.

\subsubsection{Relationship with the KL-type decomposition of function-space-valued random variables.}
Let~$\overline{\boldsymbol{q}}{}^{w}$ and~$\mathcal{C}_{\boldsymbol{q}}^{w}$ be the mean and covariance operator of~$\boldsymbol{q}^{w}$, defined in a manner consistent with the definition of second-order descriptors of function-space-valued random variables:
\begin{align}
\big\langle\overline{\boldsymbol{q}}{}^{w},\boldsymbol{p}^{w}\big\rangle_{H}&=\int_{\Theta}\big\langle\boldsymbol{q}^{w},\boldsymbol{p}^{w}\big\rangle_{H}dP,&&\forall\boldsymbol{p}^{w}\in H^{w},\\
\big\langle\mathcal{C}_{\boldsymbol{q}}^{w}\big(\boldsymbol{p}^{w}\big),\boldsymbol{r}^{w}\big\rangle_{H}&=\int_{\Theta}\big\langle\boldsymbol{q}^{w}-\overline{\boldsymbol{q}}{}^{w},\boldsymbol{p}^{w}\big\rangle_{H}\big\langle\boldsymbol{q}^{w}-\overline{\boldsymbol{q}}{}^{w},\boldsymbol{r}^{w}\big\rangle_{H}dP,&&\forall\boldsymbol{p}^{w},\boldsymbol{r}^{w}\in H^{w}.
\end{align}
Owing to the relationship~(\ref{eq:qnj}) between the random coordinates~$\boldsymbol{q}$ and the random variable~$\boldsymbol{q}^{w}$, the second-order descriptors of~$\boldsymbol{q}$ are related to those of~$\boldsymbol{q}^{w}$ as:
\begin{align}
\label{eq:relationship1}\overline{\boldsymbol{q}}{}^{\mathrm{T}}\boldsymbol{W}\boldsymbol{p}&=\big\langle\overline{\boldsymbol{q}}{}^{w},\boldsymbol{p}^{w}\big\rangle_{H},&&\forall\boldsymbol{p}\in\real^{w},\\
\label{eq:relationship2}\boldsymbol{p}^{\mathrm{T}}\boldsymbol{W}^{\mathrm{T}}\boldsymbol{C}_{\boldsymbol{q}}\boldsymbol{W}\boldsymbol{r}&=\big\langle\mathcal{C}_{\boldsymbol{q}}^{w}(\boldsymbol{p}^{w}),\boldsymbol{r}^{w}\big\rangle_{H},&&\forall\boldsymbol{p},\boldsymbol{r}\in\real^{w},
\end{align}
in which the correspondence between~$\boldsymbol{p}$ and~$\boldsymbol{p}^{w}$ and between~$\boldsymbol{r}$ and~$\boldsymbol{r}^{w}$ is such that~$\boldsymbol{p}^{w}=\sum_{j=1}^{w}p_{j}\boldsymbol{n}^{j}$ and~$\boldsymbol{r}^{w}=\sum_{j=1}^{w}r_{j}\boldsymbol{n}^{j}$, and~$\boldsymbol{W}$ is still the Gram matrix of the basis vectors.

Further, the reduced-dimensional representation~$\boldsymbol{q}^{d}$ of the random coordinates~$\boldsymbol{q}$ determines a corresponding reduced-dimensional representation~$\boldsymbol{q}^{d,w}$ of the random variable~$\boldsymbol{q}^{w}$ as follows:
\begin{equation}
\boldsymbol{q}^{d,w}=\overline{\boldsymbol{q}}{}^{w}+\sum_{j=1}^{d}\sqrt{\lambda_{j}}\eta_{j}\boldsymbol{\phi}^{j,w},\label{eq:KLzoveelste}
\end{equation}
where $\overline{\boldsymbol{q}}{}^{w}=\sum_{i=1}^{w}\overline{q}{}_{i}\boldsymbol{n}^{i}$ and $\boldsymbol{\phi}^{j,w}=\sum_{i=1}^{w}\phi_{i}^{j}\boldsymbol{n}^{i}$.
It follows from the relationship~(\ref{eq:relationship1}) between the mean vector of~$\boldsymbol{q}$ and the mean vector of~$\boldsymbol{q}^{w}$ that~$\overline{\boldsymbol{q}}{}^{w}$ in~(\ref{eq:KLzoveelste}) is precisely the mean vector of~$\boldsymbol{q}^{w}$.
Moreover, it follows from the relationship~(\ref{eq:relationship2}) between the covariance matrix of~$\boldsymbol{q}$ and the covariance operator of~$\boldsymbol{q}^{w}$ that the~$\lambda_{j}$ and the~$\boldsymbol{\phi}^{j,w}$ in~(\ref{eq:KLzoveelste}) are precisely the eigenvalues and eigenmodes of the eigenproblem determined by the covariance operator of $\boldsymbol{q}^{w}$:
\begin{equation}
\mathcal{C}_{\boldsymbol{q}}^{w}\big(\boldsymbol{\phi}^{j,w}\big)=\lambda_{j}\boldsymbol{\phi}^{j,w}.\label{eq:constr1}
\end{equation}
Hence, the reduced dimensional representation $\boldsymbol{q}^{d,w}$ of $\boldsymbol{q}^{w}$ deduced from the reduced-dimensional representation $\boldsymbol{q}^{d}$ of $\boldsymbol{q}$ by~(\ref{eq:KLzoveelste}), is precisely the reduced-dimensional representation of~$\boldsymbol{q}^{w}$ that would obtained if the methodology for the construction of a KL-type decomposition of function-space-valued random variables of Sec~\ref{sec:klfunction} were applied to~$\boldsymbol{q}^{w}$.

\subsection{Concluding remarks}
In this section, we described an adaptation of the KL decomposition which permits random variables with values in a finite-dimensional Euclidean space to be reduced in a manner that is optimal in a weighted mean-square norm.
This adaptation is well-suited for the reduction of a random variable that solves a space-time discretized stochastic model, because by appropriately choosing the weighting matrix as the Gram matrix of the discretization basis, a reduced-dimensional representation is obtained which is consistent with the function-analytic structure that the stochastic model exhibited before its discretization.   

\section{Implementation}\label{sec:sec5}
In this section, we provide details on the implementation of the proposed dimension-reduction methodology using stochastic expansion methods.

\subsection{Discretization using stochastic expansion methods}\label{sec:sec51}
Within the context of the model problem of Sec.~\ref{sec:sec4}, let~$\{\psi_{\boldsymbol{\alpha}},\boldsymbol{\alpha}\in\integer^{m+n}\}$ be a Hilbertian basis for the Hilbert space of~$P_{(\boldsymbol{\xi},\boldsymbol{\zeta})}$ square-integrable functions from~$\real^{m+n}$ into~$\real$.
Let this Hilbertian basis be indexed by multi-indices~$\boldsymbol{\alpha}=(\alpha_{1},\ldots,\alpha_{m+n})$ in~$\integer^{m+n}$.
We will employ throughout this work a Hilbertian basis that is constituted of polynomials of increasing total degree, and we will thus refer to it as a Polynomial Chaos (PC) basis. 
Note that this PC basis can be obtained by Gram-Schmidt orthonormalization of a collection of multivariate monomials, or read from tables in the literature~\citep{ghanem1998,ghanem1998b,ghiocel2002,ghanem2003,soize2004,xiu2005,babuska2007} if~$P_{(\boldsymbol{\xi},\boldsymbol{\zeta})}$ is a ``labeled" probability distribution.

The iterative methods~(\ref{eq:couplingAGS12}) and~(\ref{eq:couplingAGSred12}) require at each iteration that the subproblems be solved to obtain updated representations of the solution variables. 
Stochastic expansion methods involve the approximation of the solution variables of the subproblems by PC expansions.
The PC basis provides the approximate representation of the solution to~(\ref{eq:couplingAGS12}) as follows: 
\begin{equation}
\label{eq:pceup}
\begin{aligned}\boldsymbol{u}^{\ell,p}&=\sum_{|\boldsymbol{\alpha}|=0}^{p}\boldsymbol{u}^{\ell}_{\boldsymbol{\alpha}}\psi_{\boldsymbol{\alpha}}(\boldsymbol{\xi},\boldsymbol{\zeta}),\quad\quad\quad\boldsymbol{u}^{\ell}_{\boldsymbol{\alpha}}\in\real^{r},\\
\boldsymbol{v}^{\ell,p}&=\sum_{|\boldsymbol{\alpha}|=0}^{p}\boldsymbol{v}^{\ell}_{\boldsymbol{\alpha}}\psi_{\boldsymbol{\alpha}}(\boldsymbol{\xi},\boldsymbol{\zeta}),\quad\quad\quad\boldsymbol{v}_{\boldsymbol{\alpha}}^{\ell}\in\real^{s},
\end{aligned}
\end{equation}
and the approximate representation of the solution to~(\ref{eq:couplingAGSred12}) as follows:
\begin{equation}
\label{eq:pceupp}
\begin{aligned}\hat{\boldsymbol{u}}{}^{\ell,p}&=\sum_{|\boldsymbol{\alpha}|=0}^{p}\hat{\boldsymbol{u}}{}^{\ell}_{\boldsymbol{\alpha}}\psi_{\boldsymbol{\alpha}}(\boldsymbol{\xi},\boldsymbol{\zeta}),\quad\quad\quad\hat{\boldsymbol{u}}{}^{\ell}_{\boldsymbol{\alpha}}\in\real^{r},\\
\hat{\boldsymbol{v}}{}^{\ell,p}&=\sum_{|\boldsymbol{\alpha}|=0}^{p}\hat{\boldsymbol{v}}{}^{\ell}_{\boldsymbol{\alpha}}\psi_{\boldsymbol{\alpha}}(\boldsymbol{\xi},\boldsymbol{\zeta}),\quad\quad\quad\hat{\boldsymbol{v}}{}^{\ell}_{\boldsymbol{\alpha}}\in\real^{s},
\end{aligned}
\end{equation}
in which the summation from~$|\boldsymbol{\alpha}|=0$ to~$p$ means the summation over the finite subset of multi-indices~$\boldsymbol{\alpha}=(\alpha_{1},\ldots,\alpha_{m+n})$ in~$\integer^{m+n}$ with~$|\boldsymbol{\alpha}|=\alpha_{1}+\ldots+\alpha_{m+n}\leq p$.  

Then, the task of the solution algorithm is to compute the coordinates in these expansions, for which several methods are available, such as embedded projection~\citep{ghanem2003,soize2004}, nonintrusive projection~\citep{soize2004}, and collocation~\citep{ghanem1998,ghanem1998b,ghiocel2002,xiu2005,babuska2007}.
We use the nonintrusive projection method in this work.
Although we implement our dimension-reduction methodology using the nonintrusive projection method, it should be noted that our methodology can readily be adapted to other methods such as embedded projection and collocation.

\subsection{KL decomposition}\label{sec:sec52}
Our methodology hinges on the representation of information by a truncated KL decomposition as it passes from subproblem to subproblem and from iteration to iteration.
The implementation of this KL decomposition requires a probabilistic representation of the random variables to be reduced.
In this work, we elected to represent random variables systematically by PC expansions.
Here, we thus describe the implementation of the KL decomposition of a random variable that is represented by a PC expansion, and we show how this implementation naturally provides in turn a PC expansion of the reduced random variables.

Adopting the notations used in Secs.~\ref{sec:kl} and~\ref{sec:klb}, we consider the reduction of a second-order random variable~$\boldsymbol{q}^{p}$ that takes its values in~$\real^{w}$ and is represented by a PC expansion of the form:
\begin{equation}
\boldsymbol{q}^{p}=\sum_{|\boldsymbol{\alpha}|=0}^{p}\boldsymbol{q}_{\boldsymbol{\alpha}}\psi_{\boldsymbol{\alpha}}(\boldsymbol{\xi},\boldsymbol{\zeta}),\quad\boldsymbol{q}_{\boldsymbol{\alpha}}\in \real^{w}.\label{eq:upcekl}
\end{equation}
We specifically consider the construction of the $\boldsymbol{W}$-weighted KL decomposition presented in Sec.~\ref{sec:klb}; the standard KL decomposition presented in Sec.~\ref{sec:kl} can be recovered easily from the KL decomposition presented in Sec.~\ref{sec:klb} by simply setting the weighting matrix~$\boldsymbol{W}$ equal to the identity matrix.
Because of the orthonormality of the~$\psi_{\boldsymbol{\alpha}}$, the mean $\overline{\boldsymbol{q}}$ and the covariance $\boldsymbol{C}_{\boldsymbol{q}}$ of~$\boldsymbol{q}^{p}$ follow immediately from the PC coordinates:
\begin{align}
\overline{\boldsymbol{q}}&=\boldsymbol{q}_{\boldsymbol{0}},\\
\boldsymbol{C}_{\boldsymbol{q}}&=\sum_{|\boldsymbol{\alpha}|=1}^{p}\boldsymbol{q}_{\boldsymbol{\alpha}}\boldsymbol{q}_{\boldsymbol{\alpha}}^{\mathrm{T}}.
\end{align}   
Then, the solution of the generalized eigenproblem~$\boldsymbol{W}^{\mathrm{T}}\boldsymbol{C}_{\boldsymbol{q}}\boldsymbol{W}\boldsymbol{\Phi}^{j}=\lambda_{j}\boldsymbol{W}\boldsymbol{\Phi}^{j}$ provides the eigenvalues~$\lambda_{1}\geq\lambda_{2}\geq\ldots\lambda_{w}\geq 0$ and the associated eigenmodes~$\boldsymbol{\phi}^{1},\ldots,\boldsymbol{\phi}^{w}$ required to construct a reduced-dimensional representation~$\boldsymbol{q}^{p,d}$ of $\boldsymbol{q}^{p}$ as follows:
\begin{equation}
\boldsymbol{q}^{p}=\overline{\boldsymbol{q}}+\sum_{j=1}^{w}\sqrt{\lambda_{j}}\eta^{p}_{j}\boldsymbol{\phi}^{j},
\end{equation}
where the~$\eta_{j}^{p}$ are random variables with values in~$\real$ such that
\begin{equation}
\eta_{j}^{p}=\frac{1}{\sqrt{\lambda_{j}}}\big(\boldsymbol{q}^{p}-\overline{\boldsymbol{q}}\big)^{\mathrm{T}}\boldsymbol{W}\boldsymbol{\phi}^{j}\label{eq:upcekleta}
\end{equation}
and are zero-mean and uncorrelated.
By substituting~(\ref{eq:upcekl}) in~(\ref{eq:upcekleta}), a representation of each reduced random variable as a PC expansion is immediately obtained: 
\begin{equation}
\eta^{p}_{j}=\sum_{|\boldsymbol{\alpha}|=1}^{p}\eta_{j,\boldsymbol{\alpha}}\psi_{\boldsymbol{\alpha}}(\boldsymbol{\xi},\boldsymbol{\zeta})\quad\text{with}\quad\eta_{j,\boldsymbol{\alpha}}=\frac{1}{\sqrt{\lambda_{j}}}\boldsymbol{q}_{\boldsymbol{\alpha}}^{\mathrm{T}}\boldsymbol{W}\boldsymbol{\phi}^{j},\label{eq:upcerveta}
\end{equation}
thus indicating that the KL decomposition of a PC expansion naturally provides a complete probabilistic characterization of the reduced random variables as a PC expansion.

\subsection{Selection of the reduced dimension}
In this work, we let the solution algorithm automatically adjust the reduced dimension within each subproblem and at each iteration in such a way that a prescribed accuracy level is maintained. 
This iteration-dependent adjustment of the dimension reduction ensures that a persistent accuracy level is maintained as the successive iterations evolve. 
Naturally, and even though the truncation errors introduced owing to the dimension reductions
can be made arbitrarily small by 
retaining a sufficiently large number of terms systematically, these truncation errors will most likely have an effect on the solution of the subproblems and will also propagate to the subsequent approximations generated by the iterative method.
Appendix~I provides a concise formal analysis of the effect of these truncation errors.

\section{Realization for a stochastic multiphysics problem}\label{sec:sec6a}
We will now demonstrate the effectiveness of the proposed methodology through an illustration problem relevant to nuclear reactors.

\subsection{Problem formulation}
\begin{figure}[htp]
  \begin{center}
    \begin{picture}(250,100)(0,0)
      \put(50,50){\line(1,1){40}}
      \put(50,50){\line(-1,-1){40}}
      \put(200,50){\line(1,1){40}}
      \put(200,50){\line(-1,-1){40}}
      \put(10,10){\line(1,0){150}}
      \put(90,90){\line(1,0){150}}
      \put(120,75){\makebox{heat transfer}}
      \put(50,24){\makebox{$\;\;\;\;$neutron transport}}
      \put(80,65){\makebox{with random transmittivity}}
      \put(115,47){\vector(1,1){10}}
      \put(115,47){\vector(-1,-1){10}}
      \put(50,50){\circle{2}}
      \put(200,50){\circle{2}}
      \put(21,48){\makebox{$x=0$}}
      \put(205,48){\makebox{$x=L$}}
      \put(15,70){\makebox{$\frac{dT}{dx}\Big|_{x=0}=0$}}
      \put(235,70){\makebox{$\frac{dT}{dx}\Big|_{x=L}=0$}}
      \put(-30,25){\makebox{$\frac{d\Phi}{dx}\Big|_{x=0}=0$}}
      \put(190,25){\makebox{$\frac{d\Phi}{dx}\Big|_{x=L}=0$}}
    \end{picture}
  \end{center}
  \caption{Schematic representation of the problem.}\label{fig:figure0}
\end{figure}
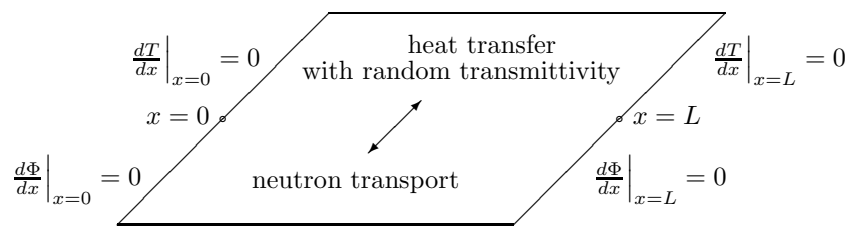
We consider the stationary transport of neutrons in a one-dimensional reactor with temperature feedback~\citep{lamarsh2002}.
Let the reactor occupy an open interval~$]0,L[$ (Fig.~\ref{fig:figure0}).
The problem then involves finding the temperature~$T$ and neutron flux~$\Phi$ such that
\begin{equation}
\label{eq:neutron1c}
\begin{aligned}
&\frac{d}{dx}\left(k\frac{dT}{dx}\right)-h(T-T_{\infty})=-E_{\text{f}}\Sigma_{\text{f}}(T)\Phi,\\
&\frac{d}{dx}\left(D(T)\frac{d\Phi}{dx}\right)-\Big(\Sigma_{\text{a}}(T)-\nu\Sigma_{\text{f}}(T)\Big)\Phi=-s,
\end{aligned}
\end{equation}
under homogeneous Neumann boundary conditions.
The first term on the left-hand side of the heat subproblem represents heat conduction, and the second term represents the transmission of heat to the surroundings; further, the right-hand side represents a distributed heat source proportional to the neutron flux.  
The first term on the left-hand side of the neutronics subproblem represents neutron diffusion, and the second term represents the net effect of the absorption and generation of neutrons; further, the right-hand side represents a distributed neutron source.  
The coefficients~$k$ and~$h$ are the heat conductivity and heat transmittivity, respectively; the temperature $T_{\infty}$ is the ambient temperature; and $\nu$ and~$E_{\text{f}}$ are the number of neutrons and the energy released per fission reaction, respectively.
The coefficients~$D$, $\Sigma_{\text{a}}$ and~$\Sigma_{\text{f}}$ are the neutron diffusion constant, fission cross section, and absorption cross section, respectively; these coefficients depend on the reactor temperature as follows:
\begin{equation}
D\big(T(x)\big)=D_{\text{ref}}\sqrt{\frac{T(x)}{T_{\text{ref}}}},\quad\Sigma_{\text{a}}\big(T(x)\big)=\Sigma_{\text{a,ref}}\sqrt{\frac{T_{\text{ref}}}{T(x)}},\quad\Sigma_{\text{f}}\big(T(x)\big)=\Sigma_{\text{f,ref}}\sqrt{\frac{T_{\text{ref}}}{T(x)}}.\label{eq:couplingmechanism}
\end{equation}

\subsection{Deterministic weak formulation}
Let~$H=H^{1}(]0,L[)$ be the space of functions that are sufficiently regular to describe the solutions.
The weak formulation then involves finding~$T$ and~$\Phi$ in~$H$ such that
\begin{equation}
\label{eq:neutron2}\begin{aligned}
&\int_{0}^{L}k\frac{dT}{dx}\frac{dS}{dx}dx+\int_{0}^{L}h(T-T_{\infty})Sdx=\int_{0}^{L}E_{\text{f}}\Sigma_{\text{f}}(T)\Phi Sdx,&&\forall S\in H,\\
&\int_{0}^{L}D(T)\frac{d\Phi}{dx}\frac{d\Psi}{dx}dx+\int_{0}^{L}\Big(\Sigma_{\text{a}}(T)-\nu\Sigma_{\text{f}}(T)\Big)\Phi\Psi dx=\int_{0}^{L}s\Psi dx,&&\forall \Psi\in H.\\
\end{aligned}
\end{equation}

\subsection{Random thermal transmittivity}
Uncertainties are incorporated by modeling the thermal transmittivity as a random field~$\{h(x,\cdot),1\leq x\leq L\}$ such that
\begin{equation}
h(x,\boldsymbol{\xi})=\overline{h}\bigg(1+\delta\sum_{j=1}^{m}\sqrt{\lambda_{j}}\sqrt{3}\xi_{j}\phi^{j}(x)\bigg),\label{eq:hN}
\end{equation}
where the~$\xi_{j}$ are statistically independent uniform random variables defined on a probability triple~$(\Theta,\mathcal{T},P)$ with values in $[-1,1]$ and the~$\sqrt{3}\xi_{j}$ are thus uniform random variables with unit standard deviation; further, the~$\lambda_{j}$ and~$\phi^{j}$ are the eigenvalues and eigenmodes, respectively, of the eigenproblem~$\mathcal{C}(\phi^{j})=\lambda_{j}\phi^{j}$, where~$\mathcal{C}$ is the covariance integral operator with
\begin{equation}
C(x,y)=\frac{4a^{2}}{\pi^{2}(x-y)^{2}}\sin^{2}\left(\frac{\pi(x-y)}{2a}\right)\label{eq:covar}
\end{equation}
as the kernel; here, the parameter~$a$ is the spatial correlation length of~$\{h(x,\cdot),1\leq x\leq L\}$.
Clearly, the random field~$\{h(x,\cdot),1\leq x\leq L\}$ thus obtained is such that the random variable~$h(x,\cdot)$ has the mean~$\overline{h}$ and coefficient of variation~$\delta$ at every position~$x$, at least when the approximation error introduced owing to the truncation of the expansion after $m$ terms is not taken into account.  

\subsection{Stochastic weak formulation} 
The weak formulation of the stochastic problem involves finding random variables~$T$ and~$\Phi$ defined on~$(\Theta,\mathcal{T},P)$, with values in~$H$, such that
\begin{equation}
\label{eq:neutron2b}
\begin{aligned}
&\int_{0}^{L}k\frac{dT}{dx}\frac{dS}{dx}dx+\int_{0}^{L}h(\boldsymbol{\xi})(T-T_{\infty})Sdx=\int_{0}^{L}E_{\text{f}}\Sigma_{\text{f}}(T)\Phi Sdx,&&\forall S\in H,\\
&\int_{0}^{L}D(T)\frac{d\Phi}{dx}\frac{d\Psi}{dx}dx+\int_{0}^{L}\Big(\Sigma_{\text{a}}(T)-\nu\Sigma_{\text{f}}(T)\Big)\Phi\Psi dx=\int_{0}^{L}s\Psi dx,&&\forall \Psi\in H.
\end{aligned}
\end{equation}
 
\subsection{Discretization of space}
The finite element~(FE) method is used for the discretization of space.
The domain~$[0,L]$ is meshed using $r-1$ elements of equal length.  
Let~$N_{1},\ldots,N_{r}$ then be a basis of element-wise linear shape functions such that~$N_{j}$ takes value~$1$ at the~$j$-th node and~0 at other nodes. 
Using this basis, the random temperature~$T$ and neutron flux~$\Phi$ are approximated as follows:
\begin{equation}
\label{eq:Vh1}
\begin{aligned}
T^{r}(x)=\sum_{j=1}^{r}T_{j}N_{j}(x),\quad\quad\quad T_{j}\in\real,\\
\Phi^{r}(x)=\sum_{j=1}^{r}\Phi_{j}N_{j}(x),\quad\quad\quad \Phi_{j}\in\real.\\
\end{aligned}
\end{equation} 
The FE discretization of the stochastic weak formulation~(\ref{eq:neutron2b}) then involves finding random variables~$\boldsymbol{T}=(T_{1},\ldots,T_{r})$ and~$\boldsymbol{\Phi}=(\Phi_{1},\ldots,\Phi_{r})$ defined on~$(\Theta,\mathcal{T},P)$, with values in~$\real^{r}$, which collect the nodal values of the random temperature and neutron flux such that
\begin{equation}
\label{eq:discrfE3b}
\begin{aligned}
&[\boldsymbol{K}+\boldsymbol{H}(\boldsymbol{\xi})]\boldsymbol{T}=\boldsymbol{q}(\boldsymbol{\Phi},\boldsymbol{T}),\\
&[\boldsymbol{D}(\boldsymbol{T})+\boldsymbol{M}(\boldsymbol{T})]\boldsymbol{\Phi}=\boldsymbol{s}.
\end{aligned}
\end{equation}
Here, $\boldsymbol{K}$, $\boldsymbol{H}$, $\boldsymbol{D}(\boldsymbol{T})$, and $\boldsymbol{M}(\boldsymbol{T})$ are $r$-dimensional matrices, and $\boldsymbol{q}(\boldsymbol{\Phi},\boldsymbol{T})$ and $\boldsymbol{s}$ are $r$-dimensional vectors such that
\begin{align}
\boldsymbol{S}_{1}^{\mathrm{T}}\boldsymbol{K}\boldsymbol{S}_{2}&=\int_{0}^{L}k\frac{dS^{r}_{1}}{dx}\frac{dS^{r}_{2}}{dx}dx,\\
\boldsymbol{S}_{1}^{\mathrm{T}}\boldsymbol{H}\boldsymbol{S}_{2}&=\int_{0}^{L}hS^{r}_{1}S^{r}_{2}dx,\\
\boldsymbol{\Psi}_{1}^{\mathrm{T}}\boldsymbol{D}(\boldsymbol{T})\boldsymbol{\Psi}_{2}&=\int_{0}^{L}D\big(T^{r}\big)\frac{d\Psi^{r}_{1}}{dx}\frac{d\Psi^{r}_{2}}{dx}dx,\\
\boldsymbol{\Psi}_{1}^{\mathrm{T}}\boldsymbol{M}(\boldsymbol{T})\boldsymbol{\Psi}_{2}&=\int_{0}^{L}\Big(\Sigma_{\text{a}}\big(T^{r}\big)-\nu\Sigma_{\text{f}}\big(T^{r}\big)\Big)\Psi^{r}_{1}\Psi^{r}_{2}dx,\\
\boldsymbol{S}^{\mathrm{T}}\boldsymbol{q}(\boldsymbol{T},\boldsymbol{\Phi})&=\int_{0}^{L}E_{\text{f}}\Sigma_{\text{f}}\big(T^{r}\big)\Phi^{r}S^{r}dx+\int_{0}^{L}hT_{\infty}S^{r}dx,\\
\boldsymbol{S}^{\mathrm{T}}\boldsymbol{s}&=\int_{0}^{L}s\Psi^{r}dx.
\end{align}

\subsection{Reformulation as a realization of the model problem}
The aforementioned illustration problem can be reformulated as a particular realization of the general model problem introduced in Sec.~\ref{sec:sec4} as follows:
\begin{equation}
\begin{aligned}
&\boldsymbol{T}=\boldsymbol{a}(\boldsymbol{T},\boldsymbol{\Phi},\boldsymbol{\xi}),&&\quad\quad\quad\boldsymbol{a}:\real^{r}\times\real^{r}\times\real^{m}\rightarrow\real^{r},\\
&\boldsymbol{\Phi}=\boldsymbol{b}(\boldsymbol{T}),&&\quad\quad\quad\boldsymbol{b}:\real^{r}\rightarrow\real^{r},
\end{aligned}
\end{equation}
where $\boldsymbol{a}(\boldsymbol{T},\boldsymbol{\Phi},\boldsymbol{\xi})=[\boldsymbol{K}+\boldsymbol{H}(\boldsymbol{\xi})]^{-1}\boldsymbol{q}(\boldsymbol{\Phi},\boldsymbol{T})$ and $\boldsymbol{b}(\boldsymbol{T})=[\boldsymbol{D}(\boldsymbol{T})+\boldsymbol{M}(\boldsymbol{T})]^{-1}\boldsymbol{s}$.
This reformulation indicates that the illustration problem is a simplified realization of the model problem, for three reasons.
First, the data of the neutronics subproblem are not affected by their own sources of uncertainty~$\boldsymbol{\zeta}$.
Second, the neutronics subproblem admits a direct solution that does not require iteration.
Lastly, the heat and neutronics subproblems are coupled directly through their solution variables rather than through intermediate coupling variables.

\subsection{Discretization of the random dimension}
Because the sources of uncertainty~$\xi_{1},\ldots,\xi_{m}$ are statistically independent uniform random variables with values in $[-1,1]$, the PC basis~$\{\psi_{\boldsymbol{\alpha}},\;\boldsymbol{\alpha}\in\integer^{m}\}$ consists here of the normalized Legendre polynomials in $m$ variables of increasing total degree with $\psi_{\boldsymbol{0}}=1$.

The PC basis provides approximate representations of the successive approximations determined by the iterative method that does not dimension reduction as follows: 
\begin{equation}
\begin{aligned}
\boldsymbol{T}^{\ell,p}=\sum_{|\boldsymbol{\alpha}|=0}^{p}\boldsymbol{T}_{\boldsymbol{\alpha}}^{\ell}\psi_{\boldsymbol{\alpha}}(\boldsymbol{\xi}),&&\quad\quad\quad\boldsymbol{T}_{\boldsymbol{\alpha}}^{\ell}\in\real^{r},\\
\boldsymbol{\Phi}^{\ell,p}=\sum_{|\boldsymbol{\alpha}|=0}^{p}\boldsymbol{\Phi}_{\boldsymbol{\alpha}}^{\ell}\psi_{\boldsymbol{\alpha}}(\boldsymbol{\xi}),&&\quad\quad\quad\boldsymbol{\Phi}_{\boldsymbol{\alpha}}^{\ell}\in\real^{r}.\\
\end{aligned}
\end{equation}

\subsection{Dimension reduction by KL decomposition of the temperature} \label{sec:sec56}
Now, we will demonstrate the proposed methodology by approximating the random temperature by a truncated KL decomposition as it is communicated from the heat to the neutronics subproblem.

The PC basis provides approximate representations of the successive approximations determined by the iterative method involving dimension reduction as follows: 
\begin{equation}
\begin{aligned}
\label{eq:upceklb}\widehat{\boldsymbol{T}}{}^{\ell,p}&=\sum_{|\boldsymbol{\alpha}|=0}^{p}\widehat{\boldsymbol{T}}{}^{\ell}_{\boldsymbol{\alpha}}\psi_{\boldsymbol{\alpha}}(\boldsymbol{\xi}),&&\quad\quad\quad\widehat{\boldsymbol{T}}{}^{\ell}_{\boldsymbol{\alpha}}\in\real^{r},\\
\widehat{\boldsymbol{\Phi}}{}^{\ell,p}&=\sum_{|\boldsymbol{\alpha}|=0}^{p}\widehat{\boldsymbol{\Phi}}{}_{\boldsymbol{\alpha}}^{\ell}\psi_{\boldsymbol{\alpha}}(\boldsymbol{\xi}),&&\quad\quad\quad\widehat{\boldsymbol{\Phi}}_{\boldsymbol{\alpha}}^{\ell}\in\real^{r}.\\
\end{aligned}
\end{equation}
The mean and covariance of~$\widehat{\boldsymbol{T}}{}^{\ell,p}$ are given by:
\begin{align}
\overline{\boldsymbol{T}}{}^{\ell}&=\widehat{\boldsymbol{T}}{}^{\ell}_{\boldsymbol{0}},\\
\boldsymbol{C}_{\boldsymbol{T}}^{\ell}&=\sum_{|\boldsymbol{\alpha}|=1}^{p}\widehat{\boldsymbol{T}}{}^{\ell}_{\boldsymbol{\alpha}}\big(\widehat{\boldsymbol{T}}{}_{\boldsymbol{\alpha}}^{\ell}\big)^{\mathrm{T}}.
\end{align}   
Further, let the $r$-dimensional square matrix~$\boldsymbol{W}$ be the Gram matrix of the FE basis, i.e.,
\begin{equation}
\boldsymbol{W}=\begin{bmatrix}
\langle N_{1},N_{1}\rangle_{H} & \ldots & \langle N_{1},N_{r}\rangle_{H}\\
\vdots & & \vdots\\
\langle N_{r},N_{1}\rangle_{H} & \ldots & \langle N_{r},N_{r}\rangle_{H}
\end{bmatrix},
\end{equation}
where the inner product~$\langle\cdot,\cdot\rangle_{H}$ is such that~$\langle S_{1},S_{2}\rangle_{H}=\int_{0}^{L}S_{1}S_{2}dx+\int_{0}^{L}(dS_{1}/dx)(dS_{2}/dx)dx$ for any pair~$S_{1}$ and~$S_{2}$ of functions in~$H$. 
The solution of the generalized eigenproblem $\boldsymbol{W}^{\mathrm{T}}\boldsymbol{C}_{\boldsymbol{T}}^{\ell}\boldsymbol{W}\boldsymbol{\phi}^{j,\ell}=\lambda_{j}^{\ell}\boldsymbol{W}\boldsymbol{\phi}^{j,\ell}$ then provides the eigenvalues~$\lambda_{j}^{\ell}$ and the associated eigenmodes~$\boldsymbol{\phi}^{j,\ell}$ required to construct a reduced-dimensional representation~$\widehat{\boldsymbol{T}}{}^{\ell,p,d}$ of $\widehat{\boldsymbol{T}}{}^{\ell,p}$ as follows:
\begin{equation}
\widehat{\boldsymbol{T}}{}^{\ell,p,d}=\overline{\boldsymbol{T}}{}^{\ell}+\sum_{j=1}^{d}\sqrt{\lambda_{j}^{\ell}}\eta{}_{j}^{\ell,p}\boldsymbol{\phi}^{j,\ell},\label{eq:Tpdkl}
\end{equation}
where the~$\eta_{j}^{\ell,p}$ are random variables defined on~$(\Theta,\mathcal{T},P)$, with values in~$\real$, such that
\begin{equation}
\eta{}_{j}^{\ell,p}=\frac{1}{\sqrt{\lambda_{j}^{\ell}}}\big(\widehat{\boldsymbol{T}}{}^{\ell,p}-\overline{\boldsymbol{T}}{}^{\ell}\big)^{\mathrm{T}}\boldsymbol{W}\boldsymbol{\phi}^{j,\ell}.\label{eq:upcekletab}
\end{equation}
By substituting~(\ref{eq:upceklb}) in~(\ref{eq:upcekletab}), a representation of the~$\eta^{\ell,p}_{i}$ as a PC expansion is obtained: 
\begin{equation}
\boldsymbol{\eta}^{\ell,p}=\sum_{|\boldsymbol{\alpha}|=1}^{p}\boldsymbol{\eta}^{\ell}_{\boldsymbol{\alpha}}\psi_{\boldsymbol{\alpha}}(\boldsymbol{\xi})\quad\text{with}\quad\boldsymbol{\eta}_{\boldsymbol{\alpha}}^{\ell}=\begin{bmatrix}\frac{1}{\sqrt{\lambda_{1}^{\ell}}}\big(\widehat{\boldsymbol{T}}{}_{\boldsymbol{\alpha}}^{\ell}\big)^{\mathrm{T}}\boldsymbol{W}\boldsymbol{\phi}^{1,\ell} & \ldots & \frac{1}{\sqrt{\lambda_{d}^{\ell}}}\big(\widehat{\boldsymbol{T}}{}_{\boldsymbol{\alpha}}^{\ell}\big)^{\mathrm{T}}\boldsymbol{W}\boldsymbol{\phi}^{d,\ell}\end{bmatrix}^{\mathrm{T}},\label{eq:upcervetaaaa}
\end{equation}
thus completely characterizing the reduced random variables as a PC expansion.

It should be noted that the random neutron flux, in principle, could also be reduced as it passes from the neutronics subproblem to the heat subproblem.
However, because the data of the neutronics subproblem are not affected by their own sources of uncertainty, a reduction of the random neutron flux would not lower the number of sources of uncertainty that enter the heat subproblem, and thus would not pave the way for a solution of the heat subproblem in a reduced-dimensional space.
This extension is therefore not demonstrated.

\subsection{Selection of the reduced dimension}
At each iteration, we select  the number of terms retained in~(\ref{eq:Tpdkl}) by the KL decomposition~$\widehat{\boldsymbol{T}}{}^{\ell,p,d}$ of~$\widehat{\boldsymbol{T}}{}^{\ell,p}$ as the smallest dimension $d$ that satisfies the following condition:
\begin{equation}
\int_{\Theta}\vectornorm{\widehat{\boldsymbol{T}}{}^{\ell,p}-\widehat{\boldsymbol{T}}{}^{\ell,p,d}}^{2}_{\boldsymbol{W}}dP\leq tol,\quad\quad\quad\forall\ell\in\integer,\label{eq:criterion}
\end{equation}
where~$tol$ is a prescribed tolerance level.
Clearly, this criterion may result in the dependence of the reduced dimension~$d$ on the iteration~$\ell$.

\subsection{Concluding remarks}
Algorithm~\ref{algo:algo6} outlines an implementation wherein the nonintrusive projection method is adopted for the discretization of the random dimension.
This algorithm requires a quadrature rule for integration with respect to the multidimensional uniform probability distribution of the input random variables. 
In this work, we use a sparse-grid Gauss-Legendre quadrature rule~\citep{holtz2010}.

\begin{algorithm}[htp]
\SetKwInOut{Input}{Input}\SetKwInOut{Output}{Output}
\SetKw{KwAnd}{and}
\SetKwBlock{FirstMonoProblem}{neutronics subproblem}{end}
\SetKwBlock{SecondMonoProblem}{heat subproblem}{end}
\SetKwBlock{NexusRegion}{dimension reduction}{end}
\Input{KL truncation tolerance $tol$;\\
PC basis $\small{\big\{\psi_{\boldsymbol{\alpha}},\;0\leq|\boldsymbol{\alpha}|\leq p\big\}}\normalsize$ up to total degree $\small{p}\normalsize$ w.r.t. $\small{P_{\boldsymbol{\xi}}}\normalsize$;\\
Quadrature rule $\small{\{(\boldsymbol{\xi}_{k},w_{k}),\;1\leq k\leq\nu\}}\normalsize$ of level $p+1$ w.r.t. $\small{P_{\boldsymbol{\xi}}}\normalsize$\;}
$\small{\ell=1}\normalsize$\;
\Repeat{$($convergence$)$}{
\SecondMonoProblem{
\For{$k=1$ \KwTo $\,\nu$}{
Solve $\small{\Big[\boldsymbol{K}+\boldsymbol{H}\big(\boldsymbol{\xi}_{k}\big)\Big]\widehat{\boldsymbol{T}}{}^{\ell}\big(\boldsymbol{\xi}_{k}\big)=\boldsymbol{q}\Big(\widehat{\boldsymbol{T}}{}^{\ell-1,p}(\boldsymbol{\xi}_{k}),\widehat{\boldsymbol{\Phi}}{}^{\ell-1,p}(\boldsymbol{\xi}_{k})\Big);}\normalsize$\vspace{0mm}}
Compute PC coordinates of $\small{\widehat{\boldsymbol{T}}{}^{\ell,p}}\normalsize$ by nonintrusive projection:\\
\vspace{-1mm}
\small{
\begin{equation*}
\widehat{\boldsymbol{T}}{}^{\ell}_{\boldsymbol{\alpha}}=\sum_{k=1}^{\nu}\widehat{\boldsymbol{T}}{}^{\ell}\big(\boldsymbol{\xi}_{k}\big)\psi_{\boldsymbol{\alpha}}\big(\boldsymbol{\xi}_{k}\big)w_{k};
\end{equation*}}\normalsize\\
\vspace{-3mm}
}
\NexusRegion{  
Compute mean $\small{\overline{\boldsymbol{T}}{}^{\ell}=\widehat{\boldsymbol{T}}{}_{\boldsymbol{0}}^{\ell}}\normalsize$ and covariance $\small{\boldsymbol{C}{}_{\widehat{\boldsymbol{T}}}^{\ell}=\sum_{|\boldsymbol{\alpha}|=1}^{p}\widehat{\boldsymbol{T}}{}_{\boldsymbol{\alpha}}^{\ell}(\widehat{\boldsymbol{T}}{}_{\boldsymbol{\alpha}}^{\ell})^{\mathrm{T}}}\normalsize$\;
Solve eigenproblem $\small{\boldsymbol{W}^{\mathrm{T}}\boldsymbol{C}_{\widehat{\boldsymbol{T}}}^{\ell}\boldsymbol{W}\boldsymbol{\phi}^{j,\ell}=\lambda_{j}^{\ell}\boldsymbol{W}\boldsymbol{\phi}^{j,\ell}}\normalsize$\;
Choose $\small{d}\normalsize$ such that $\small{\sum_{|\boldsymbol{\alpha}|=1}^{p}(\widehat{\boldsymbol{T}}{}_{\boldsymbol{\alpha}}^{\ell})^{\mathrm{T}}\boldsymbol{W}\widehat{\boldsymbol{T}}{}_{\boldsymbol{\alpha}}^{\ell}-\sum_{j=1}^{d}\lambda_{j}^{\ell}\leq tol}\normalsize$\;
Compute PC coordinates of $\small{\eta_{j}^{\ell,p}}\normalsize$ by $\small{\eta_{j,\boldsymbol{\alpha}}^{\ell}=(\widehat{\boldsymbol{T}}{}_{\boldsymbol{\alpha}}^{\ell})^{\mathrm{T}}\boldsymbol{W}\boldsymbol{\phi}^{j,\ell}}\normalsize$ for $\small{j=1}\normalsize$ to $\small{d}\normalsize$\;
}
\FirstMonoProblem{
\For{$k=1$ \KwTo $\,\nu$}{
Solve $\small{\Big[\boldsymbol{D}\Big(\widehat{\boldsymbol{T}}{}^{\ell,p,d}(\boldsymbol{\xi}_{k})\Big)+\boldsymbol{M}\Big(\widehat{\boldsymbol{T}}{}^{\ell,p,d}(\boldsymbol{\xi}_{k})\Big)\Big]\widehat{\boldsymbol{\Phi}}{}^{\ell}(\boldsymbol{\xi}_{k})=\boldsymbol{s},}\normalsize$\\
with $\small{\widehat{\boldsymbol{T}}{}^{\ell,p,d}(\boldsymbol{\xi}_{k})=\overline{\boldsymbol{T}}{}^{\ell}+\sum_{j=1}^{d}\sqrt{\lambda_{j}^{\ell}}\bigg(\sum_{|\boldsymbol{\alpha}|=1}^{p}\eta_{j,\boldsymbol{\alpha}}^{\ell}\psi_{\boldsymbol{\alpha}}(\boldsymbol{\xi}_{k})\bigg)\boldsymbol{\phi}^{j,\ell}}\normalsize$\;}
Compute PC coordinates of $\small{\widehat{\boldsymbol{\Phi}}{}^{\ell,p}}\normalsize$ by nonintrusive projection:\\
\vspace{-1mm}
\small{
\begin{equation*}
\widehat{\boldsymbol{\Phi}}{}^{\ell}_{\boldsymbol{\alpha}}=\sum_{k=1}^{\nu}\widehat{\boldsymbol{\Phi}}{}^{\ell}(\boldsymbol{\xi}_{k})\psi_{\boldsymbol{\alpha}}(\boldsymbol{\xi}_{k})w_{k};
\end{equation*}}\normalsize\\
\vspace{-3mm}
}
$\small{\ell=\ell+1}\normalsize$\;}
\caption{Implementation of the illustration problem.}\label{algo:algo6}
\end{algorithm}

\section{Numerical results}\label{sec:sec6b}
We obtained results using the following properties.  
We assumed the reactor to have a length of~$L=100\,[\text{cm}]$.
Further, we assumed a deterministic and position-independent neutron-diffusion constant~$D_{\text{ref}}=2.2\,[\text{cm}]$; absorption cross section $\Sigma_{\text{a,ref}}=0.0195\,[\text{cm}^{-1}]$; fission cross section~$\Sigma_{\text{a,ref}}=0.0075\,[\text{cm}^{-1}]$; multiplication factor~$\nu=2.2$; neutron source~$s=5.0E11\,[\text{neutrons/s/cm}^{3}$]; ambient temperature~$T_{\infty}=390\,[\text{K}]$; fission energy~$E_{\text{f}}=3.0E\text{-}11\,[\text{J/neutrons}]$; and temperatures~$T_{\text{ref}}=390\,[\text{K}]$, $T_{\text{min}}=390\,[\text{K}]$, and~$T_{\text{max}}=1000\,[\text{K}]$.

\begin{figure}[htp]
  \begin{center}
    \subfigure[Samples.]{\includegraphics[width=0.8\textwidth]{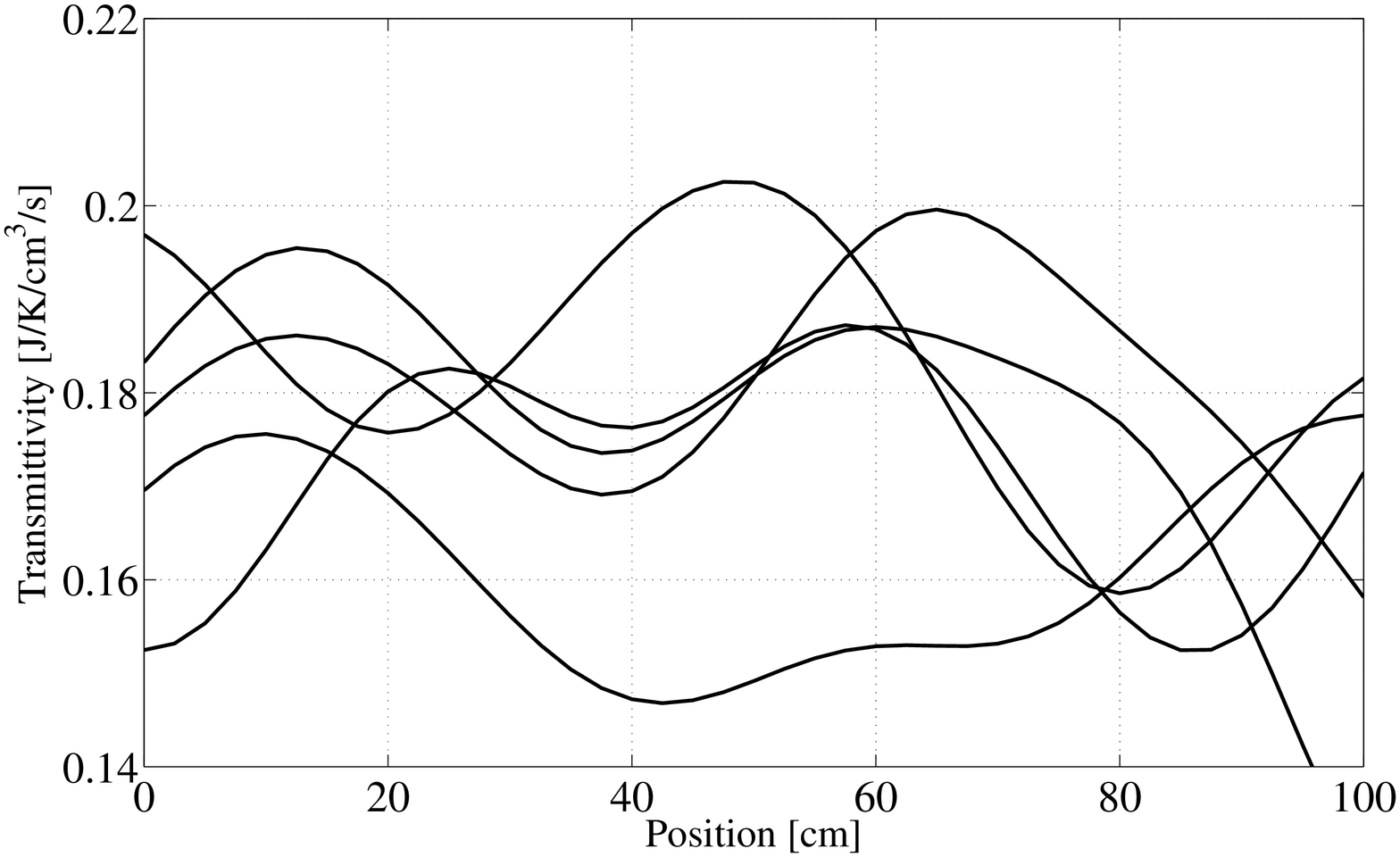}}
    \hfill
    \subfigure[Eigenvalues.]{\includegraphics[width=0.8\textwidth]{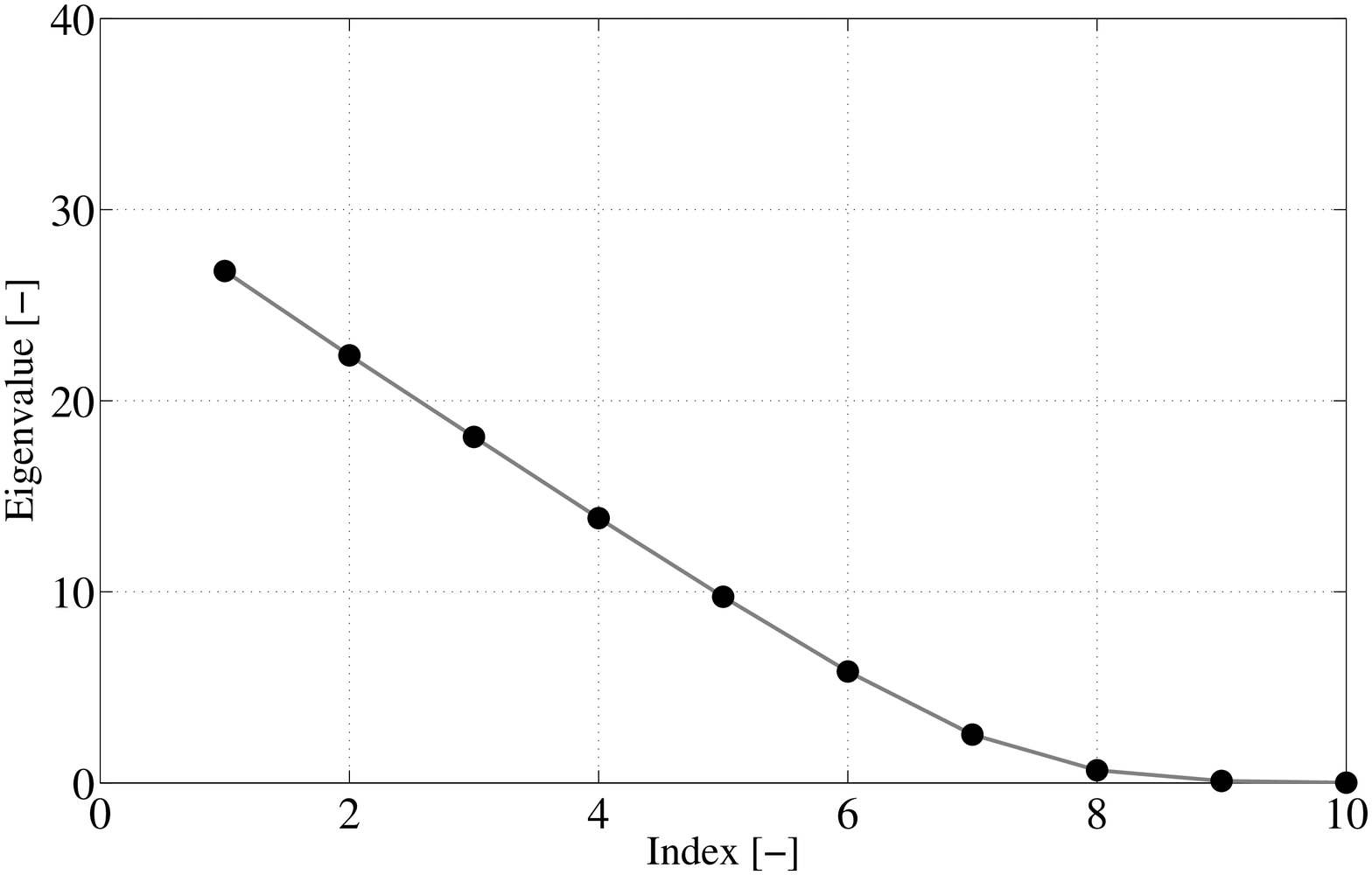}}
    \caption{Thermal transmittivity random field: (a)~five samples and~(b)~ten largest magnitude eigenvalues of the covariance integral operator.}\label{fig:figure1}
  \end{center}
\end{figure}

In addition, we used a thermal transmittivity random field with position-independent mean~$\overline{h}=0.17\,[\text{J/K/cm}^{3}\text{/s}]$, spatial correlation length~$a=15\,[\text{cm}]$, and coefficient of variation~$\delta=10\,\%$.
We retained~$m=10$ terms in expansion~(\ref{eq:hN}).
Figure~\ref{fig:figure1}(a) shows a few sample paths of the random field~$\{h(x,\cdot),\;0\leq x\leq L\}$ thus obtained.
Figure~\ref{fig:figure1}(b) shows the 10 largest magnitude eigenvalues of the covariance integral operator.  

\subsection{Monte Carlo sampling implementation}
\begin{figure}[htp]
  \begin{center}
    \subfigure[Temperature for $k=100\,\text{$[\text{J/K/cm/s}]$}$.]{\includegraphics[width=0.49\textwidth]{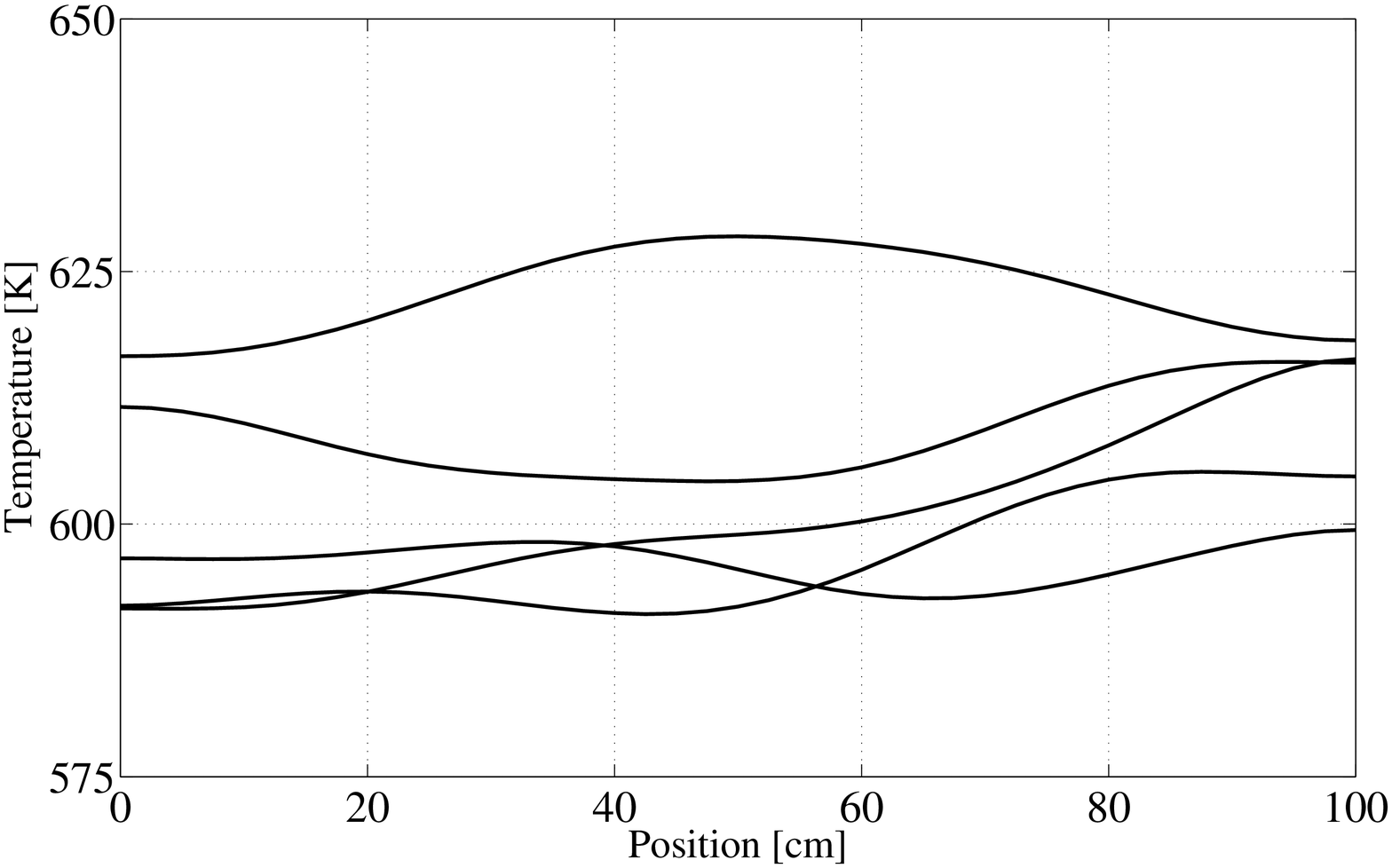}}
    \hfill
    \subfigure[Temperature for $k=1\,\text{$[\text{J/K/cm/s}]$}$.]{\includegraphics[width=0.49\textwidth]{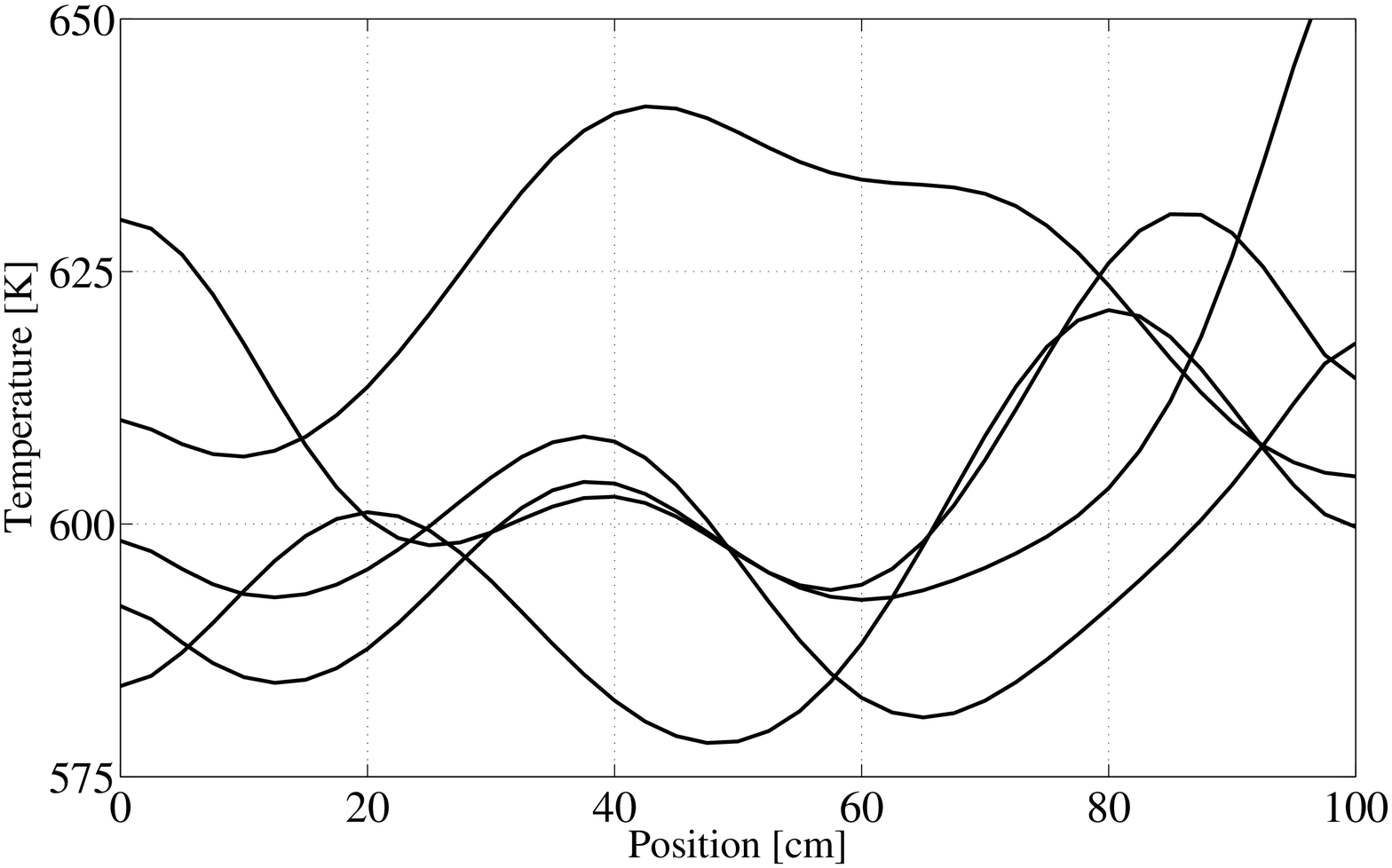}}
    \subfigure[Neutron flux for $k=100\,\text{$[\text{J/K/cm/s}]$}$.]{\includegraphics[width=0.49\textwidth]{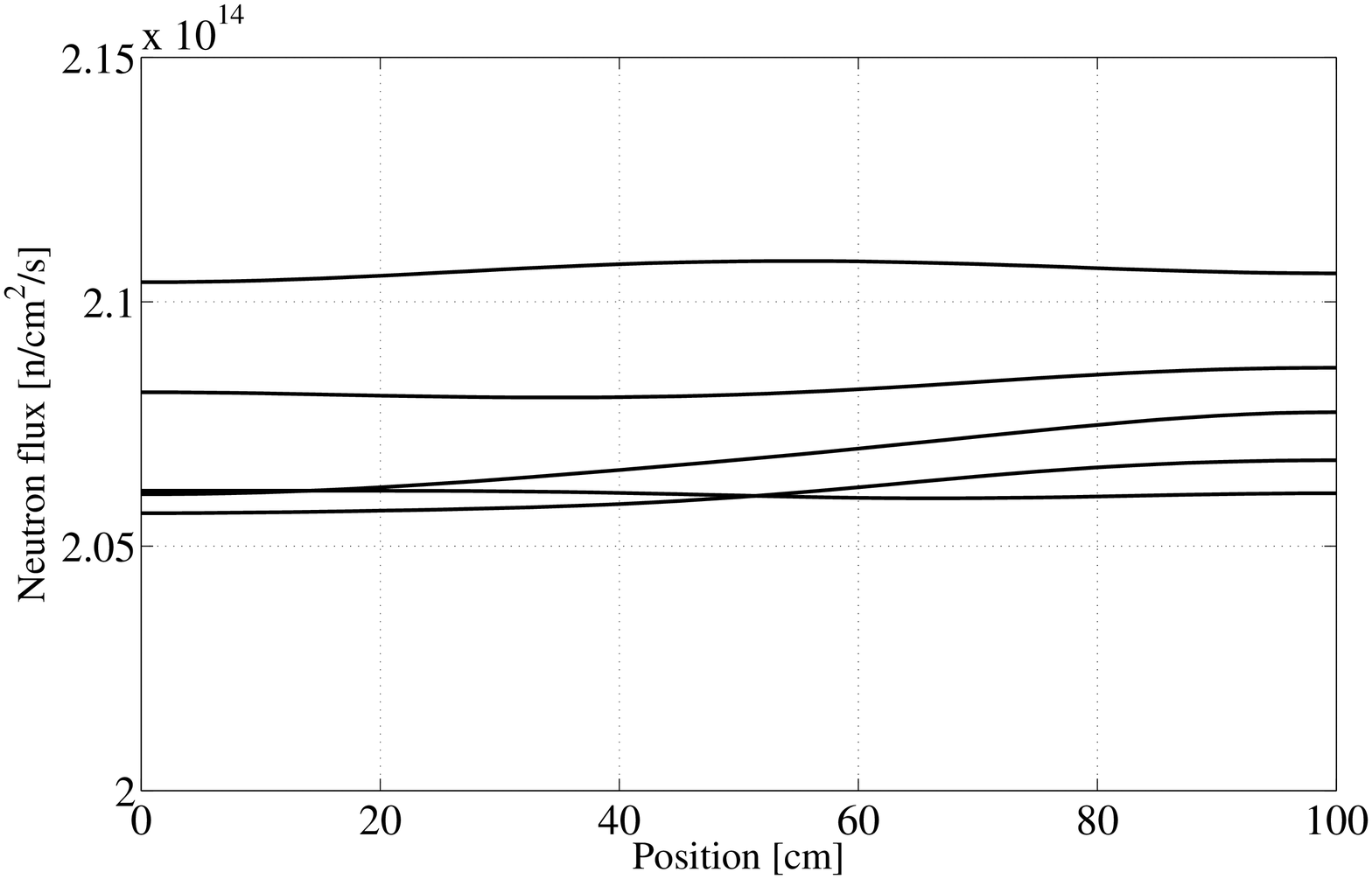}}
    \hfill
    \subfigure[Neutron flux for $k=1\,\text{$[\text{J/K/cm/s}]$}$.]{\includegraphics[width=0.49\textwidth]{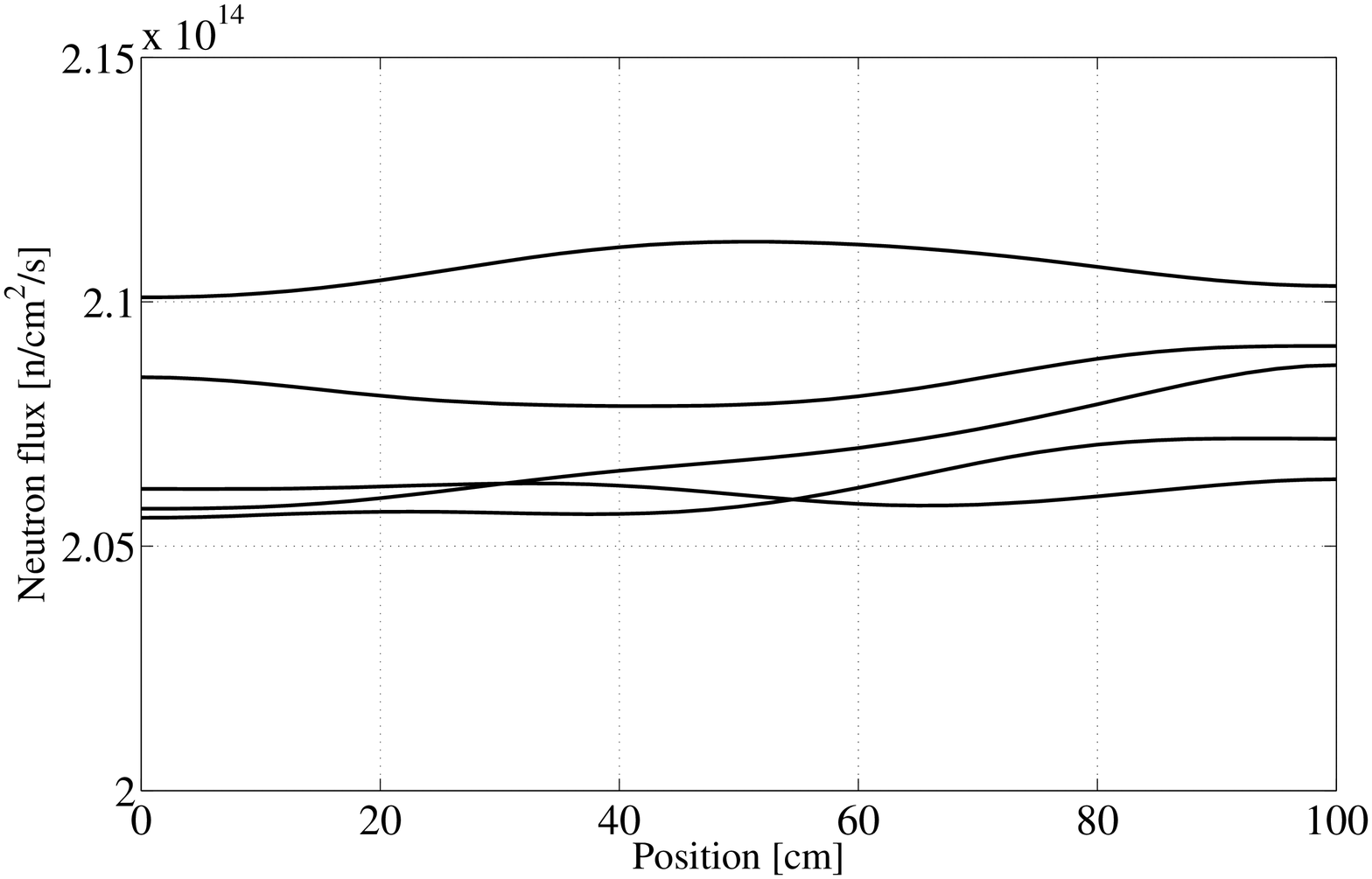}}
    \caption{Monte Carlo simulation: five samples of the solution.}\label{fig:figure2}
  \end{center}
\end{figure}
First, we carried out a Monte Carlo simulation.
We generated $MC=100,000$ sample paths of the thermal transmittivity random field.
Then, for each of these sample paths, we constructed the associated deterministic multiphysics model, each of which we solved using the FE method for the spatial discretization and Gauss-Seidel iteration as the iterative method.
We systematically obtained converged results for~$r-1=40$ finite elements and~$20$ iterations.  

Figure~\ref{fig:figure2} shows a few samples of the random temperature and neutron flux thus obtained for the values of~$k=100\,\text{$[\text{J/K/cm/s}]$}$ and~$k=1\,\text{$[\text{J/K/cm/s}]$}$.
This figure shows that the samples of the random temperature became less smooth as the value of the thermal conductivity was decreased; i.e., the samples exhibited more rapid oscillations with respect to the position in the reactor as the significance of the diffusion of heat was decreased.

\subsection{PC-based implementation not involving dimension reduction}
Next, we implemented a PC-based iterative method that did not involve dimension reduction.
We adopted the nonintrusive stochastic projection method for the discretization of the random dimension.
It should be noted that this implementation amounts to the implementation described in Algorithm~\ref{algo:algo6}, provided that the dimension-reduction step is not carried out, or, equivalently, that the reduced dimension $d$ is chosen equal to $r$ at each iteration.

\begin{figure}[htp]
  \begin{center}
     \subfigure[Convergence for $k=100\,\text{$[\text{J/K/cm/s}]$}$.]{\includegraphics[width=0.49\textwidth]{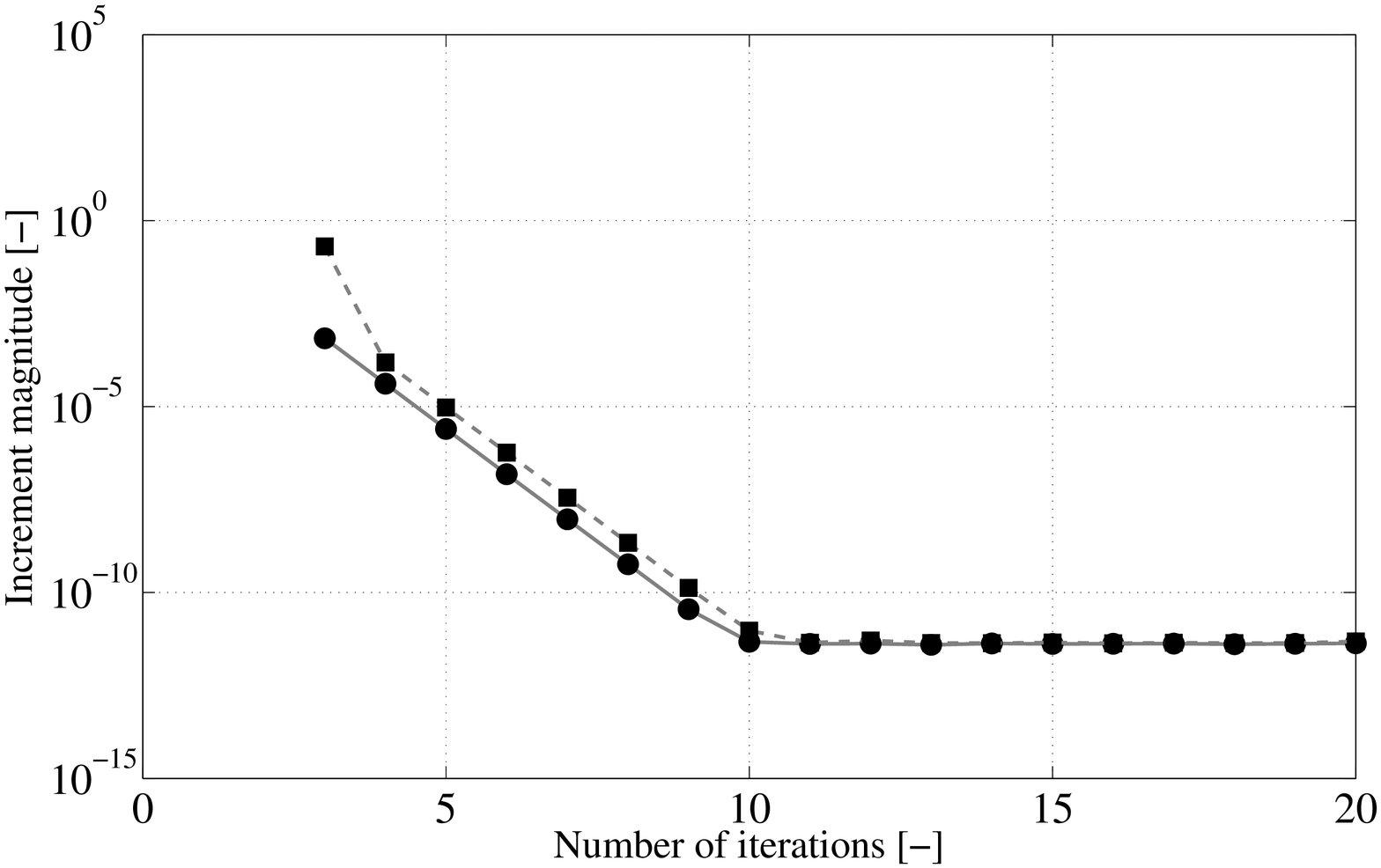}}
     \hfill
     \subfigure[Convergence for $k=1\,\text{$[\text{J/K/cm/s}]$}$.]{\includegraphics[width=0.49\textwidth]{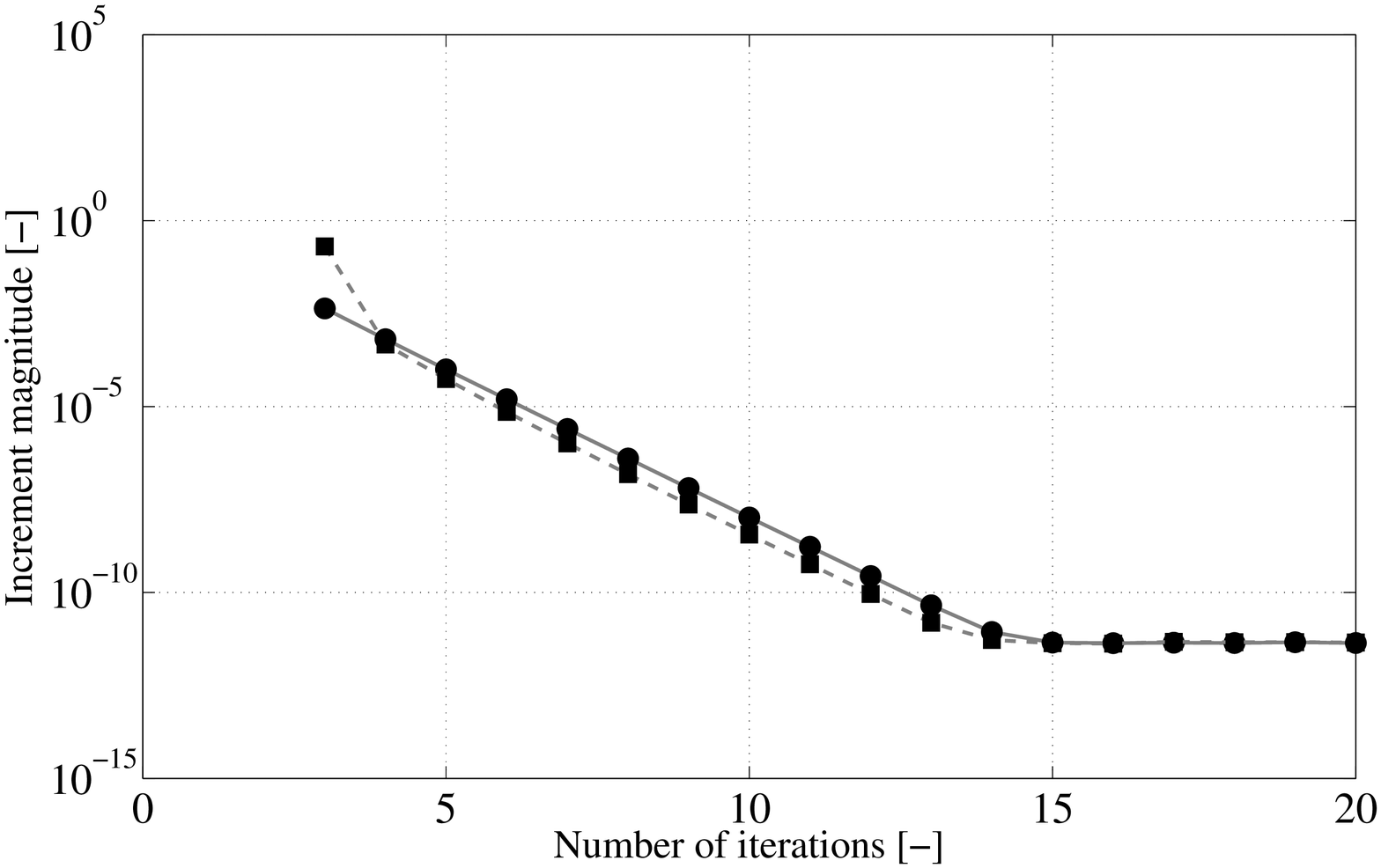}}
     \centerline{$\small{\ell\mapsto\sqrt{\sum_{|\boldsymbol{\alpha}|=0}^{p}\|\boldsymbol{T}_{\boldsymbol{\alpha}}^{\ell}-\boldsymbol{T}_{\boldsymbol{\alpha}}^{\ell-1}\|_{\boldsymbol{W}}^{2}}\Big/\sqrt{\sum_{|\boldsymbol{\alpha}|=0}^{p}\|\boldsymbol{T}_{\boldsymbol{\alpha}}^{\infty}\|_{\boldsymbol{W}}^{2}}}\normalsize$ (circles).}
     \centerline{$\small{\ell\mapsto\sqrt{\sum_{|\boldsymbol{\alpha}|=0}^{p}\|\boldsymbol{\Phi}_{\boldsymbol{\alpha}}^{\ell}-\boldsymbol{\Phi}_{\boldsymbol{\alpha}}^{\ell-1}\|_{\boldsymbol{W}}^{2}}\Big/\sqrt{\sum_{|\boldsymbol{\alpha}|=0}^{p}\|\boldsymbol{\Phi}_{\boldsymbol{\alpha}}^{\infty}\|_{\boldsymbol{W}}^{2}}}\normalsize$ (squares).}
    \caption{Simulation not involving dimension reduction: convergence with respect to the number of iterations.}\label{fig:figure3}
  \end{center}
\end{figure}

\begin{figure}[htp]
  \begin{center}
    \subfigure[Convergence for $k=100\,\text{$[\text{J/K/cm/s}]$}$.]{\includegraphics[width=0.49\textwidth]{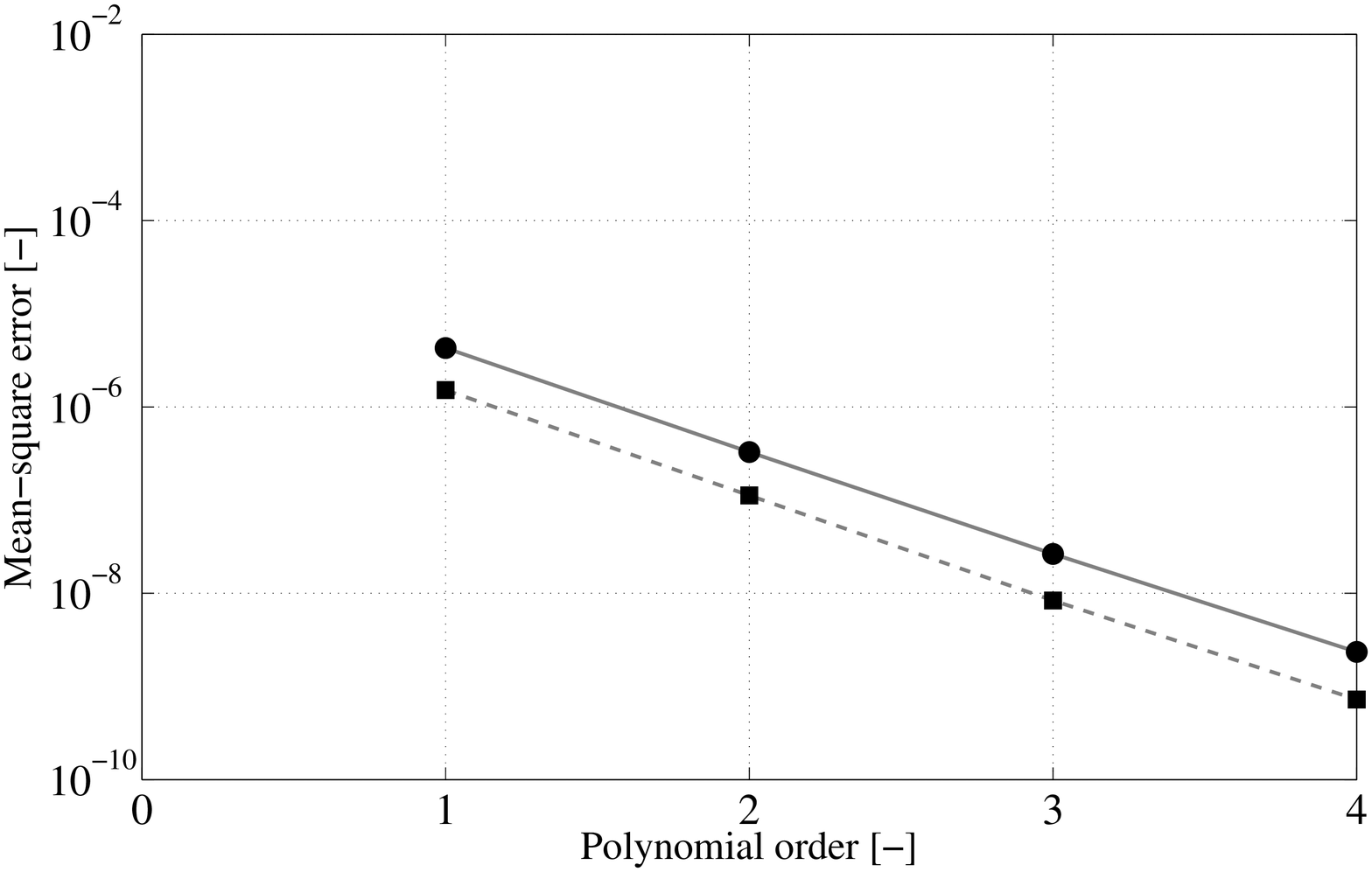}}
    \hfill
    \subfigure[Convergence for $k=1\,\text{$[\text{J/K/cm/s}]$}$.]{\includegraphics[width=0.49\textwidth]{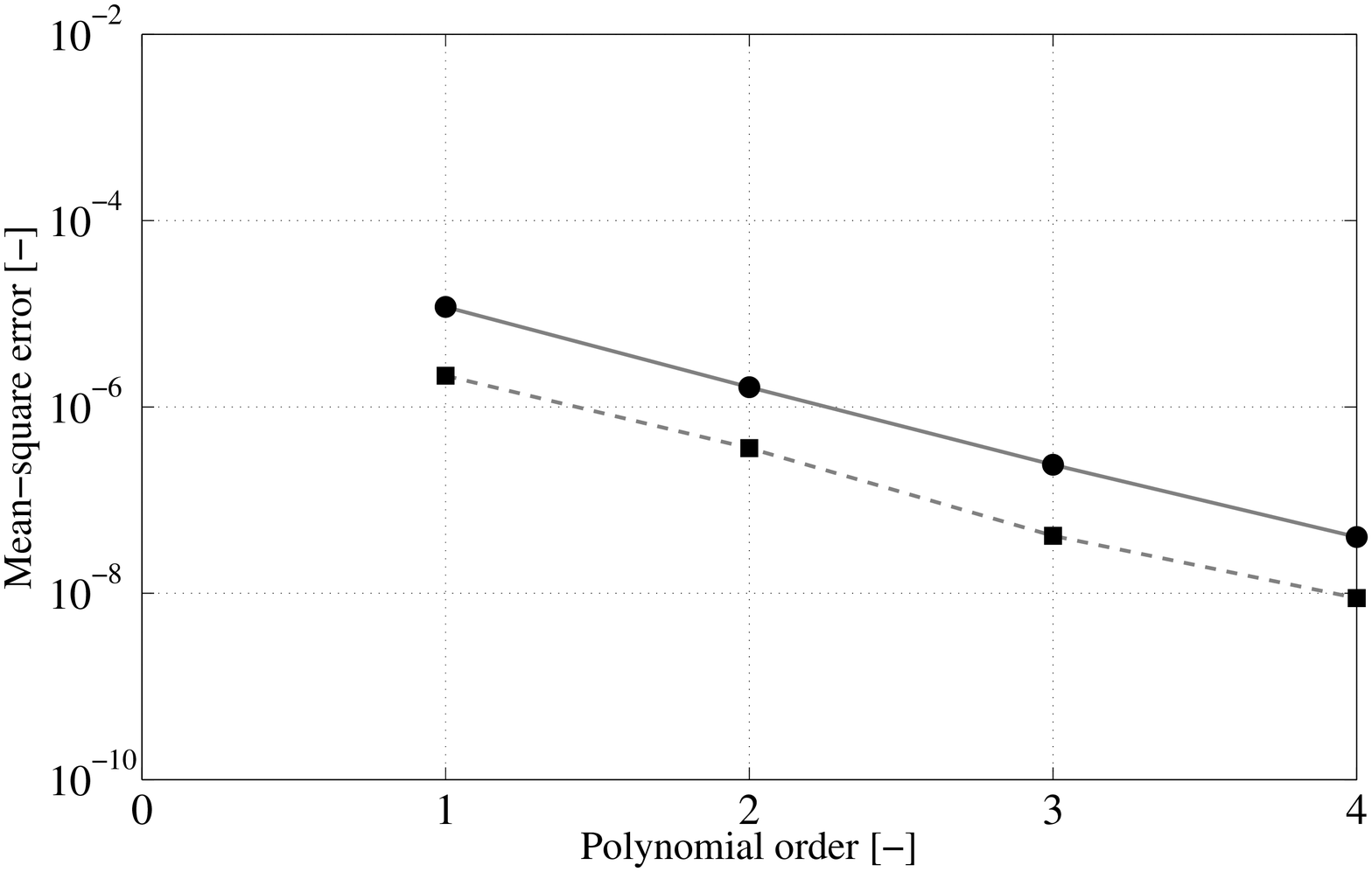}}
    \centerline{$\small{p\mapsto\sqrt{\frac{1}{MC}\sum_{k=1}^{MC}\|\boldsymbol{T}^{\infty}(\boldsymbol{\xi}_{k})-{\boldsymbol{T}}{}^{\infty,p}(\boldsymbol{\xi}_{k})\|_{\boldsymbol{W}}^{2}}\Big/\sqrt{\frac{1}{MC}\sum_{k=1}^{MC}\|\boldsymbol{T}^{\infty}(\boldsymbol{\xi}_{k})\|_{\boldsymbol{W}}^{2}}}\normalsize$ (circles).} 
    \centerline{$\small{p\mapsto\sqrt{\frac{1}{MC}\sum_{k=1}^{MC}\|\boldsymbol{\Phi}^{\infty}(\boldsymbol{\xi}_{k})-{\boldsymbol{\Phi}}{}^{\infty,p}(\boldsymbol{\xi}_{k})\|_{\boldsymbol{W}}^{2}}\Big/\sqrt{\frac{1}{MC}\sum_{k=1}^{MC}\|\boldsymbol{\Phi}^{\infty}(\boldsymbol{\xi}_{k})\|_{\boldsymbol{W}}^{2}}}\normalsize$ (squares).}
    \caption{Simulation not involving dimension reduction: convergence with respect to the polynomial order.}
    \label{fig:figure6}
  \end{center}
\end{figure}

\begin{figure}[htp]
  \begin{center}
    \subfigure[Convergence for $k=100\,\text{$[\text{J/K/cm/s}]$}$.]{\includegraphics[width=0.49\textwidth]{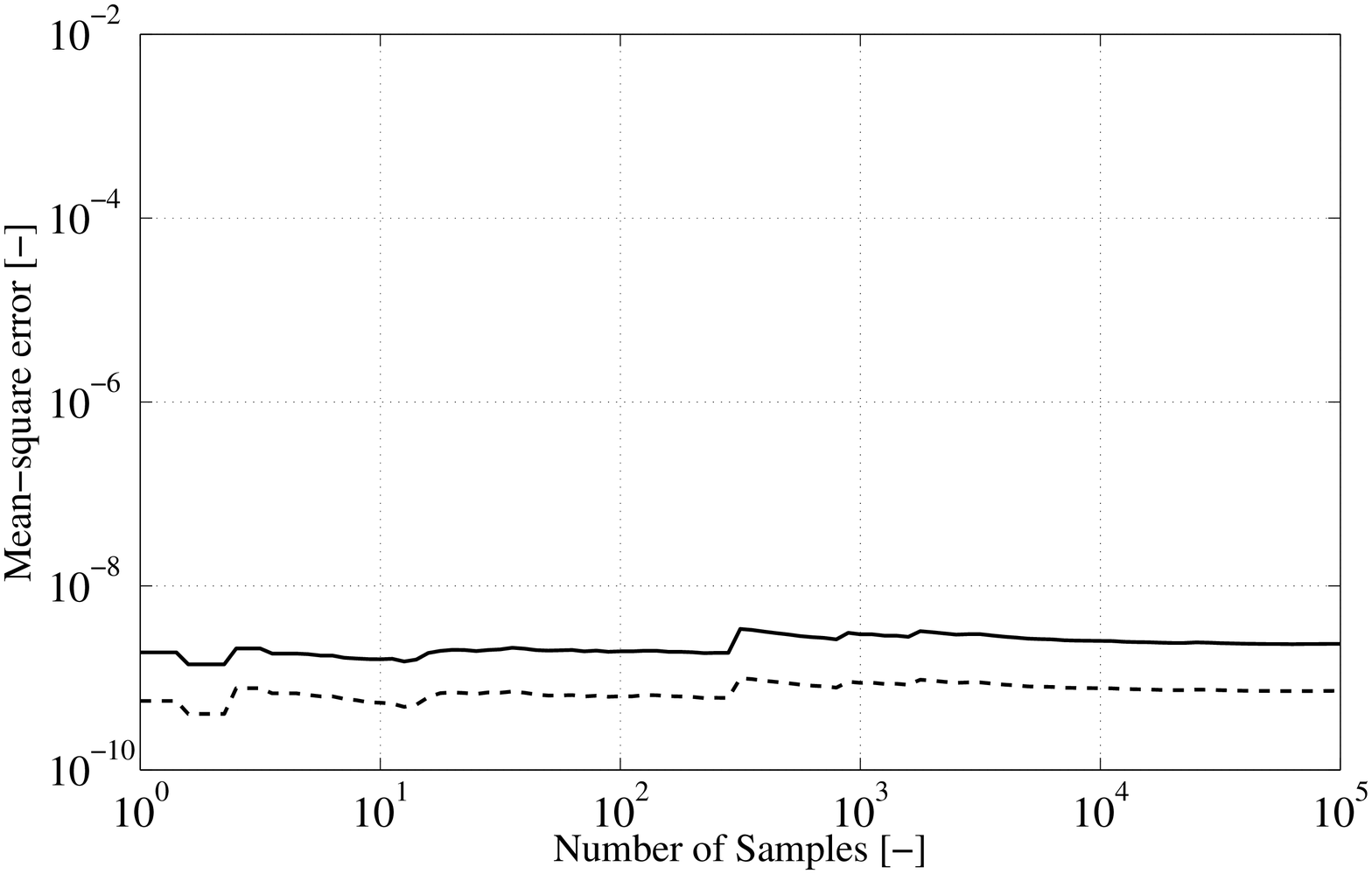}}
    \hfill
    \subfigure[Convergence for $k=1\,\text{$[\text{J/K/cm/s}]$}$.]{\includegraphics[width=0.49\textwidth]{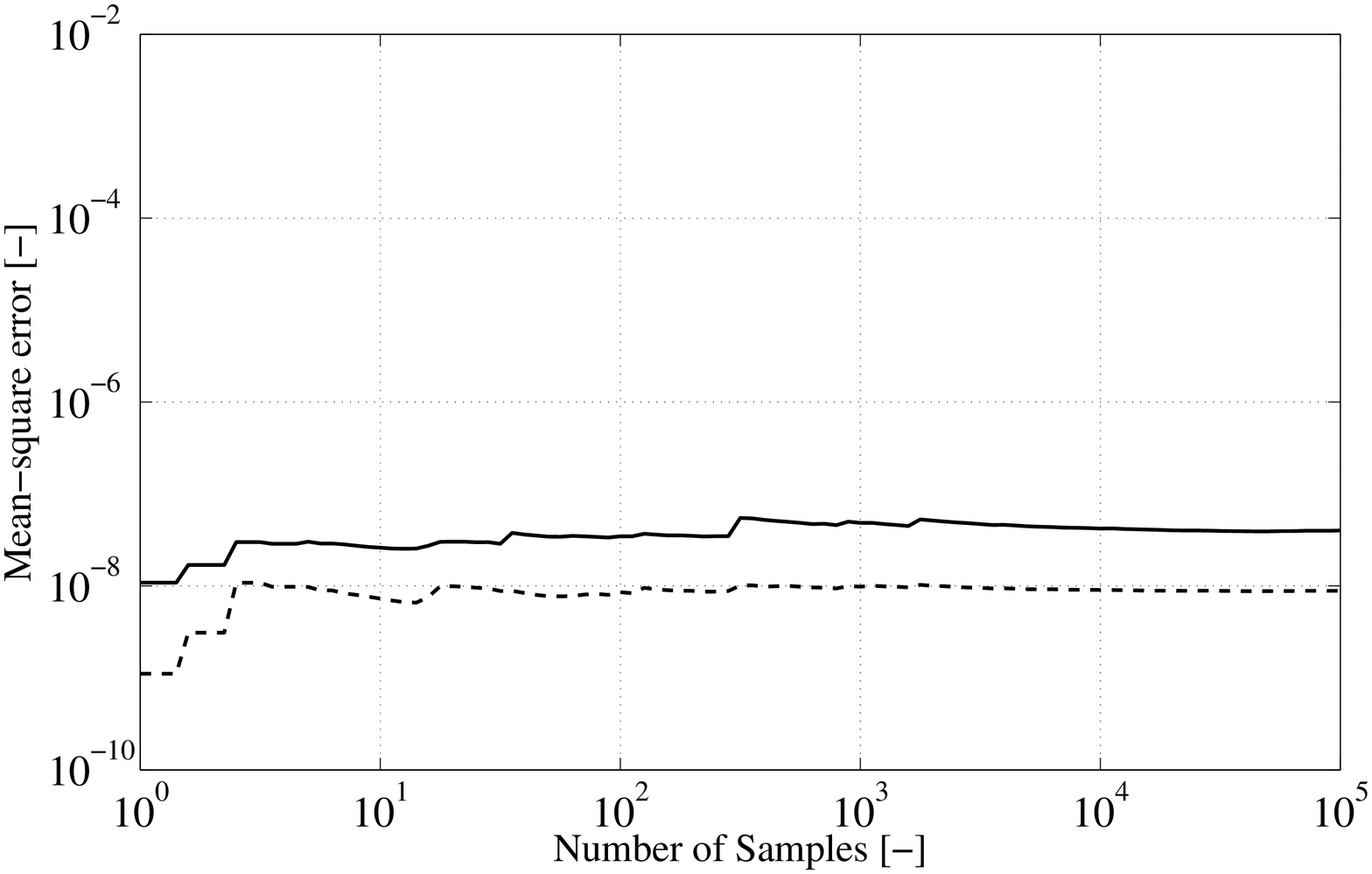}}
    \centerline{$\small{MC\mapsto\sqrt{\frac{1}{MC}\sum_{k=1}^{MC}\|\boldsymbol{T}^{\infty}(\boldsymbol{\xi}_{k})-{\boldsymbol{T}}{}^{\infty,p}(\boldsymbol{\xi}_{k})\|_{\boldsymbol{W}}^{2}}\hspace{-0.5mm}\Big/\hspace{-0.5mm}\sqrt{\frac{1}{MC}\sum_{k=1}^{MC}\|\boldsymbol{T}^{\infty}(\boldsymbol{\xi}_{k})\|_{\boldsymbol{W}}^{2}}}\normalsize$ (solid).}
    \centerline{$\small{MC\hspace{-1mm}\mapsto\hspace{-1mm}\sqrt{\frac{1}{MC}\sum_{k=1}^{MC}\|\boldsymbol{\Phi}^{\infty}(\boldsymbol{\xi}_{k})-{\boldsymbol{\Phi}}{}^{\infty,p}(\boldsymbol{\xi}_{k})\|_{\boldsymbol{W}}^{2}}\hspace{-1mm}\Big/\hspace{-1mm}\sqrt{\frac{1}{MC}\sum_{k=1}^{MC}\|\boldsymbol{\Phi}^{\infty}(\boldsymbol{\xi}_{k})\|_{\boldsymbol{W}}^{2}}}\normalsize$ (dashed).}
    \caption{Simulation not involving dimension reduction: convergence of the error estimates with respect to the number of samples for $p=4$.}\label{fig:figure6b}
  \end{center}
\end{figure}

Figure~\ref{fig:figure3} shows the convergence of the iterative method as a function of the number of iterations.  
The iterative method converged at a linear rate up to approximately iteration~$\ell=10$ for~$k=100\,\text{$[\text{J/K/cm/s}]$}$ and iteration~$\ell=15$ for~$k=1\,\text{$[\text{J/K/cm/s}]$}$, after which linear-solver tolerances became dominant and prevented further convergence.    

Figures~\ref{fig:figure6} and~\ref{fig:figure6b} show the convergence of the solution as a function of the total degree at which the PC expansions are truncated; note that the superscript $\infty$ is used in the figure captions to indicate convergence with respect to the number of iterations.
The distance between the solutions obtained through the Monte Carlo and the PC-based simulation that did not involve dimension reduction decreased at an exponential rate as the total degree was increased.

\subsection{PC-based implementation involving dimension reduction}
Lastly, we implemented the PC-based iterative method involving dimension reduction.
This implementation corresponded exactly to Algorithm~\ref{algo:algo6}.
We obtained the results to follow using PC expansions truncated at total degree $p=4$ and, with reference to~(\ref{eq:criterion}), using a range of values for the error tolerance level $tol$ adopted to determine the reduced dimension at each iteration.
We discuss convergence as a function of the error tolerance level later.
Now, we present detailed results obtained for~$tol=0.90\times \sigma_{\boldsymbol{T}}^{2}$, where~$\sigma_{\boldsymbol{T}}=\sqrt{\sum_{|\boldsymbol{\alpha}|=1}^{p}\|\boldsymbol{T}_{\boldsymbol{\alpha}}^{\infty}\|_{\boldsymbol{W}}^{2}}$ is the mean-square norm of the fluctuating part of the random temperature obtained by the PC-based iterative method that did not involve dimension reduction: specifically, $\sigma_{\boldsymbol{T}}=132.54\,[\text{JK/cm}^{2}/s]$ for $k=100\,\text{$[\text{J/K/cm/s}]$}$ and $\sigma_{\boldsymbol{T}}=201.18\,[\text{JK/cm}^{2}/s]$ for~$k=1\,\text{$[\text{J/K/cm/s}]$}$.

\begin{figure}[htp]
  \begin{center}
     \subfigure[Convergence for $k=100\,\text{$[\text{J/K/cm/s}]$}$.]{\includegraphics[width=0.49\textwidth]{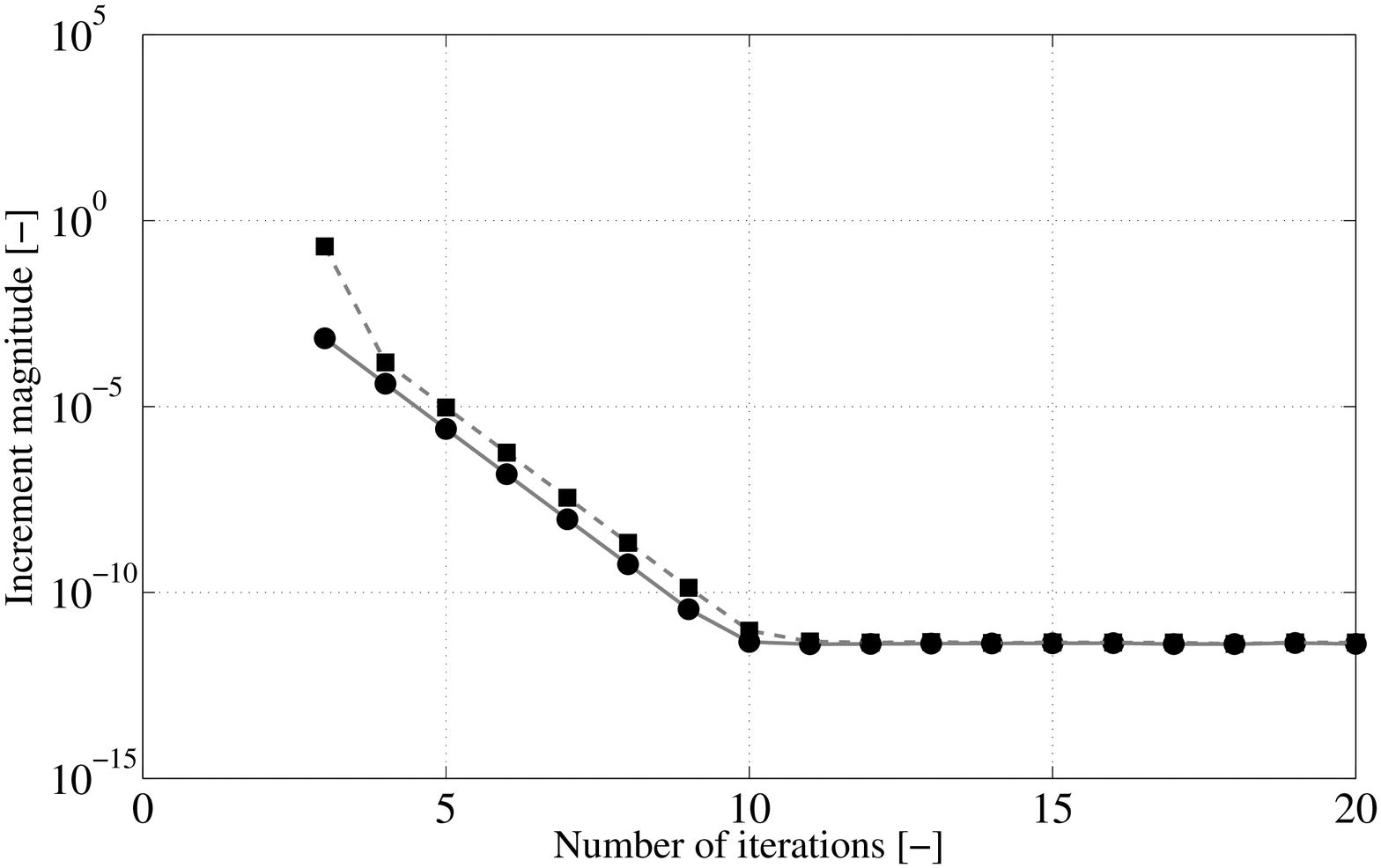}}
     \hfill
     \subfigure[Convergence for $k=1\,\text{$[\text{J/K/cm/s}]$}$.]{\includegraphics[width=0.49\textwidth]{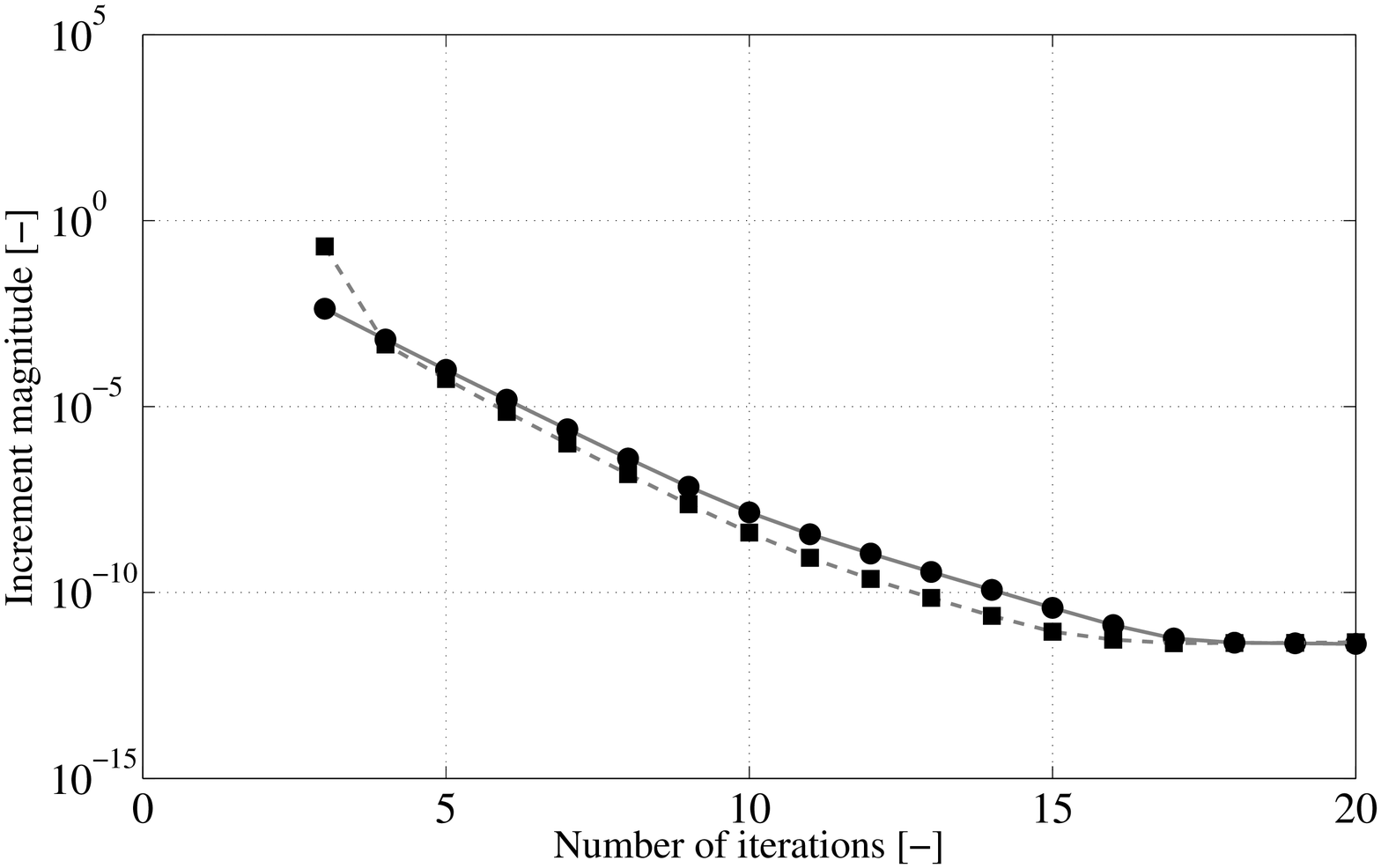}}
     \centerline{$\small{\ell\mapsto\sqrt{\sum_{|\boldsymbol{\alpha}|=0}^{p}\|\widehat{\boldsymbol{T}}{}_{\boldsymbol{\alpha}}^{\ell}-\widehat{\boldsymbol{T}}{}_{\boldsymbol{\alpha}}^{\ell-1}\|_{\boldsymbol{W}}^{2}}\Big/\sqrt{\sum_{|\boldsymbol{\alpha}|=0}^{p}\|\widehat{\boldsymbol{T}}{}_{\boldsymbol{\alpha}}^{\infty}\|_{\boldsymbol{W}}^{2}}}\normalsize$ (circles).}
     \centerline{$\small{\ell\mapsto\sqrt{\sum_{|\boldsymbol{\alpha}|=0}^{p}\|\widehat{\boldsymbol{\Phi}}{}_{\boldsymbol{\alpha}}^{\ell}-\widehat{\boldsymbol{\Phi}}{}_{\boldsymbol{\alpha}}^{\ell-1}\|_{\boldsymbol{W}}^{2}}\Big/\sqrt{\sum_{|\boldsymbol{\alpha}|=0}^{p}\|\widehat{\boldsymbol{\Phi}}{}_{\boldsymbol{\alpha}}^{\infty}\|_{\boldsymbol{W}}^{2}}}\normalsize$ (squares).}
    \caption{Simulation involving dimension reduction: convergence with respect to the number of iterations.}\label{fig:figure5}
  \end{center}
\end{figure}

Figure~\ref{fig:figure5} shows the convergence of the iterative method as a function of the number of iterations. 
As in the case shown in Fig.~\ref{fig:figure3}, the iterative method converged at a linear rate up to approximately iteration~$\ell=10$ for~$k=100\,\text{$[\text{J/K/cm/s}]$}$ and iteration~$\ell=15$ for~$k=1\,\text{$[\text{J/K/cm/s}]$}$.

\begin{figure}[htp]
  \begin{center}
    \subfigure[Mean for $k=100\,\text{$[\text{J/K/cm/s}]$}$.]{\includegraphics[width=0.49\textwidth]{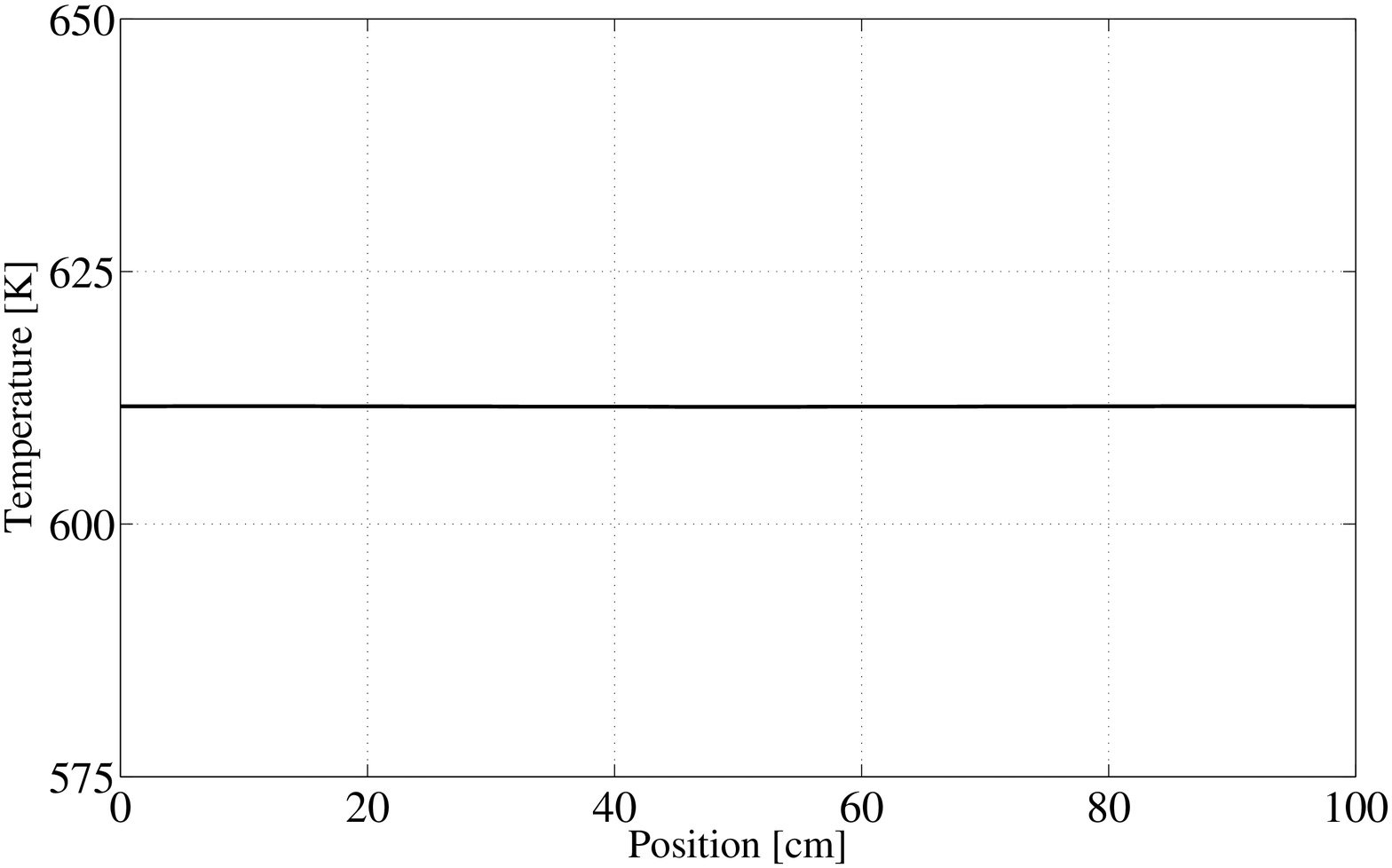}}
    \hfill
    \subfigure[Mean for $k=1\,\text{$[\text{J/K/cm/s}]$}$.]{\includegraphics[width=0.49\textwidth]{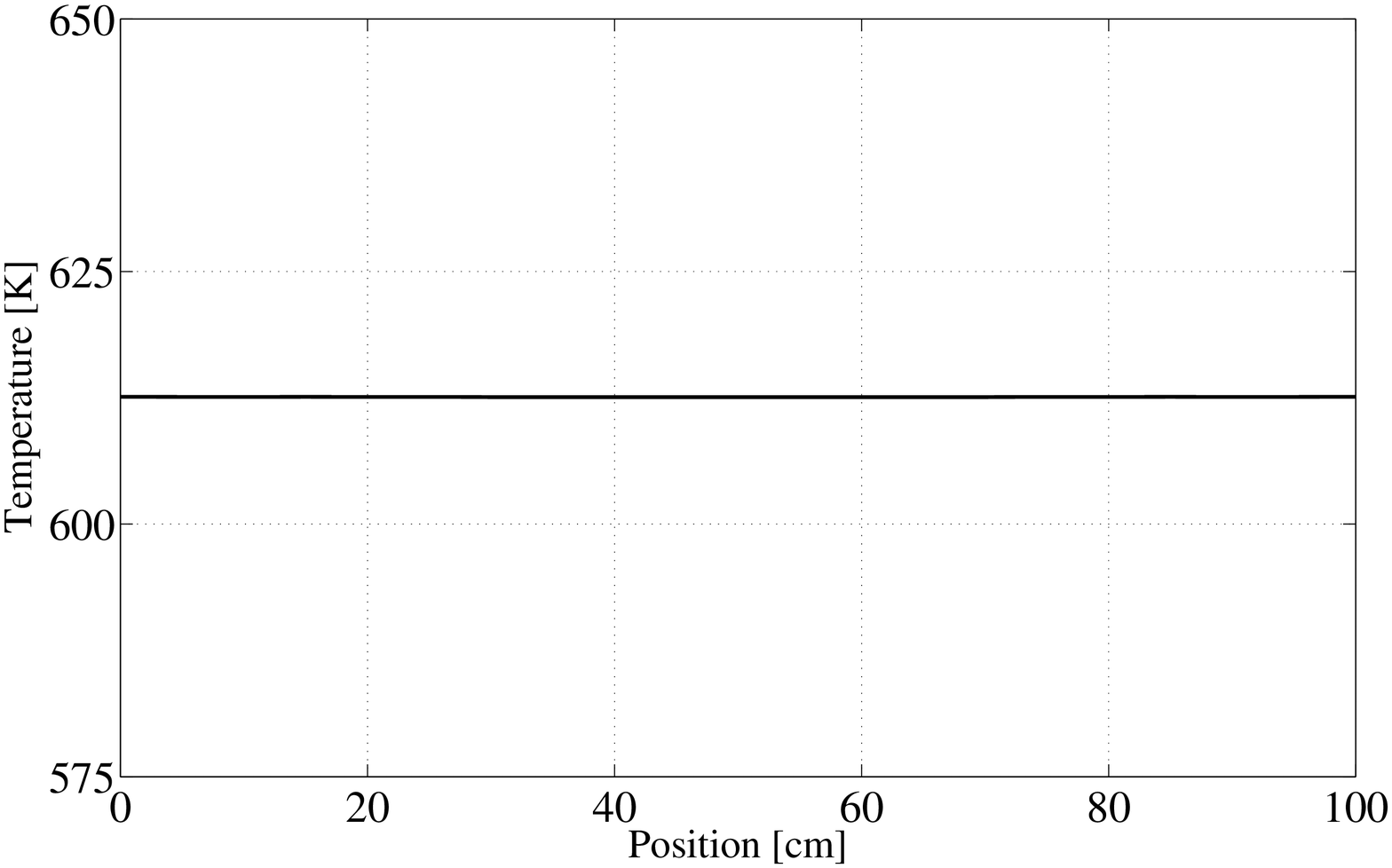}}    
    \subfigure[Eigenvalues for $k=100\,\text{$[\text{J/K/cm/s}]$}$.]{\includegraphics[width=0.49\textwidth]{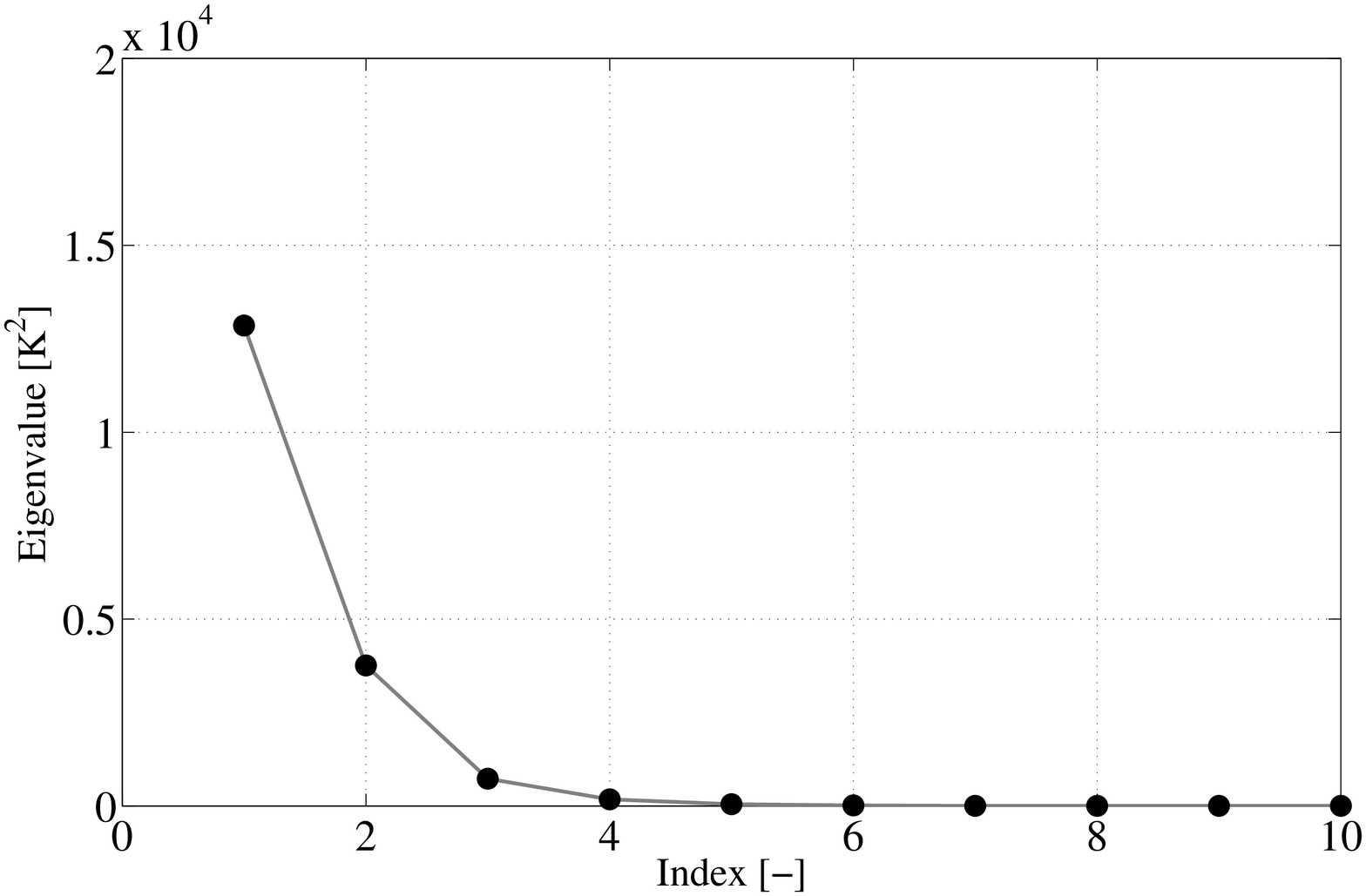}}
    \hfill
    \subfigure[Eigenvalues for $k=1\,\text{$[\text{J/K/cm/s}]$}$.]{\includegraphics[width=0.49\textwidth]{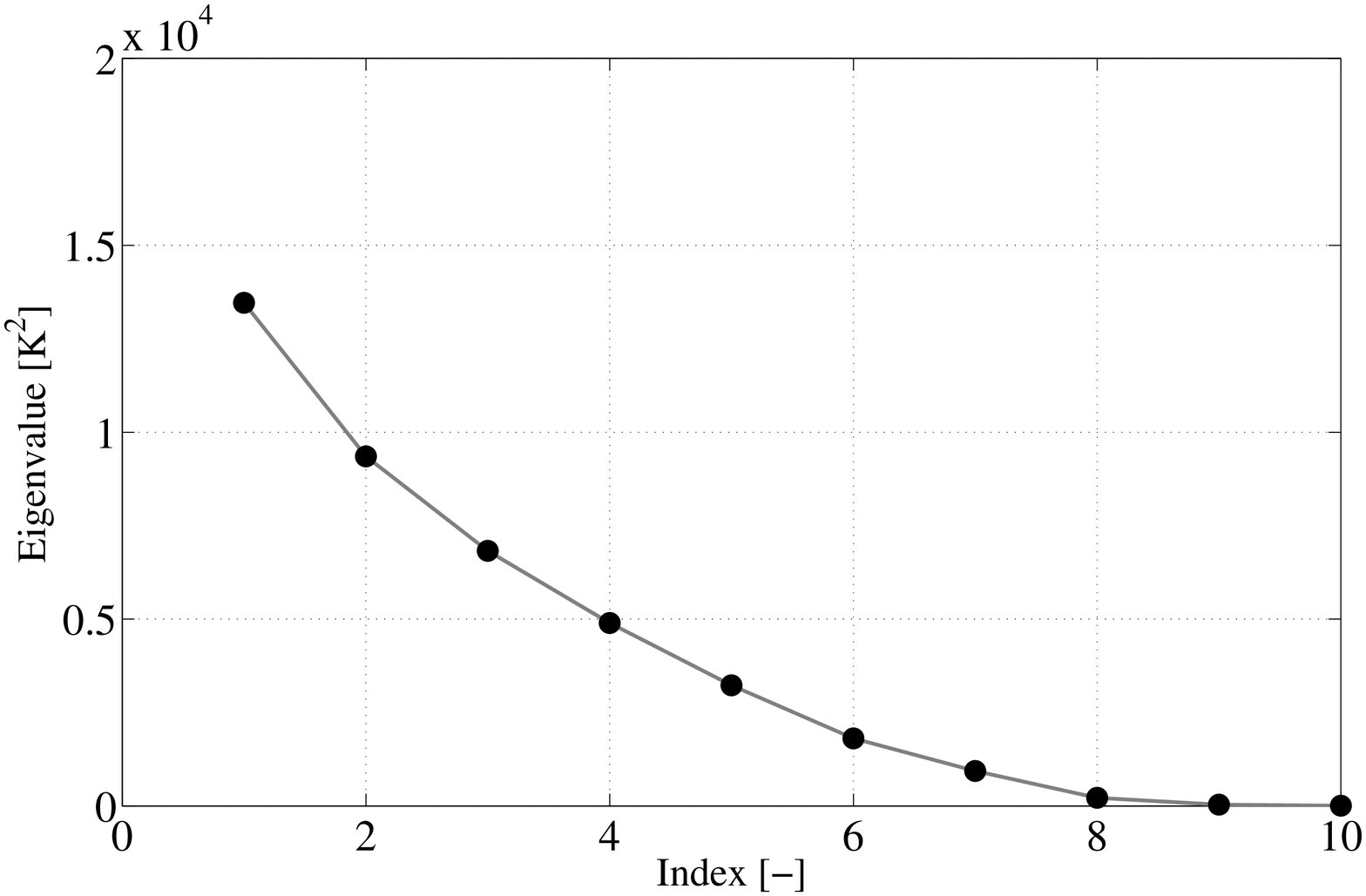}}
    \subfigure[First eigenmode for $k=100\,\text{$[\text{J/K/cm/s}]$}$.]{\includegraphics[width=0.49\textwidth]{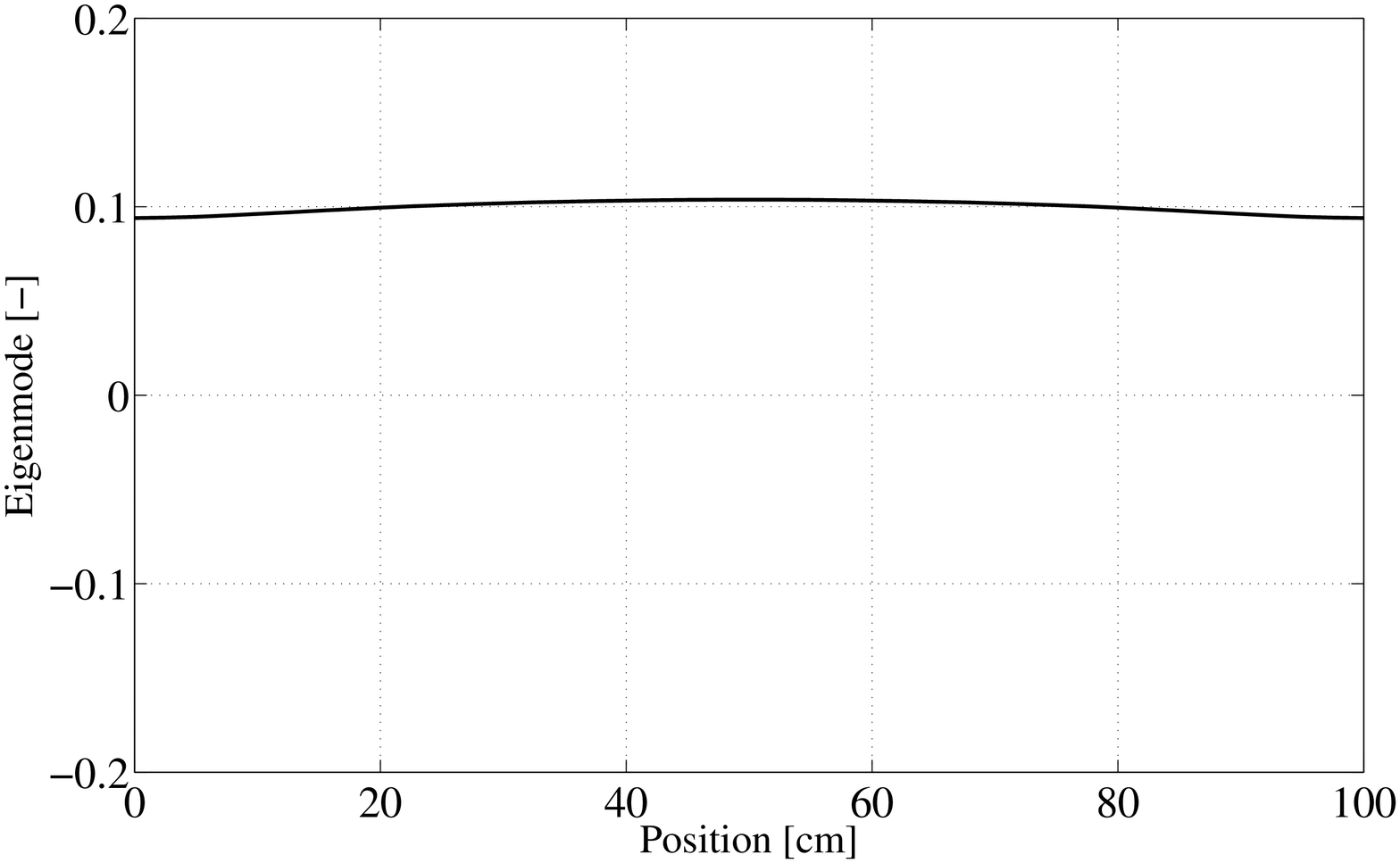}}
    \hfill
    \subfigure[First eigenmode for $k=1\,\text{$[\text{J/K/cm/s}]$}$.]{\includegraphics[width=0.49\textwidth]{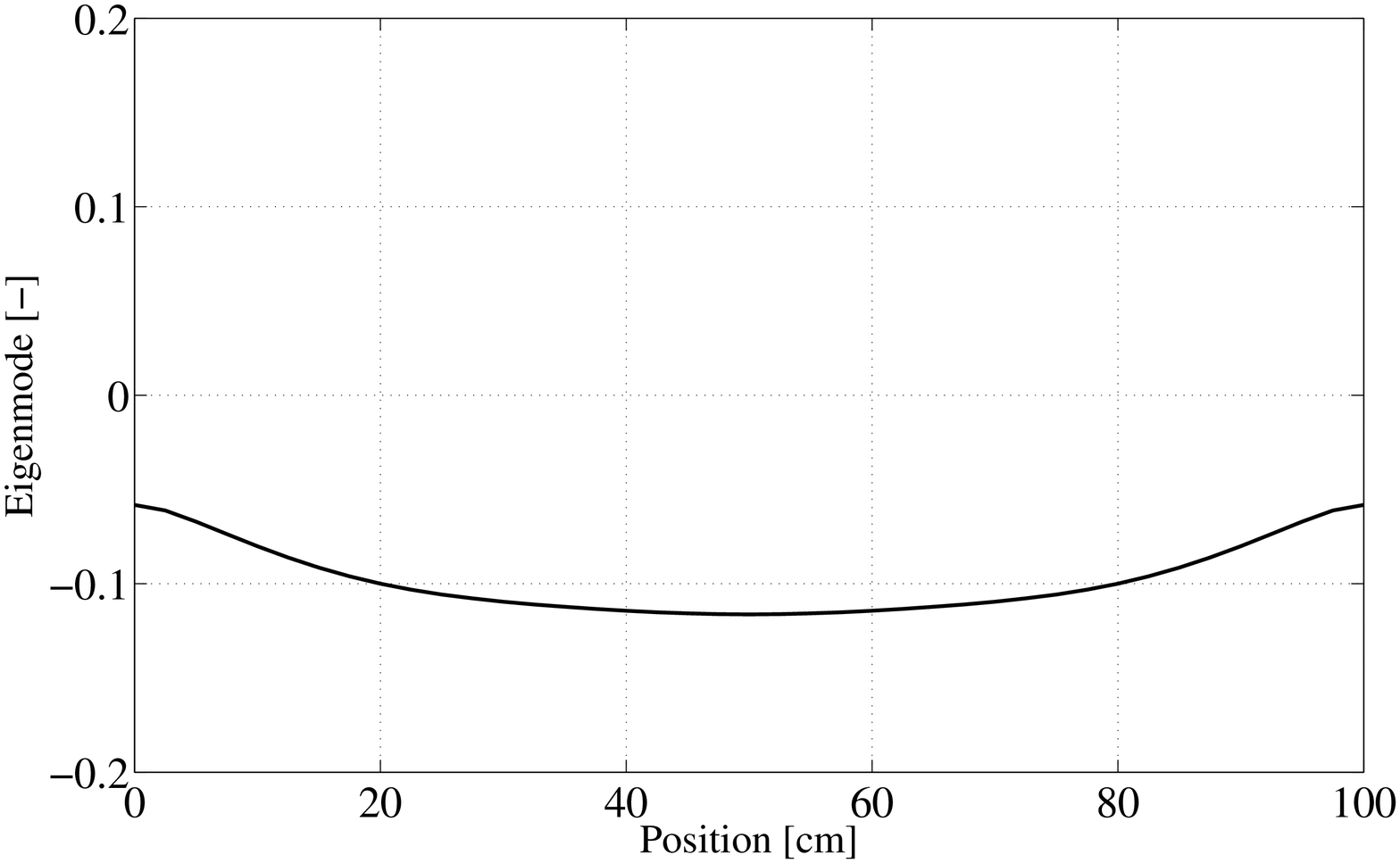}}
    \caption{Simulation involving dimension reduction: mean, ten largest magnitude eigenvalues, and first eigenmode of the KL~decomposition of the random temperature.}\label{fig:figure7}
  \end{center}
\end{figure}

\begin{figure}[htp]
  \begin{center}
    \subfigure[Second eigenmode for $k=100\,\text{$[\text{J/K/cm/s}]$}$.]{\includegraphics[width=0.49\textwidth]{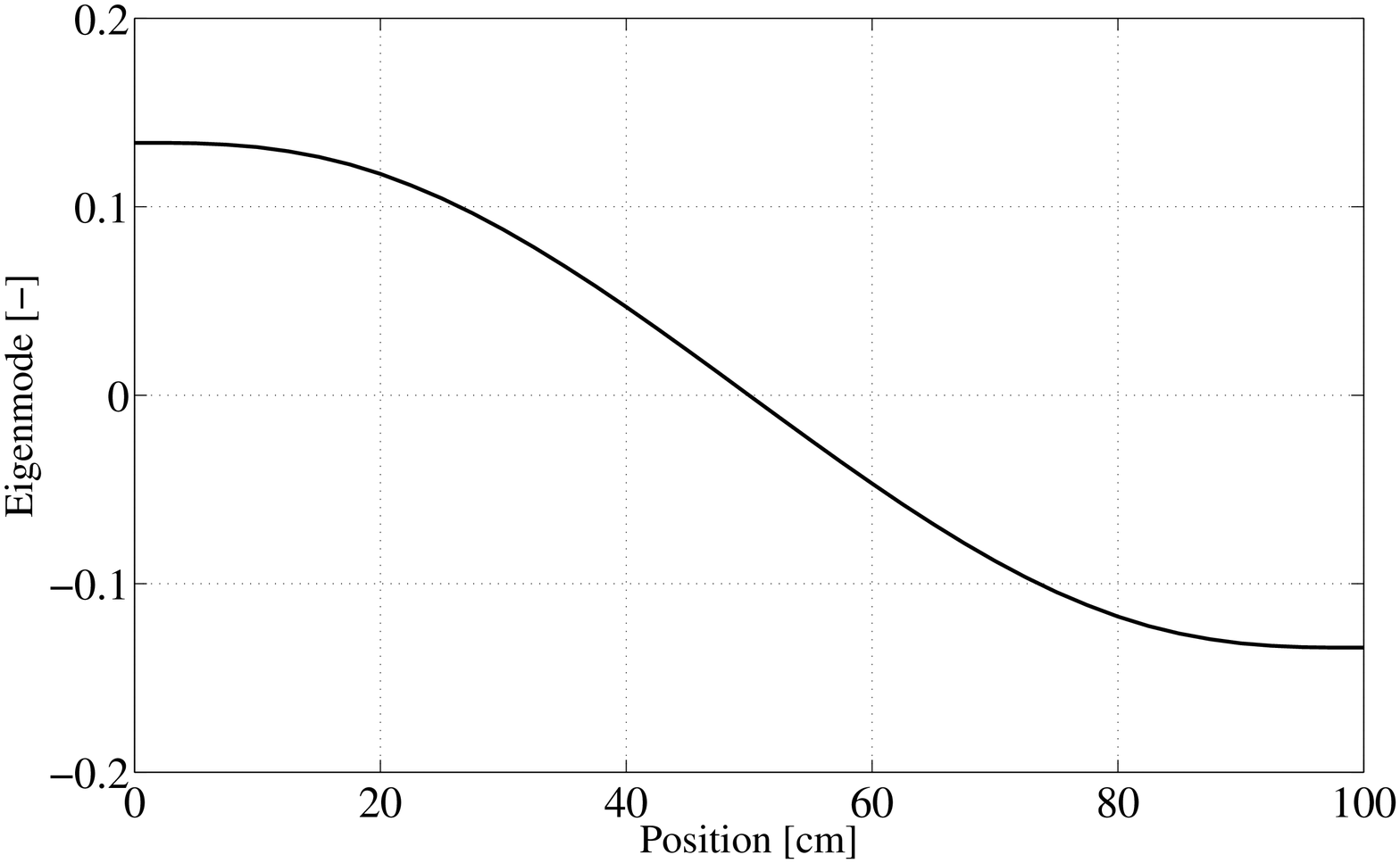}}
    \hfill
    \subfigure[Second eigenmode for $k=1\,\text{$[\text{J/K/cm/s}]$}$.]{\includegraphics[width=0.49\textwidth]{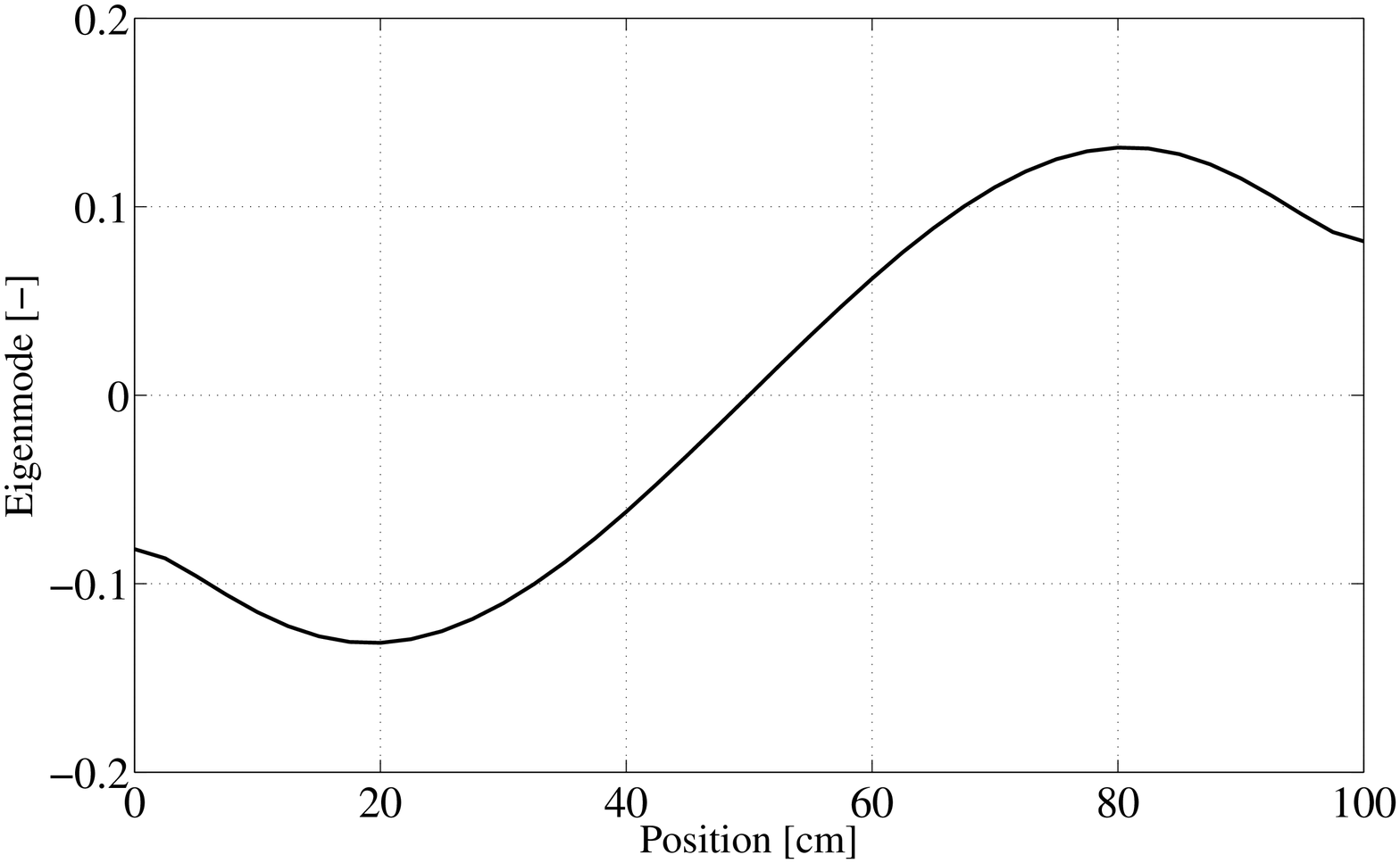}}    
    \subfigure[First reduced r.v. for $k=100\,\text{$[\text{J/K/cm/s}]$}$.]{\includegraphics[width=0.49\textwidth]{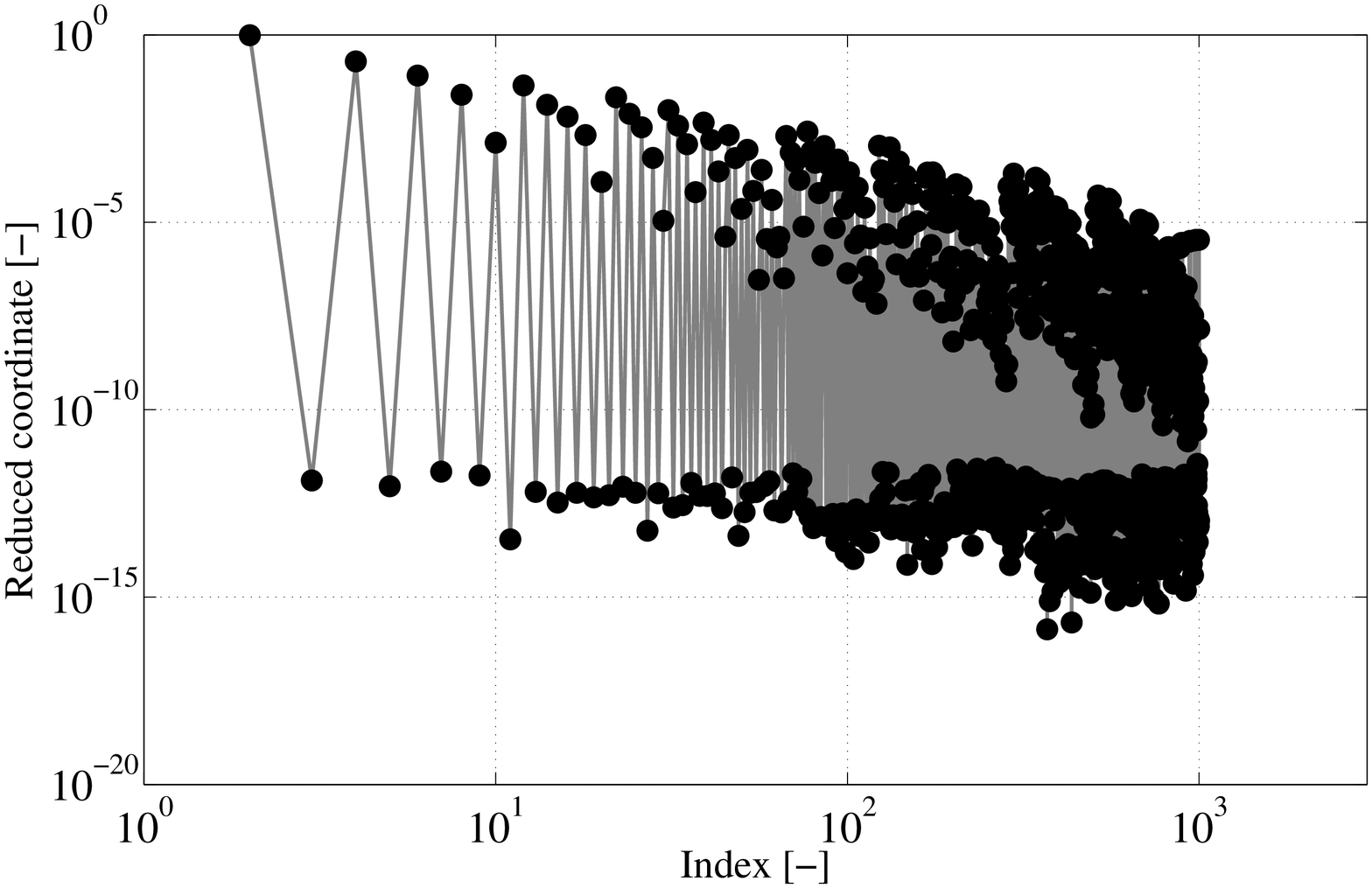}}
    \hfill
    \subfigure[First reduced r.v. for $k=1\,\text{$[\text{J/K/cm/s}]$}$.]{\includegraphics[width=0.49\textwidth]{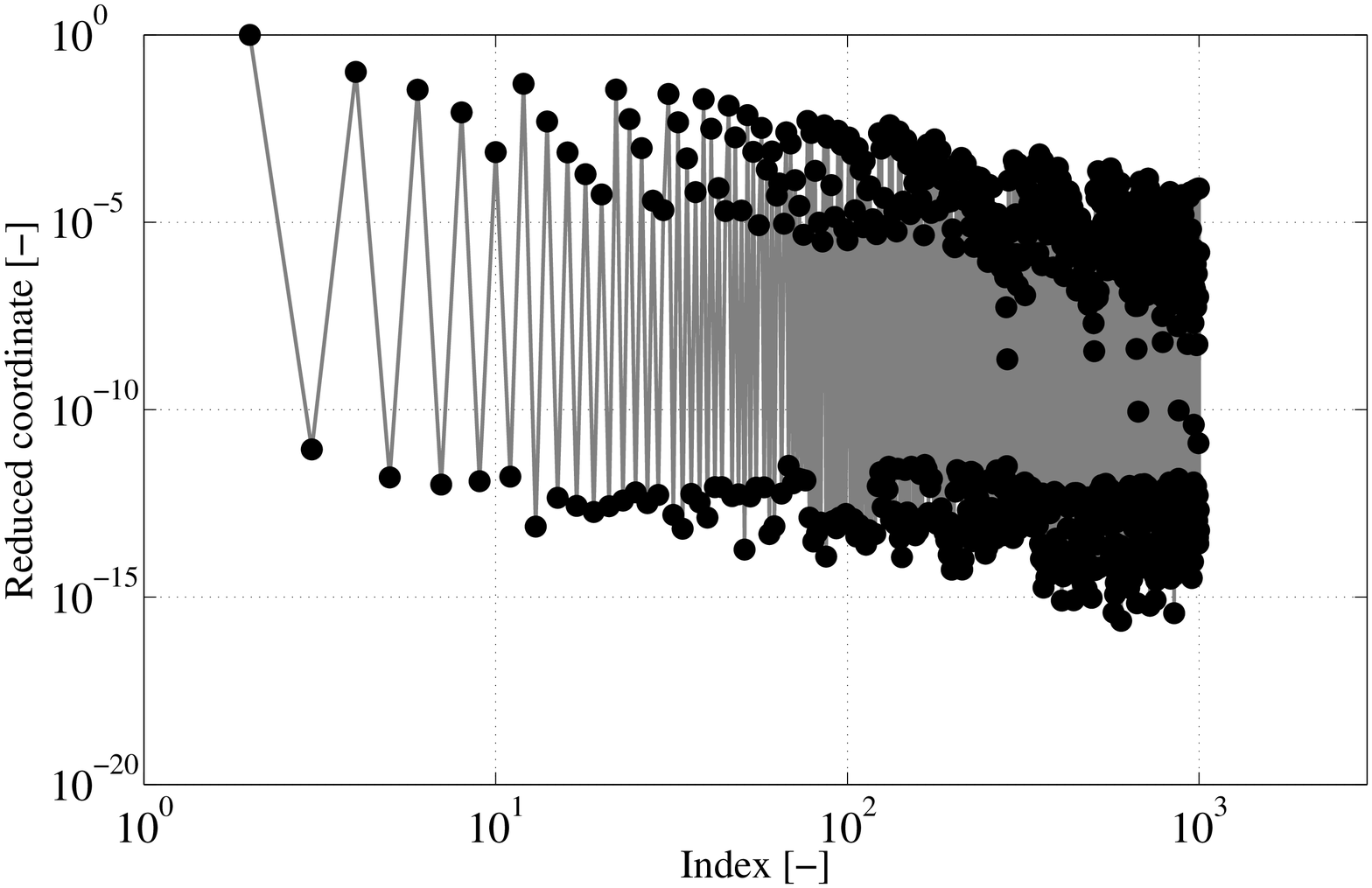}}
    \subfigure[Second reduced r.v. for $k=100\,\text{$[\text{J/K/cm/s}]$}$.]{\includegraphics[width=0.49\textwidth]{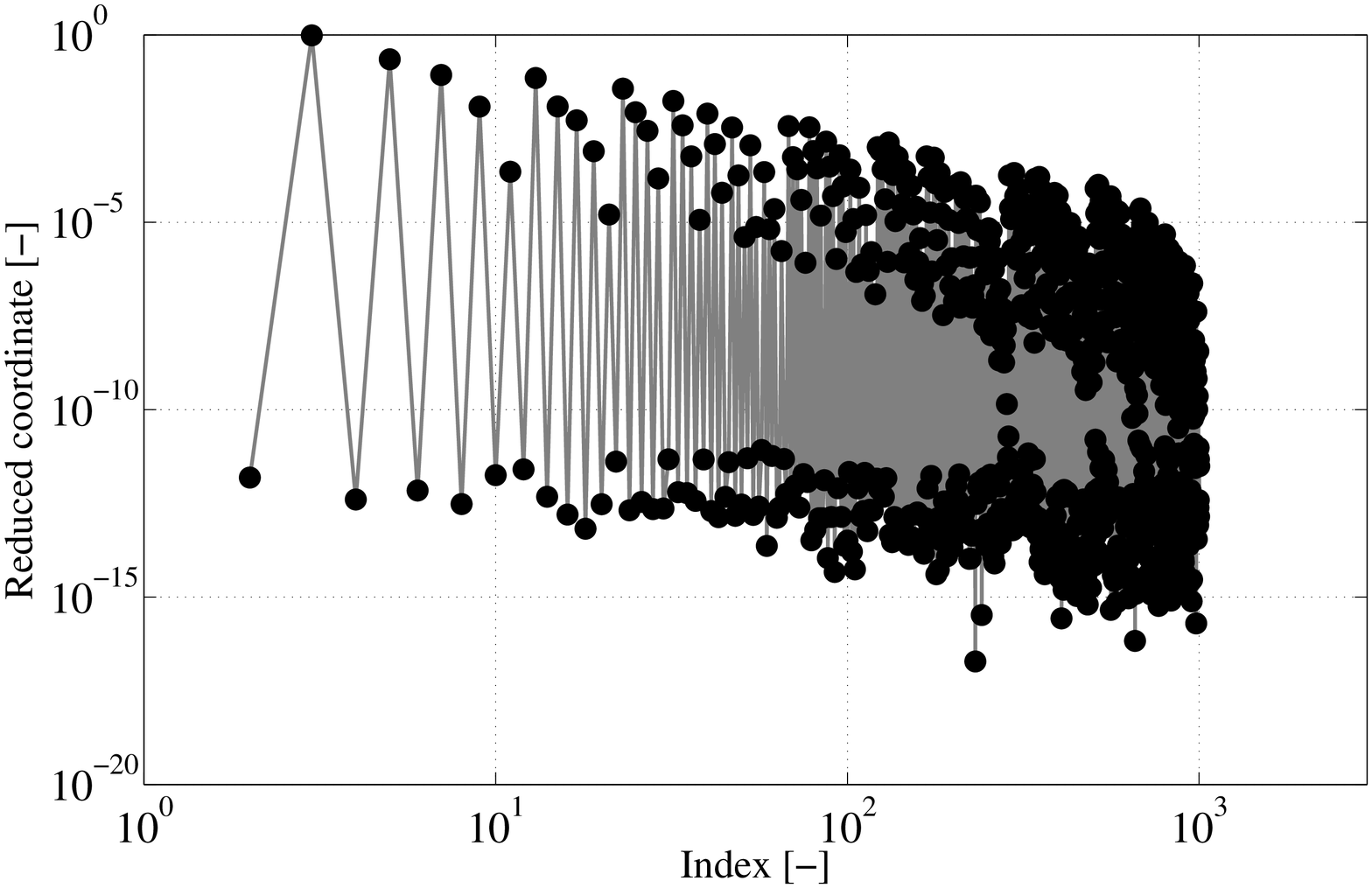}}
    \hfill
    \subfigure[Second reduced r.v. for $k=1\,\text{$[\text{J/K/cm/s}]$}$.]{\includegraphics[width=0.49\textwidth]{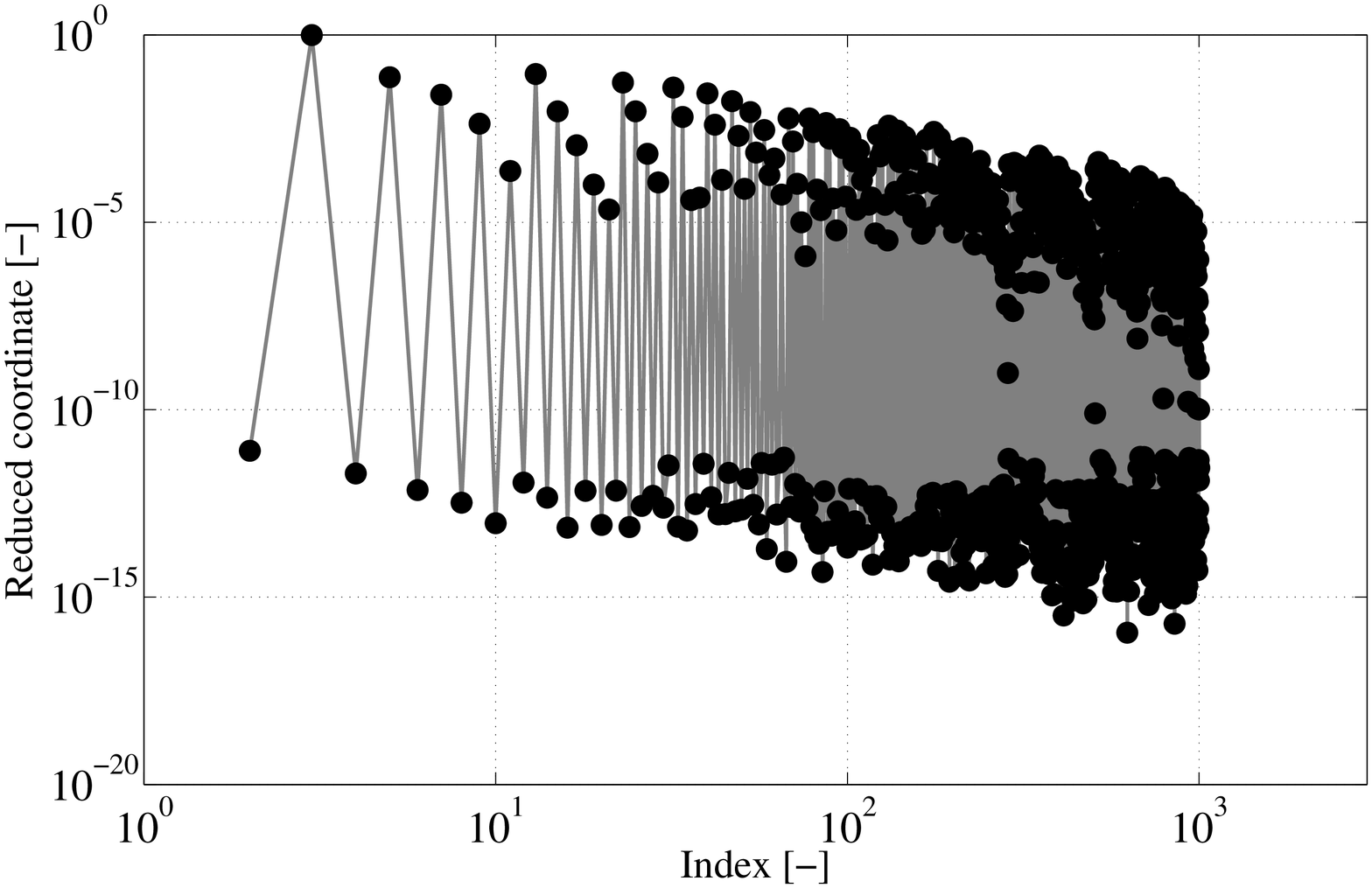}}
    \caption{Simulation involving dimension reduction: second eigenmode and PC coordinates of the first and second reduced random variables of the KL~decomposition of the random temperature.}\label{fig:figure7b}
  \end{center}
\end{figure}

Figures~\ref{fig:figure7} and~\ref{fig:figure7b} show a few components of the KL decomposition of the random temperature obtained at the last iteration, namely, the mean temperature, the 10 largest magnitude eigenvalues, and the eigenmodes and reduced random variables associated with the two largest magnitude eigenvalues.
We can observe that the eigenvalues obtained for $k=100\,\text{$[\text{J/K/cm/s}]$}$ (Fig.~\ref{fig:figure7}(c)) decay at a faster rate than those obtained for $k=1\,\text{$[\text{J/K/cm/s}]$}$ (Fig.~\ref{fig:figure7}(d)).
This behavior of the eigenvalue decay rate is consistent with our earlier observation that the samples of the random temperature become less smooth and therefore exhibits less spatial correlation as the thermal conductivity is decreased.

\begin{figure}[htp]
  \begin{center}
    \subfigure[Temperature for $k=100\,\text{$[\text{J/K/cm/s}]$}$.]{\includegraphics[width=0.49\textwidth]{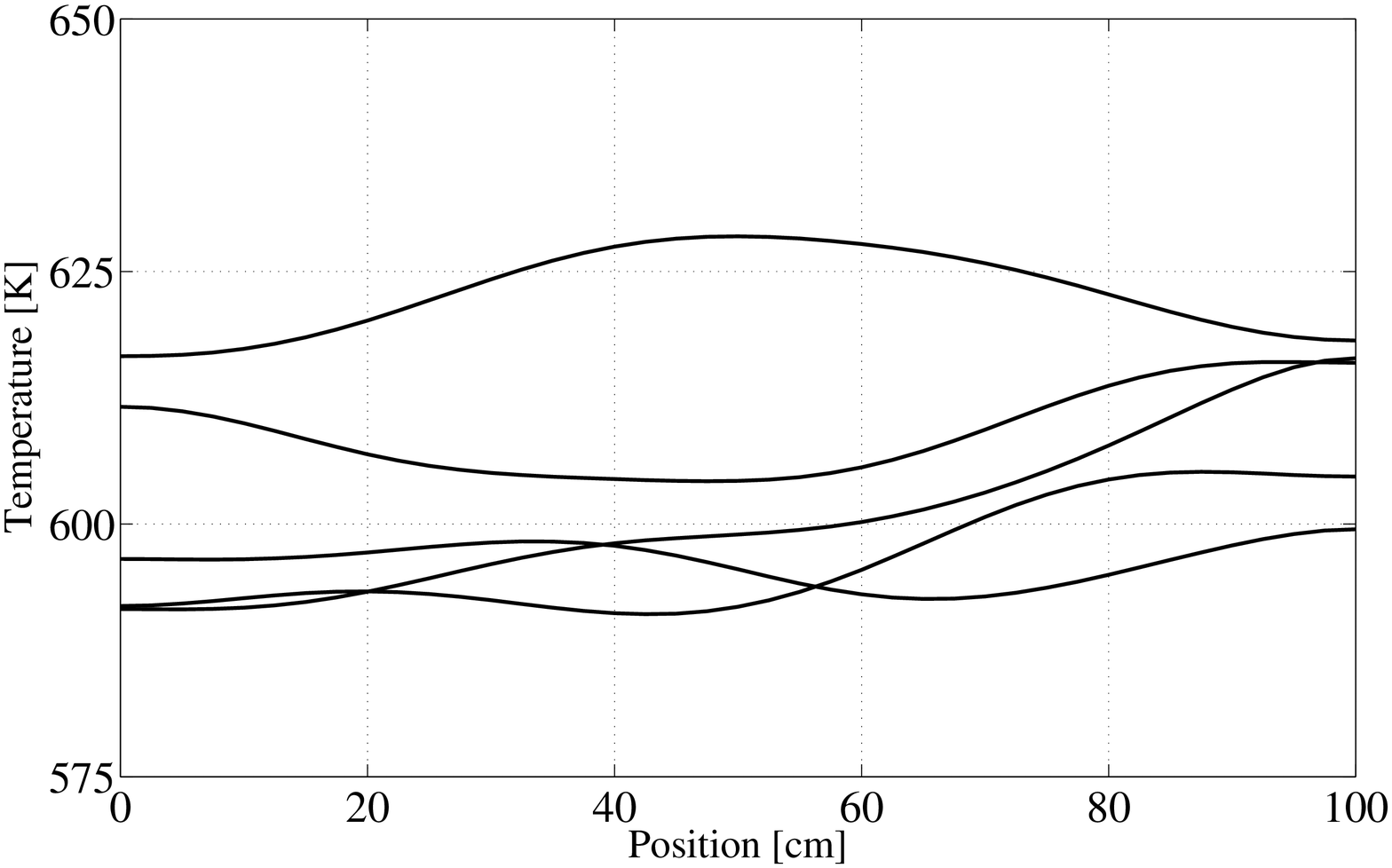}}
    \hfill
    \subfigure[Temperature for $k=1\,\text{$[\text{J/K/cm/s}]$}$.]{\includegraphics[width=0.49\textwidth]{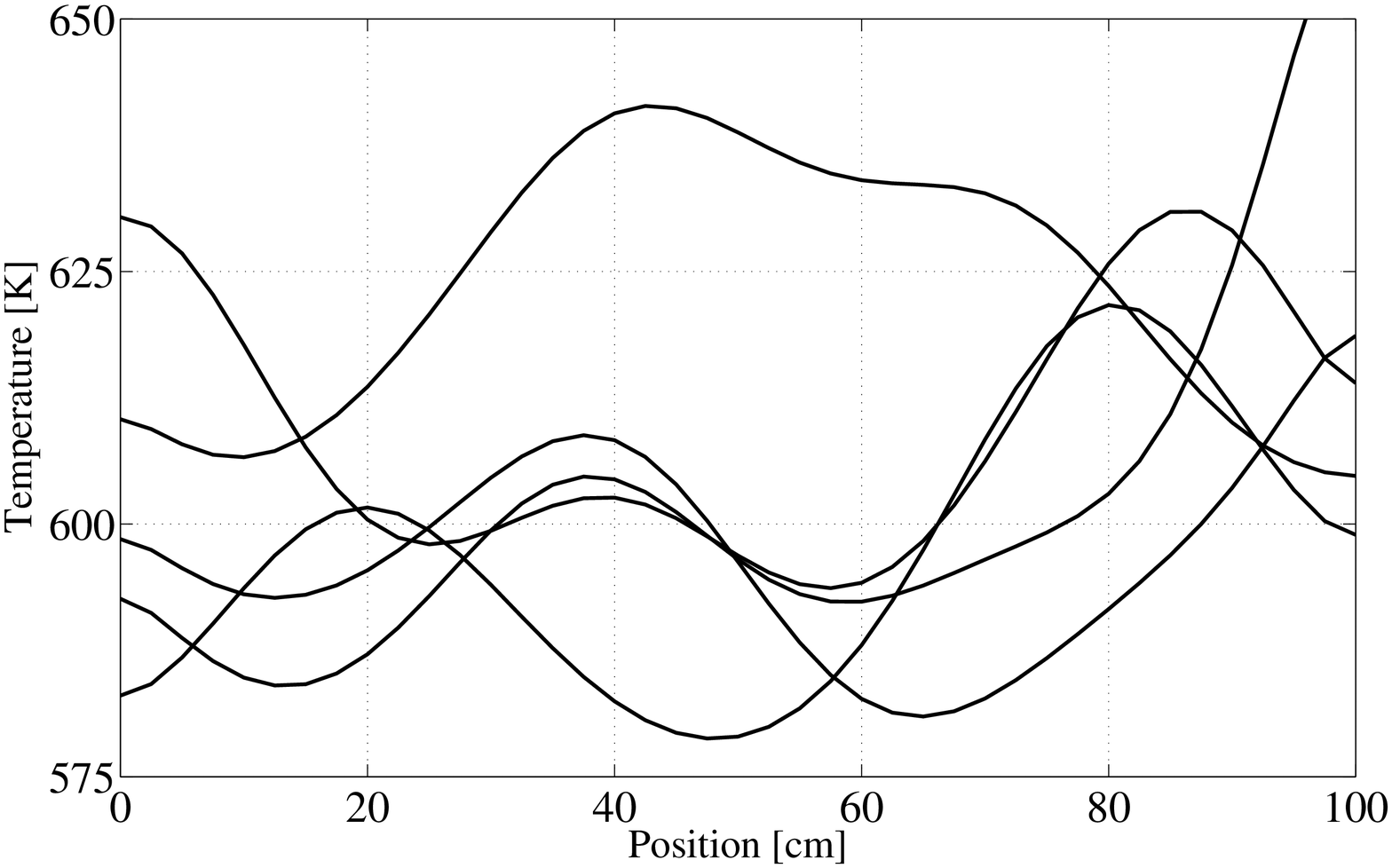}}
    \subfigure[Neutron flux for $k=100\,\text{$[\text{J/K/cm/s}]$}$.]{\includegraphics[width=0.49\textwidth]{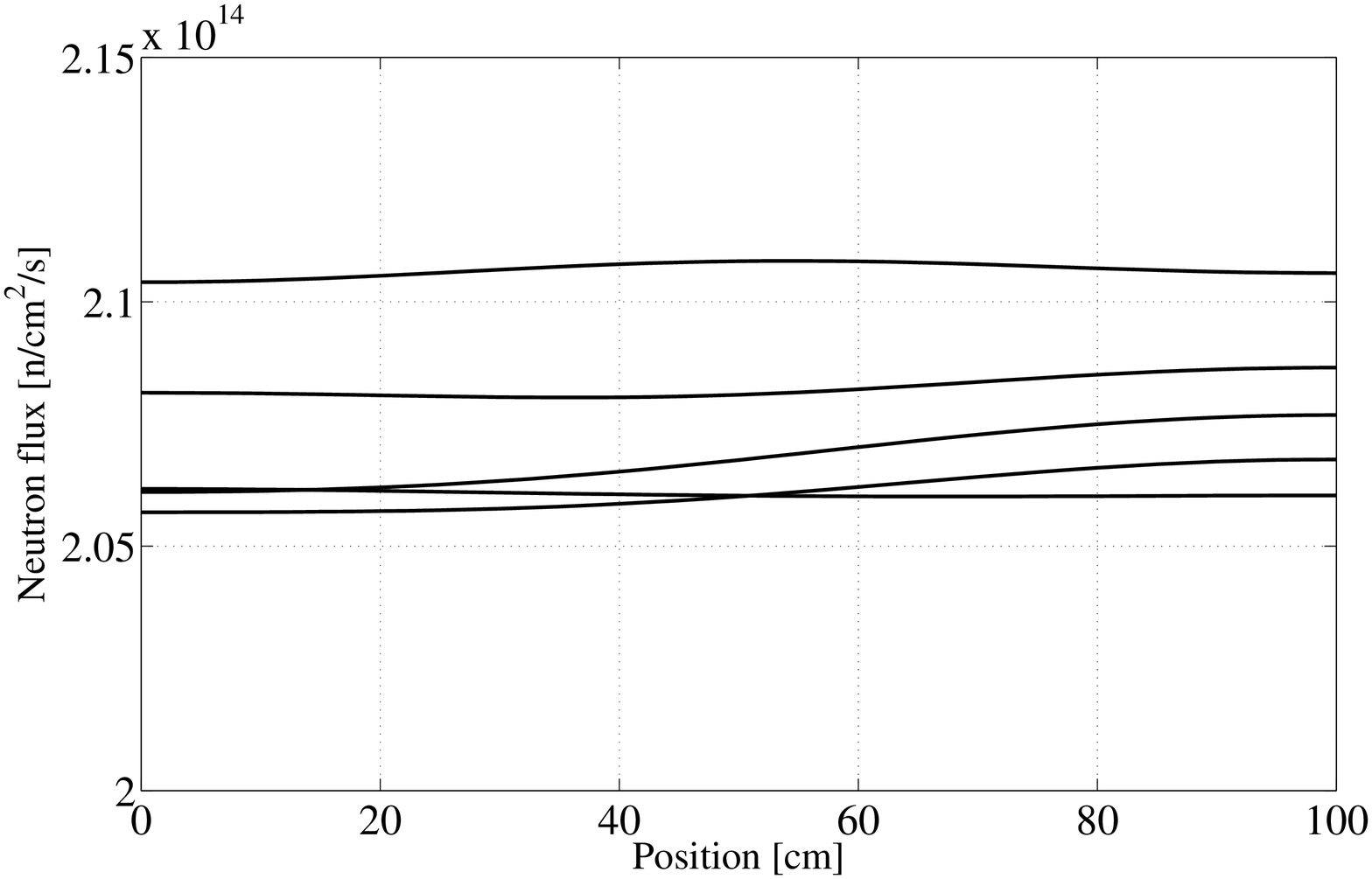}}
    \hfill
    \subfigure[Neutron flux for $k=1\,\text{$[\text{J/K/cm/s}]$}$.]{\includegraphics[width=0.49\textwidth]{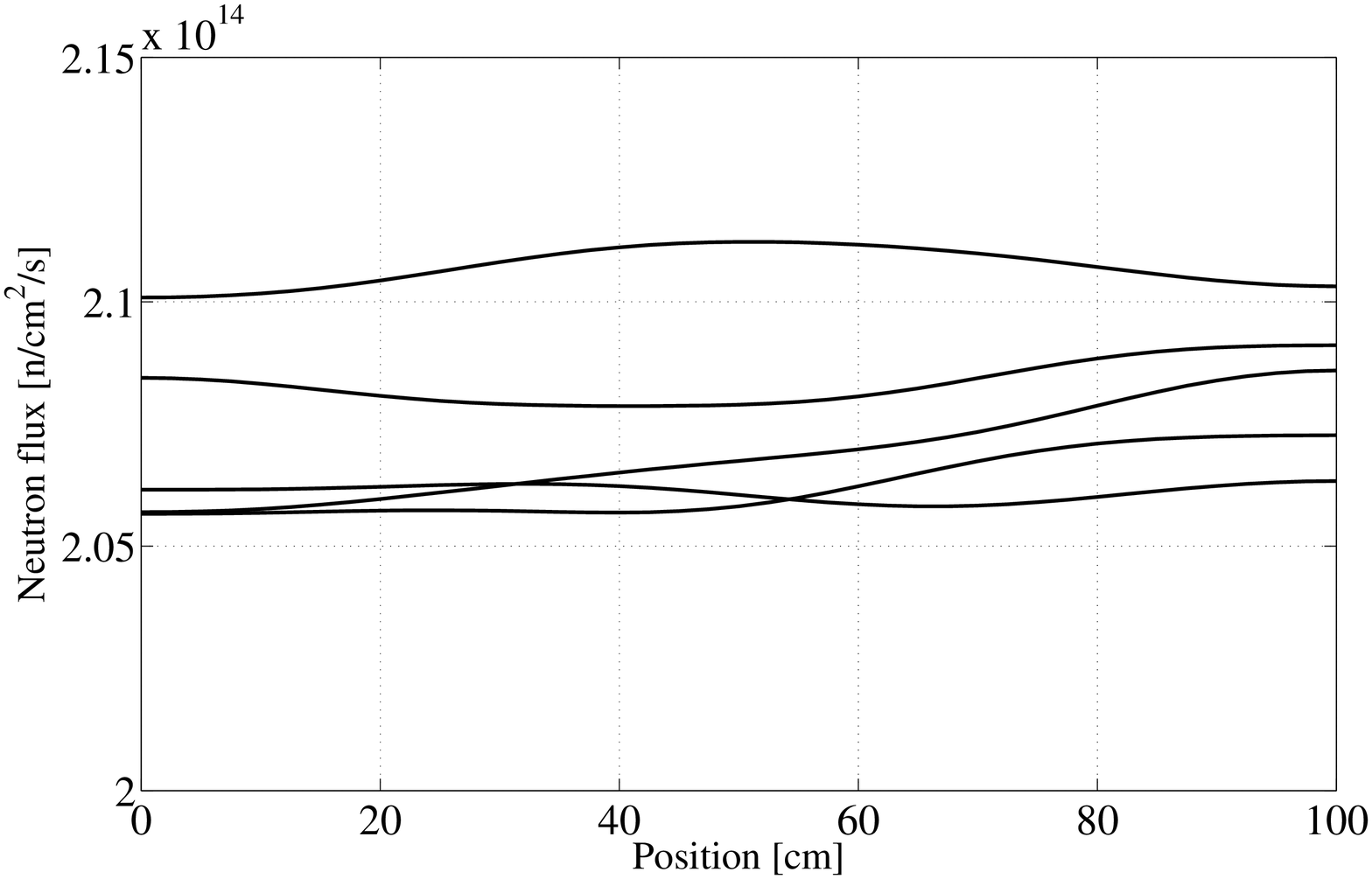}}
    \caption{Simulation involving dimension reduction: five samples of the solution.}\label{fig:figure10}
  \end{center}
\end{figure}

Figure~\ref{fig:figure10} shows a few samples of the random temperature and neutron flux deduced from the PC expansions obtained as the output of the solution algorithm.
The samples of the input random variables used to synthesize the samples of the random temperature and neutron flux shown in Fig.~\ref{fig:figure10} were identical to those used to generate the samples shown in Fig.~\ref{fig:figure2}.
The similarity of the samples in Figs.~\ref{fig:figure2} and~\ref{fig:figure12} indicates that the PC-based surrogate model not only provides an accurate \textit{global} representation of the multiphysics model but is also capable of accurately reproducing a \textit{sample-wise} response.

\subsection{Convergence analysis}

\begin{figure}[htp]
  \begin{center}
    \subfigure[Reduced dimension for $k=100\,\text{$[\text{J/K/cm/s}]$}$.]{\includegraphics[width=0.49\textwidth]{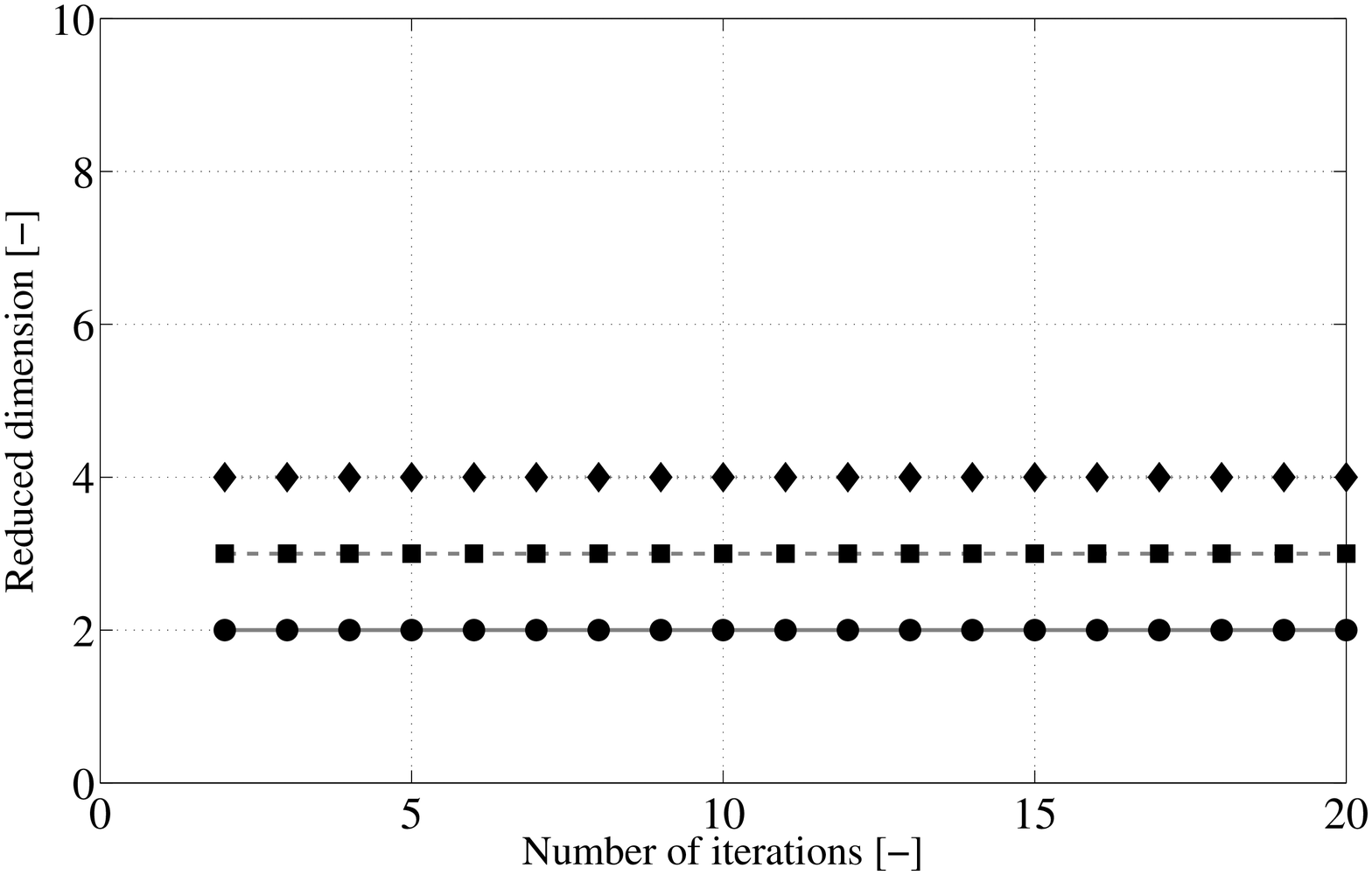}}
    \hfill
    \subfigure[Reduced dimension for $k=1\,\text{$[\text{J/K/cm/s}]$}$.]{\includegraphics[width=0.49\textwidth]{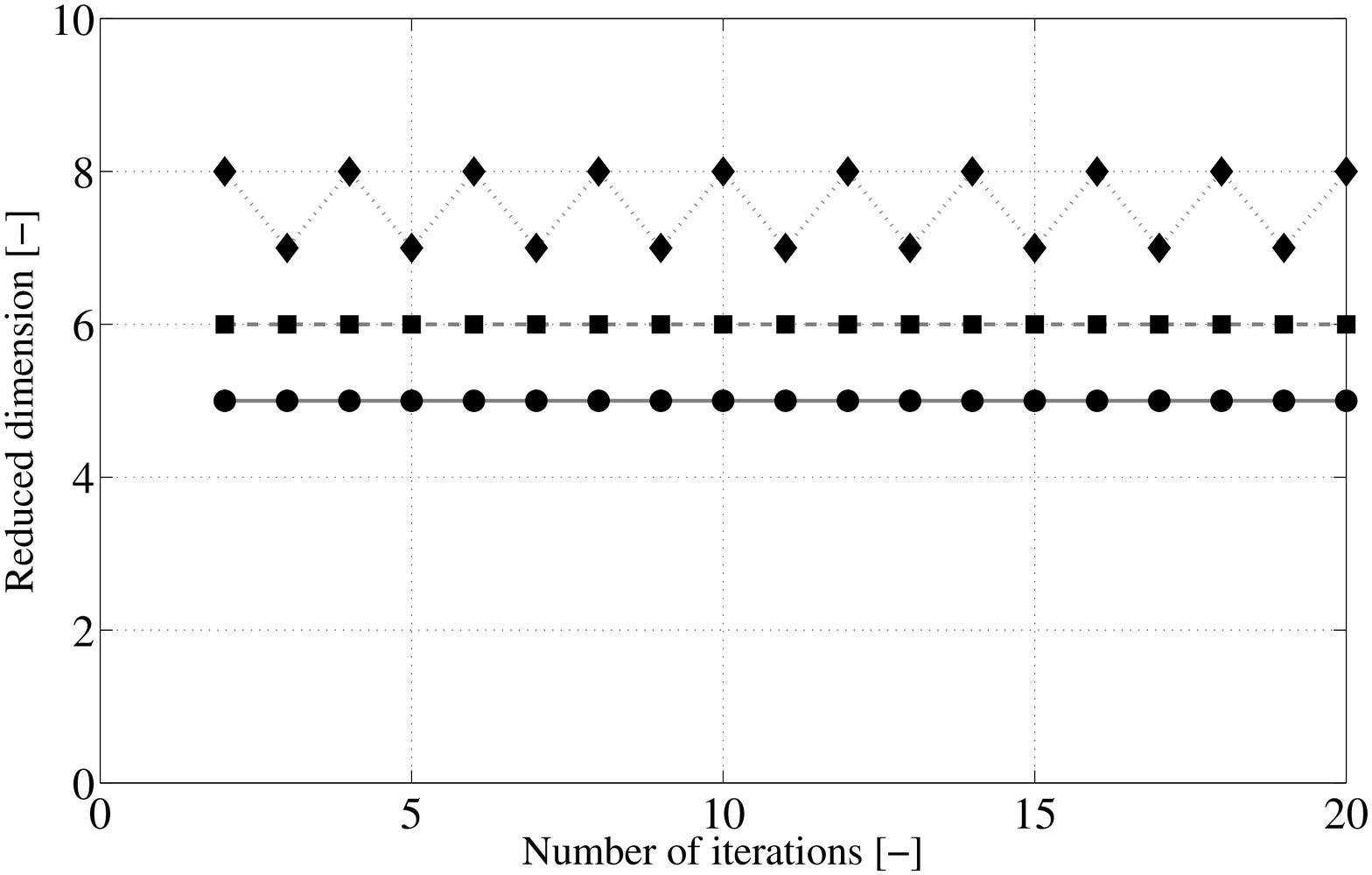}}
    \caption{Convergence analysis: reduced dimension as a function of the iteration for~$tol=0.90\times \sigma_{\boldsymbol{T}}^{2}$ (circles), $tol=0.95\times \sigma_{\boldsymbol{T}}^{2}$ (squares), and~$tol=0.99\times \sigma_{\boldsymbol{T}}^{2}$ (diamonds).}\label{fig:figure12}
  \end{center}
\end{figure}

\begin{figure}[htp]
  \begin{center}
    \subfigure[\hspace{-1mm}Temperature-based distance for $k\hspace{-1mm}=\hspace{-1mm}100\,\text{$[\text{J/K/cm/s}]$}$.]{\includegraphics[width=0.49\textwidth]{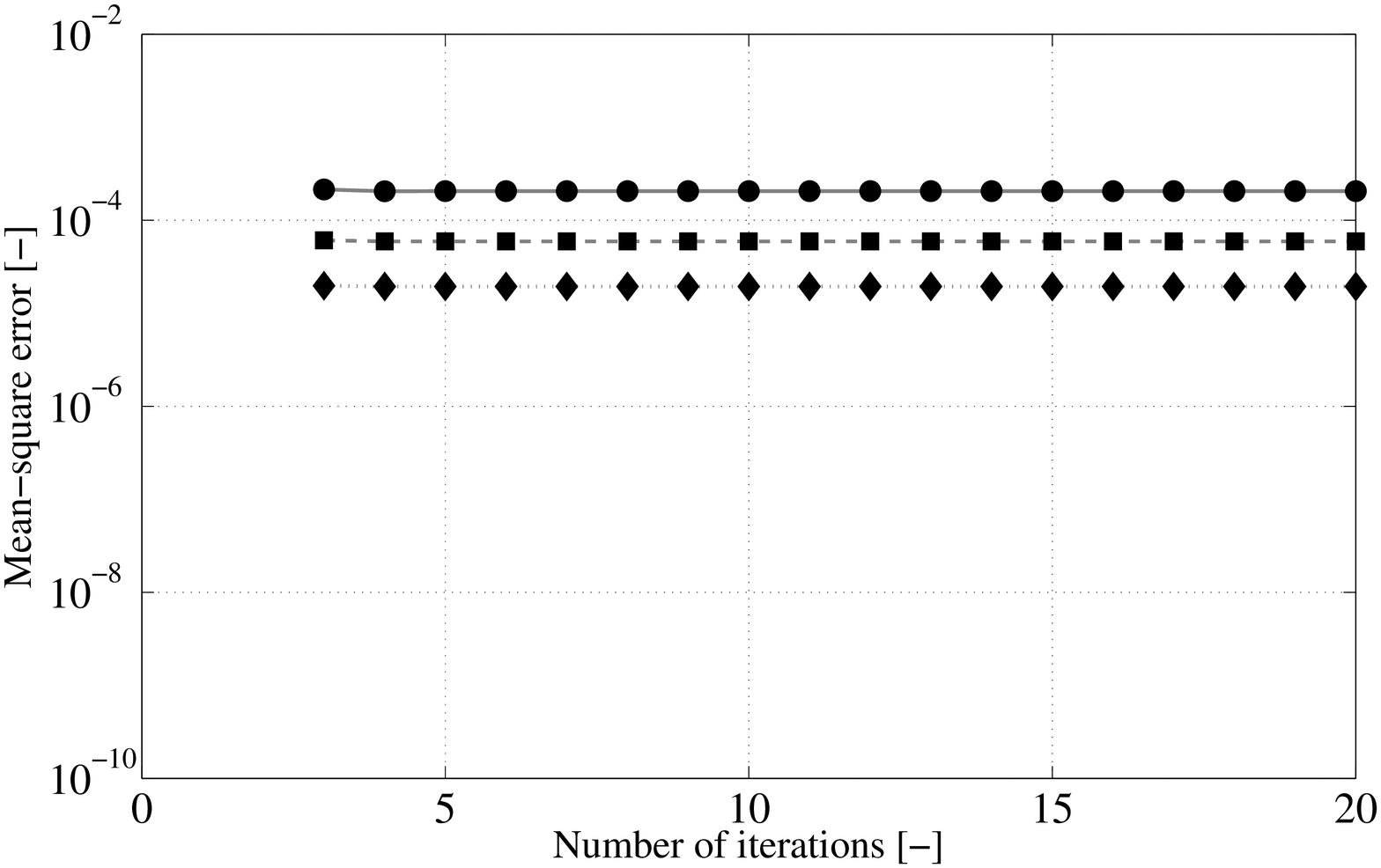}}
    \hfill
    \subfigure[Temperature-based distance for $k=1\,\text{$[\text{J/K/cm/s}]$}$.]{\includegraphics[width=0.49\textwidth]{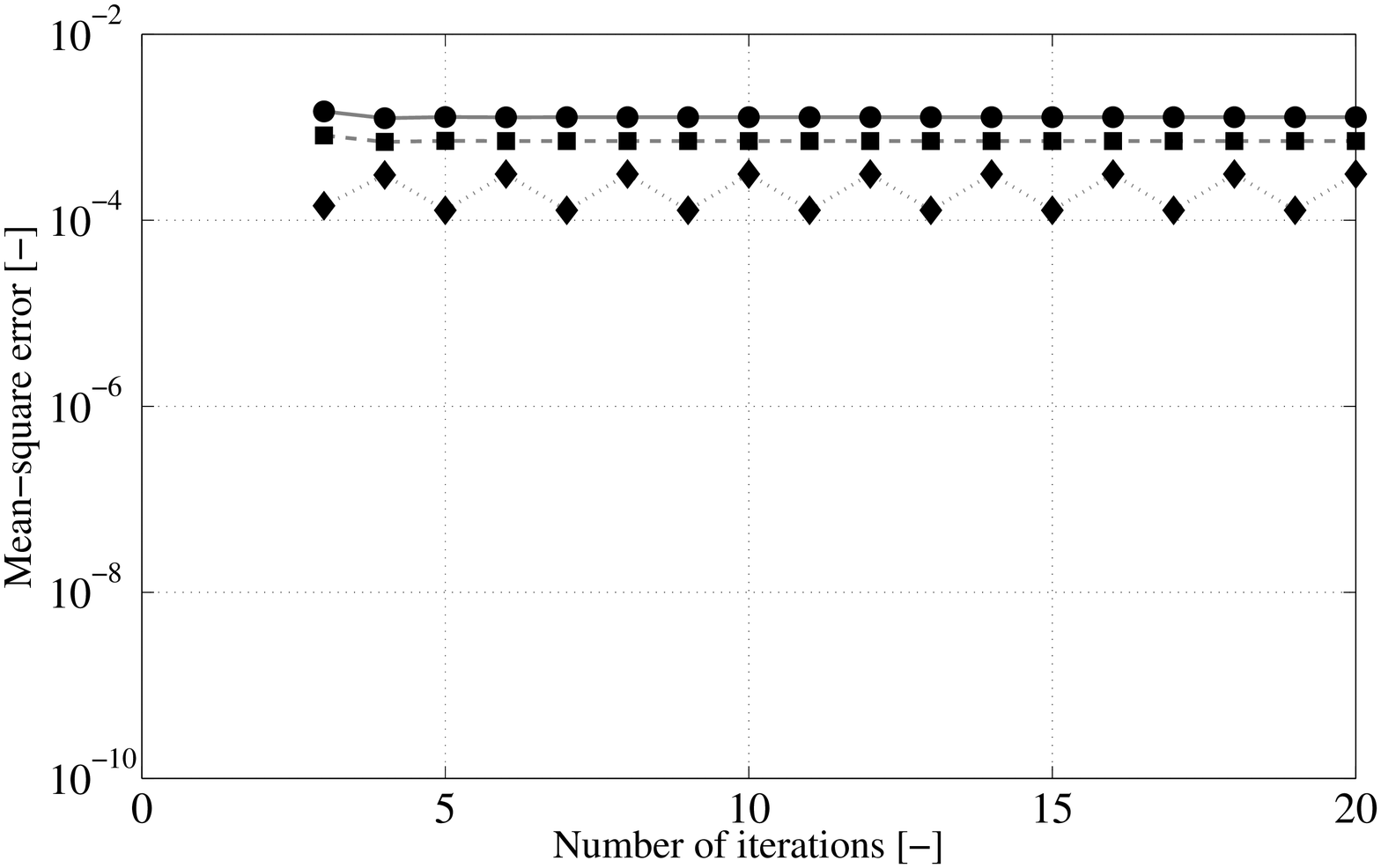}}
    \subfigure[\hspace{-1mm}Neutron-flux-based distance for $k\hspace{-1mm}=\hspace{-1mm}100\,\text{$[\text{J/K/cm/s}]$}$.]{\includegraphics[width=0.49\textwidth]{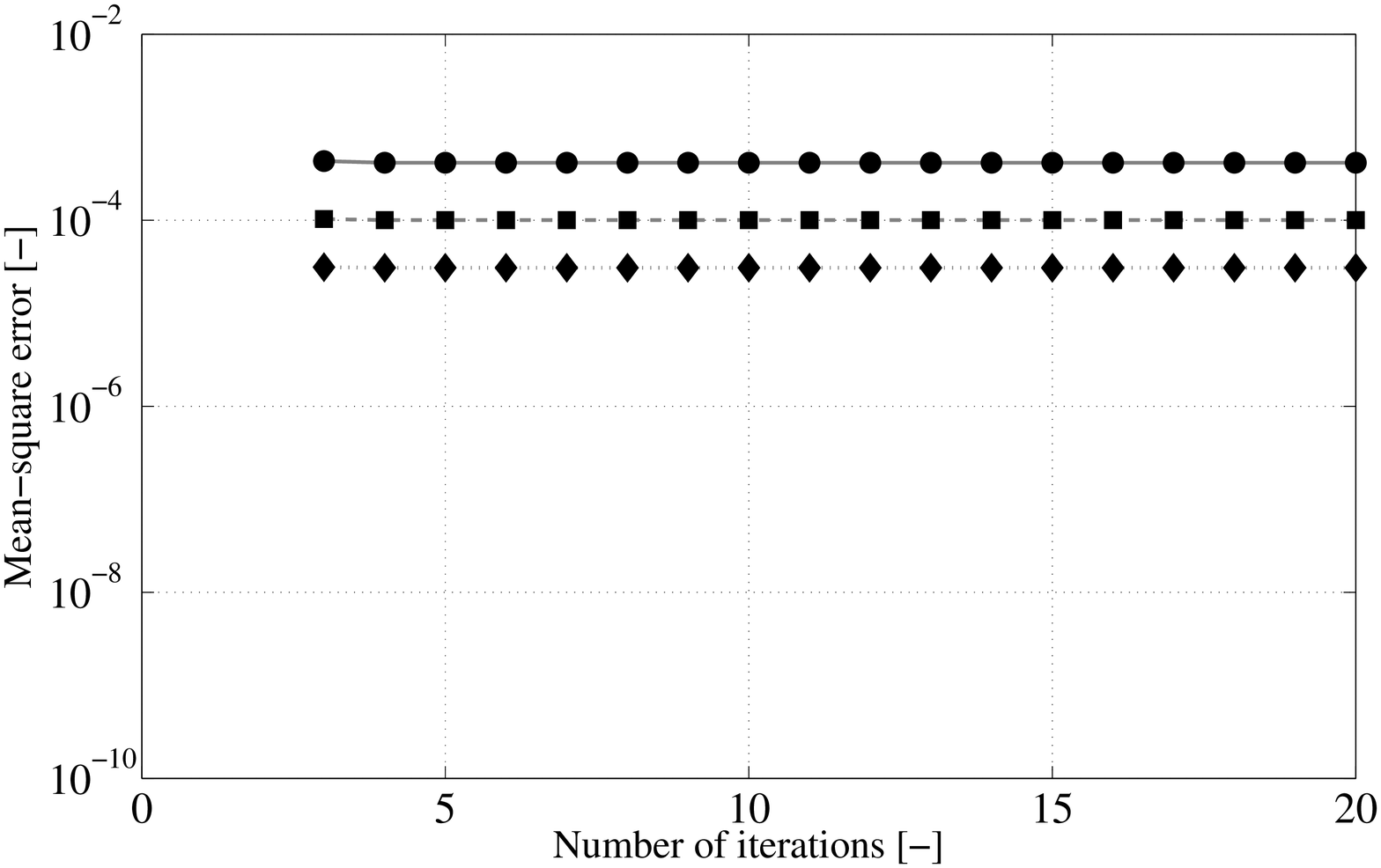}}
    \hfill
    \subfigure[Neutron-flux-based distance for $k=1\,\text{$[\text{J/K/cm/s}]$}$.]{\includegraphics[width=0.49\textwidth]{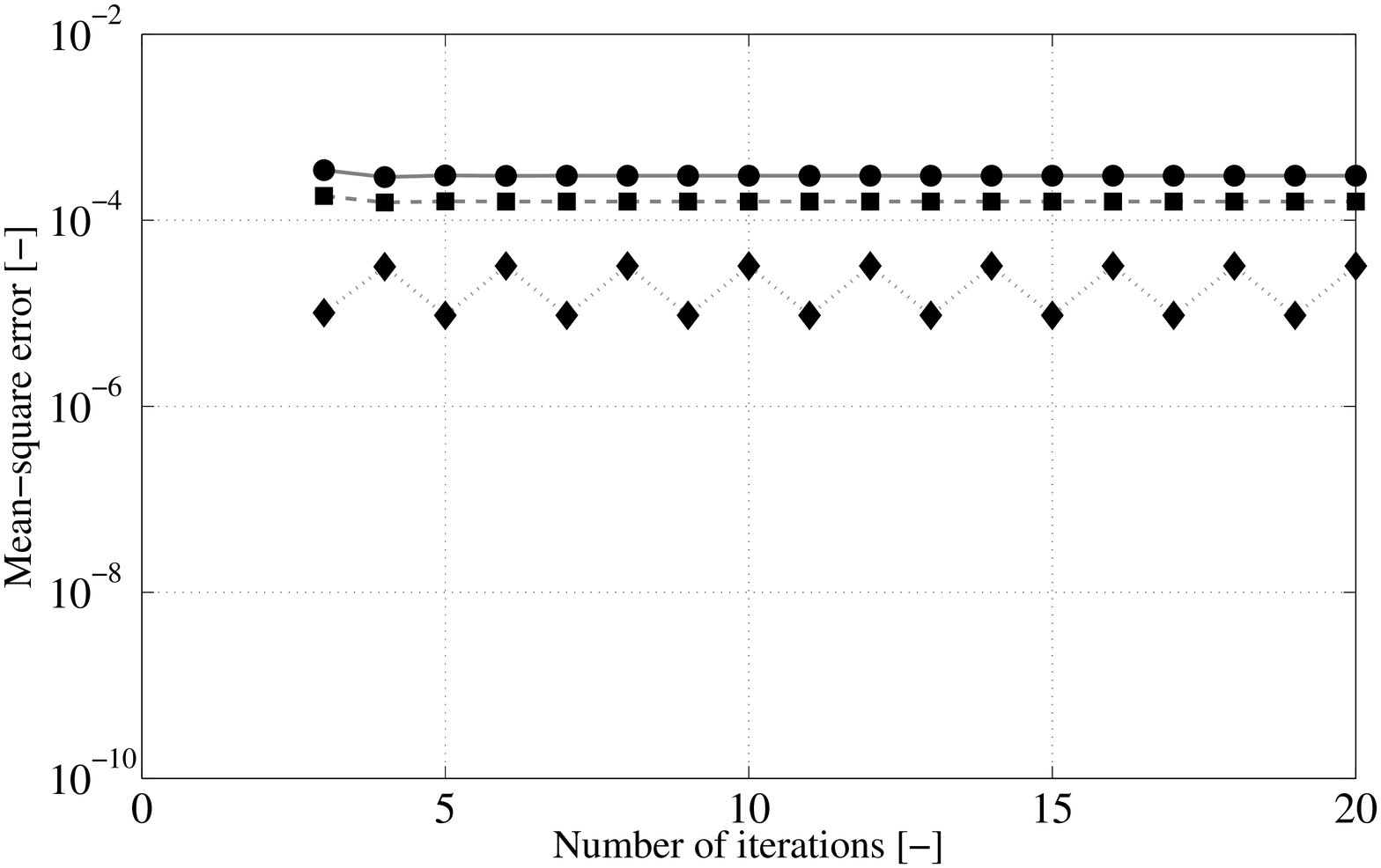}}
    \centerline{$\small{\ell\mapsto\sqrt{\sum_{|\boldsymbol{\alpha}|=0}^{p}\|\boldsymbol{T}_{\boldsymbol{\alpha}}^{\ell}-\widehat{\boldsymbol{T}}{}_{\boldsymbol{\alpha}}^{\ell}\|_{\boldsymbol{W}}^{2}}\Big/\sqrt{\sum_{|\boldsymbol{\alpha}|=0}^{p}\|\boldsymbol{T}_{\boldsymbol{\alpha}}^{\ell}\|_{\boldsymbol{W}}^{2}}}\normalsize$.}
    \centerline{$\small{\ell\mapsto\sqrt{\sum_{|\boldsymbol{\alpha}|=0}^{p}\|\boldsymbol{\Phi}_{\boldsymbol{\alpha}}^{\ell}-\widehat{\boldsymbol{\Phi}}{}_{\boldsymbol{\alpha}}^{\ell}\|_{\boldsymbol{W}}^{2}}\Big/\sqrt{\sum_{|\boldsymbol{\alpha}|=0}^{p}\|\boldsymbol{\Phi}_{\boldsymbol{\alpha}}^{\ell}\|_{\boldsymbol{W}}^{2}}}\normalsize$.}
    \caption{Convergence analysis: mean-square distance between the successive approximations generated by the simulations involving and not involving dimension reduction for~$tol=0.90\times \sigma_{\boldsymbol{T}}^{2}$ (circles), $tol=0.95\times \sigma_{\boldsymbol{T}}^{2}$ (squares), and~$tol=0.99\times \sigma_{\boldsymbol{T}}^{2}$ (diamonds).}\label{fig:figure12b}
  \end{center}
\end{figure}

We repeated the PC-based simulation involving dimension reduction for the error tolerance levels~$tol=0.90\times \sigma_{\boldsymbol{T}}^{2}$, $tol=0.95\times \sigma_{\boldsymbol{T}}^{2}$, and~$tol=0.99\times \sigma_{\boldsymbol{T}}^{2}$;
each of these levels corresponded to an increased accuracy that the KL decomposition was required to maintain at each iteration. 
Figure~\ref{fig:figure12} indicates that the KL decomposition systematically retained more terms when the error tolerance level was set to a lower value, and, thus, higher accuracy was required.
Moreover, the comparison of Figs.~\ref{fig:figure12}(a) and~\ref{fig:figure12}(b) reveals that the reduced-dimensional representations of the random temperature retained less terms for~$k=100\,\text{$[\text{J/K/cm/s}]$}$ than for~$k=1\,\text{$[\text{J/K/cm/s}]$}$, which is consistent with our earlier observation that the eigenvalues of the KL decomposition obtained for $k=100\,\text{$[\text{J/K/cm/s}]$}$ decay at a faster rate than those obtained for $k=1\,\text{$[\text{J/K/cm/s}]$}$.

Figure~\ref{fig:figure12b} shows that the distance between the successive approximations obtained by the iterative method that involved dimension reduction and the one that did not involve dimension reduction remained bounded as the iterations progressed, and that this distance can be reduced by improving the accuracy of the KL decomposition by decreasing the error tolerance level.

\subsection{Concluding remarks}
In this section, we demonstrated the effectiveness of the proposed dimension-reduction methodology through a stochastic multiphysics model of the transport of neutrons in a reactor with temperature feedback. 
A study of the parameters of the model indicated that the random temperature exchanged between the neutronics and the heat subproblem was of a low effective stochastic dimension when it was smoothed by the effect of significant heat diffusion terms.
Numerical results showed that the KL decomposition could then extract a low-dimensional representation of the random temperature as it passed between the subproblems while accuracy was maintained.
The convergence study showed that the accuracy of the solution could be increased systematically by retaining more terms in the KL decomposition.

\section{Conclusion}
While most coupled models can be expected to be affected by a large number of sources of uncertainty, information exchanged between subproblems and successive iterations often resides in a considerably lower dimensional space than the sources themselves. 
This argument calls for the use of dimension-reduction techniques for the representation of exchanged information. 
In this work, we described the formulation and implementation of an adaptation of the KL decomposition to represent information as it passes from subproblem to subproblem and from iteration to iteration during the numerical simulation of a stochastic coupled model.

Implementations can capitalize on the reduced-dimensional interface created between the subproblems and iterations to enable a more computationally efficient solution of the subproblems in a reduced-dimensional space.
This can be achieved by performing key algorithmic operations in terms of the reduced stochastic degrees of freedom of the exchanged uncertainty representations.
This is the main computational benefit of the proposed dimension-reduction methodology, and we plan to explore solution algorithms of this form in future work.

\appendix

The approximation of a random variable by a truncated KL decomposition as it is exchanged most often results in a truncation error, which will most likely have an effect on the solution of the subproblems and will also propagate to the subsequent approximations generated by the iterative method.
This appendix concisely analyzes the effect of this truncation error.

\subsection{Objective}
The probability theory offers several ways in which two random variables can be considered to be equal or a sequence of random variables can be considered to converge, namely the almost sure, mean-square, and distributional senses, among other ways. 
An error analysis can be conducted usefully from any of these perspectives.
We adopt a mean-square setting because it perfectly agrees with the convergence properties of the KL decomposition.
Specifically, in a mean-square sense, we examine the distance between the sequences~$\{(\boldsymbol{u}^{\ell},\boldsymbol{v}^{\ell)})\}_{\ell=1}^{\infty}$ and~$\{(\hat{\boldsymbol{u}}{}^{\ell},\hat{\boldsymbol{v}}{}^{\ell})\}_{\ell=1}^{\infty}$ determined by the iterative method that does not involve dimension direction and by the one that does involve dimension reduction, respectively.

\subsection{Contraction-type assumption}
We refer to the mapping that transforms any random variable~$(\boldsymbol{u},\boldsymbol{v},\boldsymbol{y},\boldsymbol{x})$ defined on $(\Theta,\mathcal{T},P)$ with values in~$\real^{r}\times\real^{s}\times\real^{r_{0}}\times\real^{s_{0}}$ into a random variable~$(\boldsymbol{a}(\boldsymbol{u},\boldsymbol{x},\boldsymbol{\xi}),\boldsymbol{b}(\boldsymbol{y},\boldsymbol{v},\boldsymbol{\zeta}),\boldsymbol{h}(\boldsymbol{u},\boldsymbol{\xi}),\boldsymbol{k}(\boldsymbol{v},\boldsymbol{\zeta}))$ defined on $(\Theta,\mathcal{T},P)$ with values in $\real^{r}\times\real^{s}\times\real^{r_{0}}\times\real^{s_{0}}$ as the \textit{iteration mapping} associated with the iterative method~(\ref{eq:couplingAGS12}) that does not involve dimension reduction.
Clearly, an error analysis requires assumptions concerning the sensitivity of this mapping with respect to its arguments.  
In this analysis, we assume the iteration mapping to be a contraction mapping.
Specifically, we assume that the iteration mapping is a contraction mapping in the block-maximum norm, in that there exists a Lipschitz continuity modulus~$\alpha$ in $[0,1)$ such that

\vspace{-3mm}
\small{
\begin{align}
\notag&\max\hspace{-0.5mm}\big(\hspace{-0.5mm}\vectornorm{\boldsymbol{a}(\boldsymbol{u},\boldsymbol{x},\boldsymbol{\xi})\hspace{-0.5mm}-\hspace{-0.5mm}\boldsymbol{a}(\boldsymbol{u}',\boldsymbol{x}',\boldsymbol{\xi})}_{r}\hspace{-0.75mm},\hspace{-0.5mm}\vectornorm{\boldsymbol{b}(\boldsymbol{y},\boldsymbol{v},\boldsymbol{\zeta})\hspace{-0.5mm}-\hspace{-0.5mm}\boldsymbol{b}(\boldsymbol{y}',\boldsymbol{v}',\boldsymbol{\zeta})}_{s}\hspace{-0.75mm},\hspace{-0.5mm}\vectornorm{\boldsymbol{h}(\boldsymbol{u},\boldsymbol{\xi})\hspace{-0.5mm}-\hspace{-0.5mm}\boldsymbol{h}(\boldsymbol{u}',\boldsymbol{\xi})}_{r_{0}}\hspace{-1.25mm},\hspace{-0.5mm}\vectornorm{\boldsymbol{k}(\boldsymbol{v},\boldsymbol{\zeta})\hspace{-0.5mm}-\hspace{-0.5mm}\boldsymbol{k}(\boldsymbol{v}',\boldsymbol{\zeta})}_{s_{0}}\hspace{-0.5mm}\big)\\
\label{eq:contraction}&\quad\quad\quad\quad\leq\alpha\max\big(\vectornorm{\boldsymbol{u}-\boldsymbol{u}'}_{r},\vectornorm{\boldsymbol{v}-\boldsymbol{v}'}_{s},\vectornorm{\boldsymbol{y}-\boldsymbol{y}'}_{r_{0}},\vectornorm{\boldsymbol{x}-\boldsymbol{x}'}_{s_{0}}\big)
\end{align}}\normalsize
for all the second-order random variables~$(\boldsymbol{u},\boldsymbol{v},\boldsymbol{y},\boldsymbol{x})$ and~$(\boldsymbol{u}',\boldsymbol{v}',\boldsymbol{y}',\boldsymbol{x}')$ defined on $(\Theta,\mathcal{T},P)$ with values in~$\real^{r}\times\real^{s}\times\real^{r_{0}}\times\real^{s_{0}}$; here, $\vectornorm{\cdot}_{r}$, $\vectornorm{\cdot}_{s}$, $\vectornorm{\cdot}_{r_{0}}$, and~$\vectornorm{\cdot}_{s_{0}}$ are suitably weighted mean-square norms on the spaces of second-order random variables on $(\Theta,\mathcal{T},P)$ valued in~$\real^{r}$, $\real^{s}$, $\real^{r_{0}}$, and $\real^{s_{0}}$.
This assumption is justified in the current context by the fact that this assumption is consistent with conditions that ensure the mean-square convergence of the sequence generated by the iterative method~(\ref{eq:couplingAGS12}) that does not involve  dimension reduction: by the Banach contraction mapping theorem, a sufficient condition for~$\{(\boldsymbol{u}^{\ell},\boldsymbol{v}^{\ell)}\}_{\ell=1}^{\infty}$ to converge in mean square is that the iteration mapping should map a non-empty closed set of random variables into itself and that~(\ref{eq:contraction}) should be fulfilled.

\subsection{Adaptive selection of the reduced dimension and induced error bound}
If the iterative method~(\ref{eq:couplingAGS12}) not involving dimension reduction and the iterative method~(\ref{eq:couplingAGSred12}) involving dimension reduction are initialized identically, if the iteration mapping is a contraction mapping in the block-maximum norm such that~(\ref{eq:contraction}) is fulfilled with continuity modulus~$\alpha$ in $[0,1)$, and if the reduced dimensions, denoted here as~$d_{\ell}$ and $e_{\ell}$, are selected at each iteration so as to achieve systematically a prescribed accuracy $tol\geq 0$ such that
\begin{equation}
\label{eq:nlbound}\begin{aligned}
&\max\big(\vectornorm{\hat{\boldsymbol{u}}{}^{\ell-1}-\hat{\boldsymbol{u}}{}^{\ell-1,e_{\ell}}}_{r}, \vectornorm{\hat{\boldsymbol{x}}{}^{\ell-1}-\hat{\boldsymbol{x}}{}^{\ell-1,e_{\ell}}}_{s_{0}}\big)\leq tol,&&\forall\ell\in\integer,\\
&\max\big(\vectornorm{\hat{\boldsymbol{y}}{}^{\ell}-\hat{\boldsymbol{y}}{}^{\ell,d_{\ell}}}_{r_{0}}, \vectornorm{\hat{\boldsymbol{v}}{}^{\ell-1}-\hat{\boldsymbol{v}}{}^{\ell-1,d_{\ell}}}_{s}\big)\leq tol,&&\forall\ell\in\integer,
\end{aligned}
\end{equation}
then the distance between the sequences of successive approximations generated by the iterative method that does not involve dimension reduction and the one that does involve dimension reduction satisfies
\begin{equation}
\max\big(\vectornorm{\boldsymbol{u}^{\ell}-\hat{\boldsymbol{u}}{}^{\ell}}_{r}, \vectornorm{\boldsymbol{v}^{\ell}-\hat{\boldsymbol{v}}{}^{\ell}}_{s}\big)\leq\frac{2\alpha}{1-\alpha}tol,\quad\forall\ell\in\integer.\label{eq:errorbound}
\end{equation}
It should be emphasized that the reduced dimensions~$d_{\ell}$ and~$e_{\ell}$ can always be chosen such that~(\ref{eq:nlbound}) is fulfilled owing to the key properties~(\ref{eq:KLaccuracy}) and~(\ref{eq:keypropertyyy}) of the KL decomposition. 

\subsection{Proof of the error bound}
From the definition of the iteration processes~(\ref{eq:couplingAGS12}) and~(\ref{eq:couplingAGSred12}), from the contraction-type assumption~(\ref{eq:contraction}), and from the triangle inequality, it is observed that
\begin{equation}
\max\big(\big\|\boldsymbol{u}^{\ell}-\hat{\boldsymbol{u}}{}^{\ell}\big\|_{r},\big\|\boldsymbol{v}^{\ell}-\hat{\boldsymbol{v}}{}^{\ell}\big\|_{s}\big)\leq\alpha\max\big(\big\|\boldsymbol{u}^{\ell-1}-\hat{\boldsymbol{u}}{}^{\ell-1}\big\|_{r},\big\|\boldsymbol{v}^{\ell-1}-\hat{\boldsymbol{v}}{}^{\ell-1}\big\|_{s}\big)+2\,\alpha\,tol.
\end{equation}
If~$(\boldsymbol{u}^{0},\boldsymbol{v}^{0})=(\hat{\boldsymbol{u}}{}^{0},\hat{\boldsymbol{v}}{}^{0})$, the repeated application of this inequality yields
\begin{equation}
\max\big(\big\|\boldsymbol{u}^{\ell}-\hat{\boldsymbol{u}}{}^{\ell}\big\|_{r},\big\|\boldsymbol{v}^{\ell}-\hat{\boldsymbol{v}}{}^{\ell}\big\|_{s}\big)\leq\sum_{k=1}^{\ell}2\alpha^{\ell-k+1}\,tol,\label{eq:keyinequality}
\end{equation}
which subsequently yields~(\ref{eq:errorbound}) because by the convergence of geometric series for $0\leq\alpha<1$, we obtain~$\lim_{\ell\rightarrow+\infty}\sum_{k=1}^{\ell}\alpha^{\ell-k+1}=\alpha/(1-\alpha)$.
This reasoning indicates that a contraction-type assumption is sufficient to obtain a result of the form given by~(\ref{eq:errorbound}) because it ensures that the effect of any approximation error introduced owing to the truncation of a KL decomposition is systematically diminished, and not amplified, as the error propagates to subsequent iterations.

\subsection{Concluding remarks}
This appendix showed that under a contraction-type assumption, the sequence of successive approximations determined by the iterative method involving dimension reduction can be brought as close as desired to the sequence determined by the iterative method not involving dimension reduction by suitably adjusting the reduced dimension at each iteration.

Although this appendix provides a useful result, it is not a comprehensive convergence and error analysis of the proposed dimension-reduction methodology.
Further research can be conducted to analyze convergence properties and error estimates under weakened assumptions such as those involving nonuniform contraction-type or local differentiability properties, for instance, by following the approaches given in~\citep{brezzi1980,fink1983,tsuchiya1997}.
Moreover, a posteriori error estimators can be developed, for instance, by following the approach given in~\citep{estep2008}.  

\section*{Acknowledgments}
This work was supported by DOE through an ASCR grant.
The authors would also like to thank the two anonymous reviewers for their suggestions and Professor Christian Soize for relevant discussions during the final stages of the preparation of this paper.   

\bibliography{multiphysics2}

\end{document}